\documentclass[11pt]{article}
\usepackage{amsfonts}
\usepackage{mathrsfs}
\usepackage{latexsym,amssymb,graphics}
\usepackage{amsmath,amssymb,epsfig,subfigure,url,color}
\usepackage{graphicx}
\usepackage{subfigure}
\usepackage{exscale}
\usepackage{relsize}
\usepackage{amsmath}
\usepackage{txfonts}
\usepackage{url}
\usepackage[flushleft]{caption2}
\usepackage{multirow} \usepackage{makecell}
\usepackage{multicol}
\usepackage[numbers,sort&compress]{natbib}
\date{}
\begin{document}

\title{Total variation with overlapping group sparsity for image deblurring under impulse noise \thanks{The work of Gang Liu, Ting-Zhu Huang and Jun Liu is by 973 Program (2013CB329404), NSFC (61370147), Sichuan Province Sci. \& Tech. Research Project (2012GZX0080). The work of Xiao-Guang Lv is supported by Postdoctoral Research Funds (2013M540454, 1301064B).}}
\author{Gang Liu\thanks{School of Mathematical Sciences, University of Electronic Science and Technology of China, Chengdu, Sichuan, 611731, P. R. China ({\em wd5577@163.com}).},
Ting-Zhu Huang\thanks{School of Mathematical Sciences, University of Electronic Science and Technology of China, Chengdu, Sichuan, 611731, P. R. China ({\em tingzhuhuang@126.com}).},
Jun Liu\thanks{School of Mathematical Sciences, University of Electronic Science and Technology of China, Chengdu, Sichuan, 611731, P. R. China ({\em junliucd@163.com}).},
Xiao-Guang Lv \thanks{School of Mathematical Sciences, Nanjing Normal University, Nanjing, Jiangsu, 210097, P.R. China({\em xiaoguanglv@126.com}).}}

\maketitle

\begin{abstract}
The total variation (TV) regularization method is an effective method for image deblurring in preserving edges.
However, the TV based solutions usually have some staircase effects. In this paper, in order to alleviate the staircase effect, we propose a new model for restoring blurred images with impulse noise. The model consists of an $\ell_1$-fidelity term and a TV with overlapping group sparsity (OGS) regularization term.
Moreover, we impose a box constraint to the proposed model for getting more accurate solutions. An efficient and effective algorithm is proposed to solve the model under the framework of the alternating direction method of multipliers (ADMM). We use an inner loop which is nested inside the majorization minimization (MM) iteration for the subproblem of the proposed method.
Compared with other methods, numerical results illustrate that the proposed method, can significantly improve the restoration quality, both in avoiding staircase effects and in terms of peak signal-to-noise ratio (PSNR) and relative error (ReE).

\end{abstract}
\emph{keywords}: Impulse noise; total variation; overlapping group sparsity; convex  optimization; image deblurring; ADMM

AMS: 94A08; 68U10; 65F22

\markboth{\textsc{TV with OGS for image deblurring under impulse noise}}{G. LIU, T.-Z. HUANG, J. LIU \and X.-G. LV}

\section{Introduction}

Image deblurring and denoising has been widely studied in last decades. In the literatures, it is widely assumed that the observed image is the convolution of a standard linear and space invariant blurring function with the true image plus some noise. Let $g$ denote the blurred and noisy
image, $h$ the blur kernel, $f$ the original image and $n$ the
noise. The image $f$ is assumed to be a real function defined
on a bounded and piecewise smooth open subset $\Gamma$ of
$\mathbb{R}^2$. In general, the image formation process can be modeled as: $g=h\star f+n$, where
``$\star$'' denotes the two-dimensional convolution operation. Image deblurring is to estimate
the true image $f$ from the blurred and noisy image $g$.
As is well known, image deblurring is a typically ill-posed problem \cite{KV1990,TA1977}. To handle this problem,
regularization technique is usually considered to obtain a stable and accurate solution. That is, we want to solve the following problem
\begin{equation}\label{REG}
  \min_f \psi(f) + \frac{\mu}{2}\int _\Gamma |h\star f-g |^2 {\rm d}x, \
\end{equation}
where the first term is called the regularization term, the second term is called the fidelity term ($\ell_{2}$-fidelity), $\mu>0$ is the
regularization parameter, and $\psi$ is the regularization functional.

Without loss of generality, we assume that a discretized image has $n\times n$ pixels, then $f$, $g$ and $n$ are vectors of length $n^2$. Let $H$ be the corresponding blurring matrix of $n^2\times n^2$ from $h$ \cite{HN2006}. Then the discretized form of the minimization problem (\ref{REG}) is equivalent to the following matrix-vector form
\begin{equation}\label{DREG}
  \min_f \psi(f) + \frac{\mu}{2} \|Hf-g\|_2^2 , \ \
\end{equation}
where $\|\cdot\|_2$ denotes the Euclidean $\ell_2$ norm. Notice that $H$ is a matrix of block circulant with circulant blocks (BCCB) structure when periodic boundary conditions are applied or other structures when other boundary conditions are applied \cite{HN2006}.

How to choose a good regularization functional is an active area of
research in the imaging science. In the early 1960s, D. L. Phillips \cite{DP1962} and A. N.
Tikhonov \cite{AT1963} proposed the definition of $\psi$ as an $\ell_2$-type norm (academically called Tikhonov regularization), that is,
$\psi = \|Lf\|_2^2$ with $L$ an identity operator or difference operator.
Although the functional $\psi$ of this type has the advantage
of facilitating the calculations, it is rarely used in current
practice because it has the drawback of penalizing discontinuities
in resulting solutions, for instance, over-smoothing edges. Therefore,
this is not a good choice since natural images have many edges.

To overcome this drawback, many different types of regularization functionals have been proposed, for instance, the regularizers introduced in \cite{LLS2011} for image denoising. Particularly, one well-known model was introduced by Rudin, Osher and Fatemi (ROF) in \cite{ROF1992}.
They proposed a total variation (TV) regularization with an $\ell_{2}$-fidelity term ($\ell_2$-TV)
for image restoration. Its corresponding minimization task is:
\begin{equation}\label{DROF}
  \min_f \|f \|_\text{TV} + \frac{\mu}{2} \|Hf-g\|_2^2 , \ \
{\rm over}\ f \in \text{BV}(\Gamma),
\end{equation}
where BV$(\Gamma)$ denotes the space of functions of bounded variation. That is, $f = \text{BV}(\Gamma)$ if and only if $f\in L^1 (\Gamma)$ and BV-seminorm (TV norm)
$$ \int_\Gamma | Df |=\sup\left\{\int_\Omega f \text{ div} \vec v: \vec v \in(C_0 ^\infty (\Gamma))^2 ,
| \vec v |\leqslant 1\right\} <\infty ,$$
where ``div'' is the divergence operator.
Here $\|f \|_\text{TV}$ in (\ref{DROF}) is the discretization form of the TV norm $\int_\Gamma | Df |$ \cite{Chamb2004,GLNG2009}. It is defined by
$ \|f \|_\text{TV}:= \sum\limits_{1\leqslant i,j \leqslant n} \|(\nabla f)_{i,j}\|_2 = \sum\limits_{1\leqslant i,j \leqslant n}$ $ \sqrt{{|(\nabla_x f)}_{i,j}|^2 +
 |(\nabla_y f)_{i,j}|^2 }$ which is called isotropic TV, or $ \|f \|_{TV}:= \sum\limits_{1\leqslant i,j \leqslant n}$    $\|(\nabla f)_{i,j}\|_1$ $= \sum\limits_{1\leqslant i,j \leqslant n} |(\nabla_x f)_{i,j}| + |(\nabla_y f)_{i,j}|$ which is named anisotropic TV, where $\|\cdot\|_1$ denotes the Euclidean $\ell_1$ norm.
Operator $\nabla: \mathbb{R}^{n^2}\rightarrow\mathbb{R}^{2\times{n^2}}$ denotes the discrete gradient
 operator (under periodic boundary conditions) which is defined by
 $(\nabla f)_{i,j}=((\nabla_x f)_{i,j},(\nabla_y f)_{i,j} ),$
with
\begin{equation*}
(\nabla_x f)_{i,j} =\left\{\begin{array}{lll}  f_{i+1,j}-f_{i,j}&{\rm if}&i<n,\\
                                            f_{1,j}-f_{n,j}&{\rm if}&i=n,\end{array}\right.
(\nabla_y f)_{i,j} =\left\{\begin{array}{lll}  f_{i,j+1}-f_{i,j}&{\rm if}&j<n,\\
                                            f_{i,1}-f_{i,n}&{\rm if}&j=n,\end{array}\right.
\end{equation*}
for $i,j=1,2,\cdots,n$, where $f_{i,j}$ refers to the $((j-1)n+i)$th entry
of the vector $f$ (it is the $(i,j)$th pixel location of the
$n\times n$ image, and this notation remains valid throughout the paper
unless otherwise specified).

Many methods have been proposed to solve the restoration model (\ref{DROF}) such as the fast TV deconvolution (FTVd)
\cite{WY2008,YYZW2009}, the augmented Lagrangian method (ALM) \cite{Esser2009,EB1992,TBS2012,CX2010}, the dual methods \cite{CGM1999,CX2010}, and
the split Bregman method
\cite{COS2009,Esser2009}. We know that these methods are designed for Gaussian noise removal. However, in many
cases, the noise does not satisfy the Gaussian
assumption, for instance, the noise may follow a Laplace distribution \cite{AR1994}. There has been a growing interest in using an $\ell_1$-fidelity term instead of the $\ell_2$-fidelity term for image restoration in many literatures \cite{GLNG2009,MILA2004,YGO2007,YWYW2009} for considering another non-Gaussian noise--impulse noise. The corresponding regularization model with an {$\ell_1$}-fidelity term (REGL1) that our work will consider is as:
 \begin{equation}\label{REGL1}
  \min \psi(f) + \mu \|Hf-g\|_1.
\end{equation}

 A well-known approach is to use the TV regularizer by $\psi (f) = \|f \|_{TV}$, which we call $\ell_1$-TV. Recently, Wang {\it et al.} \cite{YWYW2009} used the FTVd method to solve the $\ell_1$-TV model fast. Guo {\it et al.} \cite{GLNG2009} proposed a fast $\ell_1$-TV algorithm for image restoration in the $\ell_1$-TV model. Their method was to add a penalty term by using the variable substitution method, which belongs to penalty methods in optimization. They employed an alternating minimization method to solve it. They first proved the convergence of their method and second got better results and faster than FTVd. Wu {\it et al.} \cite{WZT2011} used ALM to solve the $\ell_1$-TV model. They also got better results than FTVd.
More recently, Chan {\it et al.} \cite{CTY2013} proposed a constrained total variation (TV) regularization method for image restoration for $\ell_1$-TV (\ref{DROF}) (they also considered $\ell_2$-TV but it had nothing to do with our work). They used alternating direction method of multipliers (ADMM)
to solve the models by combining augmented Lagrangian method and variable splitting method.
Their method used a box constrained projection to ensure the restored images stay in a given dynamic range. They got better results by their method than other methods such as FTVd and ALM. Their numerical results showed that for some images where there are many pixels
with values lying on the boundary of the dynamic range, the gain could get very high numerical superiority in the peak signal-to-noise ratio. It shows that the constrained projection is necessary in image restoration.

However, although the TV regularization using in the restoration problems can recover sharp edges of a degraded image, it
also gives rise to some undesired effects and transforms smooth
signal into piecewise constants, the so-called staircase effects \cite{Chamb1997,CMM2000}.
To overcome this deficiency, one effective method is to replace the original TV norm by a high-order TV norm. The high-order
TV regularization schemes have been studied so far mainly for overcoming the staircase effects while preserving the edges in the restored image. More details please refer to \cite{CMM2000,LLT2003,LT2006,Steidl2006}.
But the high-order TV usually has some other behaviors. For example, it may transforms the smooth signal to over-smoothing, and it may take more time to compute.
More recently, Selesnick and Chen \cite{SC2013} proposed an overlapping group sparsity (OGS) TV regularizer to one-dimensional signal denoising. They applied the majorization minimization (MM) method to solve their model. Their numerical experiments showed that their method can overcome staircase effects effectively. However, their method has the disadvantages of the low speed of computation and the difficulty to extend to the two-dimensional case because they did not choose a variable substitution method.

In this paper, inspired by Selesnick and Chen's work \cite{SC2013}, we consider to set $\psi$ in (\ref{REGL1}) to
be the OGS-TV functional to two-dimension images deblurring under impulse noise. Moreover, we impose a box constraints to the proposed model to obtain more accurate solutions. We propose an efficient and effective algorithm to solve the model under the framework of the alternating direction method of multipliers (ADMM). We use an inner loop which is nested inside the majorization minimization (MM) iteration for the subproblem of the proposed method by a variable substitution method.
 The numerical experiments show that our method using the OGS-TV regularizer could avoid staircase effects effectively. Moreover, the numerical results also show that our method is very effective and competitive with other methods, such as Chan's method \cite{CTY2013} and Guo's method \cite{GLNG2009}.

The outline of the rest of this paper is as follows. In Section 2
we will briefly introduce the definition of the OGS regularization functional. We will also review the MM method and ADMM, which are used in our proposed method. In Section 3,
we propose an OGS-TV based model for recoving images under blur and impulse noise and derive an efficient solving algorithm.
The numerical results are given in Section 4. Finally, we conclude this paper in Section 5. 

\section{Some preliminaries} 
\subsection{OGS-TV}

In \cite{SC2013}, the authors denoted a $K$-point group ($K$ denotes the group size) of the vector
$t\in \mathbb{R}^n$ by
\begin{equation}\label{sampK}
  t_{i,K}=[t(i),t(i+1),\cdots,t(i+K-1)] \in \mathbb{R}^K.
\end{equation}
Note that $t_{i,K}$ can be seen as a block of $K$ contiguous samplings
of $t$ staring at index $i$. With the notation (\ref{sampK}), a group sparsity
regularizer for one-dimensional case is defined in \cite{SC2013} as
\begin{equation}\label{OverK}
  \zeta(t) = \sum_{i=1}^n \|t_{i,K}\|_2.
\end{equation}
Similarly, we can define a $K$-square-point group of a two-dimensional signal such as images considered in this work $v\in \mathbb{R}^{n^2}$ (vector $v$ is obtained by rearranging the $n\times n$ entires in a matrix in column-major order, that is, the $(i,j)$th entry of a matrix is assigned to be the $((j-1)n+i)$th entry of the vector $v$)
by
\begin{equation}\label{sampKK}
\tilde{v}_{i,j,K,K} =\left[\begin{array}{cccc}  v_{i-K_l,j-K_l}  &v_{i-K_l,j-K_l+1}  &\cdots &v_{i-K_l,j+K_r}\\
                                         v_{i-K_l+1,j-K_l}  &v_{i-K_l+1,j-K_l+1}  &\cdots &v_{i-K_l+1,j+K_r}\\
                                          \vdots    &\vdots       &\ddots &\vdots\\
                                         v_{i+K_r,j-K_l}&v_{i+K_r,j-K_l+1}&\cdots &v_{i+K_r,j+K_r}\\ \end{array}\right]\in \mathbb{R}^{K\times K}
\end{equation}
where $K_l=[\frac{K-1}{2}]$, $K_r=[\frac{K}{2}]$ and [$x$] denotes the largest integer less than or equal to $x$. The group size is denoted by $K^2$. Note that $\tilde{v}_{i,j,K,K}$ can be seen as a square block of $K\times K$ contiguous samplings of $v$ with the center at index $(i,j)$. Here we choose a group entries around the objective point rather than a group following the objective point like one-dimensional in \cite{SC2013} because of the faster and easier computation in the experiments. Moreover, to the best of our knowledge, the former is much better than the later in image restoration because the pixels in image are related to or depended on all the ambient pixels rather than partial surrounding pixels.
Let $v_{i,j,K,K}$ be a $K^2$-vector obtained by arranging the $K\times K$ elements of $\tilde v_{i,k,K,K}$ in lexicographic order.
This notation also remains valid throughout the paper unless otherwise specified. Then the overlapping group sparsity functional of the two-dimensional array can be defined by
\begin{equation}\label{OverKK}
  \varphi(v) = \sum_{i=1}^n\sum_{j=1}^n \|v_{i,j,K,K}\|_2.
\end{equation}
Here we can also use $\ell_1$ norm (as anisotropic TV) instead of $\ell_2$ norm, and we choose $\ell_2$ norm here because of the better performance of $\ell_2$ norm in classic TV regularization models. From the definition above, we can easily get that this function is convex. Consequently, we define the regularization functional $\psi$ in (\ref{REGL1})
to be the form
\begin{equation}\label{OGSREG}
  \psi(f) = \varphi(\nabla_x f) +\varphi(\nabla_y f) .
\end{equation}
We call the regularizer $\psi$ in (\ref{OGSREG}) as the OGS anisotropic TV functional, and call the corresponding convex minimization model (\ref{REGL1}) L1-OGS-ATV.

\subsection{The MM method}
The MM method is an asymptotical method in solving optimization problems. That is, instead of directly solving
a difficult minimization problem $P(v)$, the MM method approach solves a sequence of easier optimization problems
$Q(v,v_{k})$ $(k=0,1,2,...)$ firstly and then manages to get the minimizer of $P(v)$.
Generally, an MM iterative
algorithm for minimizing $P(v)$ has the form
\begin{equation}\label{MM}
  v^{k+1} = \arg\min_v Q(v,v^k),
\end{equation}
where $Q(v,v')\geqslant P(v)$ for all $v$, $v'$, and $Q(v^k,v^k)= P(v^k)$, i.e., each functional $Q(v, v¡ä)$ is a majorizor of $P(v)$.
When $P(v)$ is convex, then under the former conditions, the
sequence $v^k$ produced by (\ref{MM}) converges to the minimizer of
$P(v)$ \cite{FBN2007,OBF2009}.

Before the discussion of our method, we consider a minimization problem of the form
\begin{equation}\label{MMden}
  \min_v P(v) = \left\{\frac{\alpha}{2}\|v - v_0\|_2^2 + \varphi(v) \right\}, v\in \mathbb{R}^{n^2},
\end{equation}
where $\alpha$ is a positive parameter and the functional $\varphi$ is by
the definition in (\ref{OverKK}). In \cite{CS2014}, the authors have studied this problem elaborately. However, for the sake of completeness, we briefly introduce the solving method here and fix minor bugs of \cite{CS2014} for using Matlab built-in function \texttt{conv2}. To derive an
effective and efficient algorithm with the MM scheme for solving the problem (\ref{MMden}),
we want to find a majorizor of $P(v)$. Here, we only need
to find a majorizor of $\varphi(v)$ because of the simple enough
quadratic term of the first term in (\ref{MMden}). Note that
\begin{equation}\label{Major}
  \frac{1}{2}(\frac{1}{\|u\|_2} \|v\|_2^2 + \|u\|_2) \geqslant \|v\|_2,
\end{equation}
for all $v$ and $u\neq 0$ ($u,v\in \mathbb{R}^{n^2}$) with equality when $u = v$. Substituting
each group of $\varphi(v)$ into (\ref{Major}) and summing them, we get a majorizor of
$\varphi(v)$
\begin{equation}\label{MMajor1}
  S(v,u) = \frac{1}{2}\sum_{i=1}^n\sum_{j=1}^n\left[\frac{1}{\|{u}_{i,j,K,K}\|_2}
  \|{v}_{i,j,K,K}\|_2^2+\|{u}_{i,j,K,K}\|_2\right],
\end{equation}
with
\begin{equation}\label{MMajor2}
  S(v,u)\geqslant \varphi(v), S(u,u)= \varphi(u),
\end{equation}
provided $\|{v}_{i,j,K,K}\|\neq 0$ for all $i,j$. After simple calculation,
$S(v, u)$ can be rewritten as
\begin{equation}\label{MMajor3}
  S(v,u)=\frac{1}{2} \|\Lambda(u){v}\|_2^2 +C,
\end{equation}
where $C$ is independent of $v$, and $\Lambda(u)$ is a
diagonal matrix with each diagonal component
\begin{equation}\label{MMajor4}
 [\Lambda(u)]_{m,m}=\sqrt {\sum_{i=-K_l}^{K_r}\sum_{j=-K_l}^{K_r}\left[\sum_{k_1=-K_l}^{K_r}
 \sum_{k_2=-K_l}^{K_r}|u_{m-i+k_1,m-j+k_2}|^2\right]^{-\frac{1}{2}}},
\end{equation}
with $m = 1, 2, \cdots, n^2$. The entries of $\Lambda$ can be easily computed
by using Matlab built-in function \texttt{conv2}. Then a majorizor
of $P(v)$ can be easily given by
\begin{equation}\label{MMajor5}
 Q(u,v)=\frac{\alpha}{2}\|v-v_0\|_2^2 + S(v,u)=\frac{\alpha}{2}\|v-v_0\|_2^2 + \frac{1}{2} \|\Lambda(u){v}\|_2^2 +C,
\end{equation}
with $Q(v,u)\geqslant P(v)$ for all $u$, $v$, and $Q(u,u)= P(u)$.
 To minimize $P(v)$,
the MM aims to iteratively solve
\begin{equation}\label{MM1}
  v^{k+1} = \arg\min_v \frac{\alpha}{2}\|v-v_0\|_2^2 + \frac{1}{2} \|\Lambda(v^k){v}\|_2^2, k= 1, 2, \cdots,
\end{equation}
with the solution
\begin{equation}\label{MM2}
  \hat{v}^{k+1} = \left(I + \frac{1}{\alpha}\Lambda^2(v^k)\right)^{-1}{v}_0, k= 1, 2, \cdots.
\end{equation}
where $ I $ is  an  identity  matrix  with  the  same  size  of  $\Lambda(v^k)$.  We can easily get that $\Lambda^2(v^k)$ is also a diagonal matrix with each diagonal component $[\Lambda^2(v^k)]_{m,m}$ that equals to the form of removing the out root of the right term of (\ref{MMajor4}). Moreover, the inversion of the matrix $I + \frac{1}{\alpha}\Lambda^2(v^k)$ can  be  computed  very  efficiently since it only requires simple component-wise calculation.  Therefore, we obtain the  Algorithm  1  (we  call  it  MMOdn  for convenience) for solving the problem (\ref{MMden}).\\\\
\begin{tabular}{l}
\hline
\hline
{\textsc{\textbf{Algorithm 1}}} \textup{MMOdn for solving (\ref{MMden})}\\
\hline
$1$. \textit{\textbf{initialization}}:  Starting point $v= v_0$, $\alpha$, group size $K^2$,  $K_l=[\frac{K-1}{2}]$, \\
\quad $K_r=[\frac{K}{2}]$, $\epsilon$, Maximum inner iterations $NIt$,  $k=0$.\\
$2$. \textit{\textbf{iteration}}:\\
\quad  Do
\end{tabular}\\
\begin{tabular}{l}
$\begin{array}{rcl}
\quad \quad [\Lambda^2(v^k)]_{m,m}&=& {\sum\limits_{i=-K_l}^{K_r}\sum\limits_{j=-K_l}^{K_r}\left[\sum\limits_{k_1=-K_l}^{K_r}
 \sum\limits_{k_2=-K_l}^{K_r}|v^k_{m-i+k_1,m-j+k_2}|^2\right]^{-\frac{1}{2}}},\ \ \ \\
  {v}^{k+1} &= &\left(I + \frac{1}{\alpha}\Lambda^2(v^k)\right)^{-1}{v}_0,\\
  k&=&k+1,
\end{array}$\\
\quad \textup{until}\ $\|v^{k+1}-v^k\|_2/\|v^k\|_2<\epsilon$ \textup{or} $k>NIt$.\\
$3$. \textit{\textbf{get $v^k$}}.\\
\hline
\end{tabular}\\

\subsection{Variable splitting and ADMM}

Consider an unconstrained optimization problem in which the objective function is the sum of two functions,
which is written as
\begin{equation}\label{CSCO}
   \min  \phi_1(x_1) + \phi_2(x_2),
    {\rm s.\ t.\ }  A_1 x_1 + A_2 x_2 =b,
      x_i \in \chi_i, i = 1,2,
 \end{equation}
where $\phi_i: \mathbb{R}^{n_i}\rightarrow \mathbb{R}$ are closed proper convex functions, $\chi_i\subseteq \mathbb{R}^{n_i}$ are closed convex sets, $A_i\in \mathbb{R}^{l\times n_i}$, and $b\in\mathbb{R}^l$ is a given vector.
The augmented Lagrangian function \cite{NW2006} of (\ref{CSCO})
\begin{equation}\label{CSCOAL}
  \begin{array}{cl}
   \mathcal{L}(x_1,x_2,\lambda) &=\phi_1(x_1) + \phi_2(x_2) - \lambda^T (A_1 x_1 + A_2 x_2-b)
    +\frac{\beta}{2}\|A_1 x_1+ A_2 x_2-b\|_2^2\\
    &=\phi_1(x_1) + \phi_2(x_2) + \frac{\beta}{2}\|A_1 x_1+ A_2 x_2-b -\frac{\lambda}{\beta}\|_2^2 +C
       \end{array}
\end{equation}
where $\lambda\in\mathbb{R}$ is the Lagrange multiplier, $\beta$ is a penalty
parameter which  controls  the  linear  constraint, and $C$ does not depend on $x_1,x_2$.
 The idea of the ADMM is to find a saddle point of $\mathcal{L}$.
 Usually, the ADMM consists in minimizing $\mathcal{L}$ in
an alternating way, for instance, minimizing $\mathcal{L}$ with  respect  to  $x_1$ by fixing $x_2$ and $\lambda$.
That delivers to the
following simple but powerful algorithm ADMM:\\\\
\begin{tabular}{l}
\hline
\hline
{\textsc{\textbf{Algorithm 2}}} \textup{Classic ADMM for the minimization problem (\ref{CSCO})}\\
\hline
 \textit{\textbf{initialization}}:  Starting point$x_1^0$, $x_2^0$, $\lambda^0$, $\beta$.\\
 \textit{\textbf{iteration}}:\\
\end{tabular}\\
\begin{tabular}{l}
$\begin{array}{rcl}
 x_1^{k+1}&=&\arg\min \phi_1(x_1) + \frac{\beta}{2}\|A_1 x_1+ A_2 x_2^k-b -\frac{\lambda^k}{\beta}\|_2^2,\\
 x_2^{k+1}&=&\arg\min \phi_2(x_2) + \frac{\beta}{2}\|A_1 x_1^{k}+ A_2 x_2-b -\frac{\lambda^k}{\beta}\|_2^2,\\
 \lambda^{k+1}&=&\lambda^{k} - \beta(A_1 x_1^{k+1}+ A_2 x_2^{k+1}-b),\quad\quad\quad\quad\quad\quad\quad\quad\quad\quad\ \\
  k&=&k+1,\\
\end{array}$\\
\textit{\textbf{until a stopping criterion is satisfied.}}\\
\hline
\end{tabular}\\

According to the literature \cite{EB1992}, we can see the classic ADMM is convergent because of the nonexpansive and absolute summable properties of the $x_1$ and $x_2$ subproblems.
However, the speed is not fast. In order to speed the convergence, we can introduce a step length parameter $\gamma$ for updating the multiplier \cite{RG1984,GM1976,HY1998}.  The algorithm framework is outlined as follows for the general ADMM.\\\\
\begin{tabular}{l}
\hline
\hline
{\textsc{\textbf{Algorithm 3}}} \textup{General ADMM for the minimization problem (\ref{CSCO})}\\
\hline
 \textit{\textbf{initialization}}:  Starting point$x_1^0$, $x_2^0$, $\lambda^0$, $\beta$.\\
 \textit{\textbf{iteration}}:\\
\end{tabular}\\
\begin{tabular}{l}
$\begin{array}{rcl}
 x_1^{k+1}&=&\arg\min \phi_1(x_1) + \frac{\beta}{2}\|A_1 x_1+ A_2 x_2^k-b -\frac{\lambda^k}{\beta}\|_2^2,\\
 x_2^{k+1}&=&\arg\min \phi_2(x_2) + \frac{\beta}{2}\|A_1 x_1^{k}+ A_2 x_2-b -\frac{\lambda^k}{\beta}\|_2^2,\\
 \lambda^{k+1}&=&\lambda^{k} - \gamma\beta(A_1 x_1^{k+1}+ A_2 x_2^{k+1}-b),\quad\quad\quad\quad\quad\quad\quad\quad\quad\quad\ \\
  k&=&k+1,\\
\end{array}$\\
\textit{\textbf{until a stopping criterion is satisfied.}}\\
\hline
\end{tabular}\\

Here, $\gamma>0$ is also called a relax parameter. In fact, if $\gamma=1$, the general ADMM is the classic ADMM. We do not fix $\gamma = 1$ in our work since $\gamma$ plays an important role in convergence of the general ADMM. From the literatures \cite{RG1984,GM1976,HY1998}, the general ADMM is convergent if $\gamma\in (0,(\sqrt(5)+1)/2)$. Moreover, $\gamma= 1.618$ makes it converge noticeably faster than $\gamma= 1$. Therefore, we set $\gamma= 1.618$ in our work.

\section{Proposed method}

With the definition of (\ref{OGSREG}), we will consider
a minimization problem of the form (L1-OGS-ATV)
\begin{equation}\label{L1OGSTV}
  \min_f \varphi(\nabla_x f) +\varphi(\nabla_y f) + \mu \|Hf-g\|_1.
\end{equation}
Note that for any
true digital image, its pixel value can attain only a finite number
of values. Hence, it is natural to require all pixel values of
the restored image to lie in a certain interval $[a, b]$,
see \cite{CTY2013} for more details. For example, for 8-bit images, we would like
to restore them in a dynamic range $[0, 255]$. More in general, with the easy computation and the certified results in \cite{CTY2013}, we only consider all the images located on the range $[0, 1]$. Therefore, the images we mentioned all lie in the interval $[0, 1]$.
We define a projection operator $\mathcal{P}_{\Omega}$ on the set $\Omega =\left\{f\in \mathbb{R}^{n\times n}|0\leqslant f\leqslant 1\right\}$,
\begin{equation}\label{Cpro}
  \mathcal{P}_{\Omega}(f)_{i,j}=
  \left\{\begin{array}{ll}
  0,&f_{i,j}<0,\\
  f_{i,j},&f_{i,j}\in[0,1],\\
  1,&f_{i,j}>1.
  \end{array}\right.
\end{equation}
Similarly as \cite{CTY2013}, we will solve the problem
\begin{equation}\label{L1COGSTV}
  \min_{f\in\Omega} \varphi(\nabla_x f) +\varphi(\nabla_y f) + \mu \|Hf-g\|_1.
\end{equation}
We refer to this model as CL1-OGS-ATV. Obviously, this model is also convex.

By introducing new auxiliary variables $v_x$, $v_y$, $z$, $w$, we transform
the minimization problem (\ref{L1COGSTV}) to the equivalent constrained minimization problem
\begin{equation}\label{L1COGSTVVS}
  \min_{w\in\Omega,f,z,v_x,v_y}\left\{ \varphi(v_x) +\varphi(v_y) + \mu \|z\|_1: \
  {\rm s.\ t.\ }  z=Hf-g, v_x =\nabla_x f, v_y=\nabla_y f, w=f \right\}.
\end{equation}
Note that the constraint is now imposed on $w$ instead of $f$. The augmented Lagrangian
function of (\ref{L1COGSTVVS}) is
\begin{equation}\label{L1COGSTVAL}
  \begin{array}{rl}
   \mathcal{L}(v_x,v_y,z,w,f;\lambda_1,\lambda_2,\lambda_3,\lambda_4)=&\varphi(v_x) - \lambda_1^T (v_x-\nabla_x f)
    +\frac{\beta_1}{2}\|v_x-\nabla_x f\|_2^2\\
    &+\varphi(v_y) - \lambda_2^T (v_y-\nabla_y f)
    +\frac{\beta_1}{2}\|v_y-\nabla_y f\|_2^2\\
    &+\mu \|z\|_1 - \lambda_3^T \left(z-(Hf - g)\right) + \frac{\beta_2}{2}\|z-(Hf-g)\|_2^2\\
    &-\lambda_4^T (w -f) + \frac{\beta_3}{2}\|w-f\|_2^2,
       \end{array}
\end{equation}
where $\beta_1, \beta_2, \beta_3 >0$ are penalty parameters and $\lambda_1, \lambda_2, \lambda_3, \lambda_4 \in\mathbb{R}^{n^2}$ are the Lagrange multipliers.
According to the scheme of the general ADMM mentioned above (Section 2.3), for a given $(v_x^k,v_y^k,z^k,w^k,f^k$; $\lambda_1^k,\lambda_2^k,\lambda_3^k,\lambda_4^k)$, the next iteration
$(v_x^{k+1},v_y^{k+1},z^{k+1},w^{k+1}$, $f^{k+1}$; $\lambda_1^{k+1},\lambda_2^{k+1},\lambda_3^{k+1},\lambda_4^{k+1})$ is generated as follows:

1. Fix $f = f^k$, $\lambda_1=\lambda_1^k, \lambda_2=\lambda_2^k, \lambda_3=\lambda_3^k, \lambda_4=\lambda_4^k$, $z=z^k$, $w=w^k$,
and minimize (\ref{L1COGSTVAL}) with respect to $v_x$ and $v_y$.
The minimizers are obtained by
\begin{equation}\label{COTVSUB1}\begin{array}{rl}
  v_x^{k+1}&=\arg\min \varphi(v_x) - {\lambda_1^k}^T (v_x-\nabla_x f^k)
    +\frac{\beta_1}{2}\|v_x-\nabla_x f^k\|_2^2\\
    &=\arg\min \varphi(v_x) +\frac{\beta_1}{2}\|v_x-\nabla_x f^k - \frac{\lambda_1^k}{\beta_1}\|_2^2,\\
  \end{array}
\end{equation}
\begin{equation}\label{COTVSUB11}\begin{array}{rl}
  v_y^{k+1}&=\arg\min \varphi(v_y) - {\lambda_2^k}^T (v_y-\nabla_y f^k)
    +\frac{\beta_1}{2}\|v_y-\nabla_y f^k\|_2^2\\
    &=\arg\min \varphi(v_y) +\frac{\beta_1}{2}\|v_y-\nabla_y f^k - \frac{\lambda_2^k}{\beta_1}\|_2^2.\\
  \end{array}
\end{equation}
It is obvious that problems (\ref{COTVSUB1}) and (\ref{COTVSUB11}) match
the framework of the problem (\ref{MMden}), thus the solutions of (\ref{COTVSUB1}) and (\ref{COTVSUB11})
can be obtained by using Algorithm 1 (Section 2.2), respectively.

2. Compute $z^{k+1}$ easily.
\begin{equation*}\begin{array}{rl}
  z^{k+1}&=\arg\min\mu \|z\|_1 - {\lambda_3^k}^T \left(z-(Hf^k - g)\right) + \frac{\beta_2}{2}\|z-(Hf^k-g)\|_2^2\\
    &=\arg\min\mu \|z\|_1 + \frac{\beta_2}{2}\|z-(Hf^k-g) - \frac{\lambda_3^k}{\beta_2}\|_2^2.\\
  \end{array}
\end{equation*}
The minimization with respect to $z$ can be given by the well-known Shrinkage \cite{YWYW2009} explicitly by:
\begin{equation}\label{COTVSUB2}
  z^{k+1}={\rm sgn}\left\{Hf^k - g +\frac{\lambda_3^k}{\beta_2}\right\}\circ \max\left\{|Hf^k - g +\frac{\lambda_3^k}{\beta_2}|-\frac{\mu}{\beta_2},0\right\},
\end{equation}
where $|\cdot|$, sgn and ``$\circ$'' represent the componentwise absolute value, signum function, and componentwise product, respectively. 

3. Compute $w^{k+1}$ easily.
\begin{equation*}\begin{array}{rl}
  w^{k+1}&=\arg\min -{\lambda_4^k}^T (w -f^k) + \frac{\beta_3}{2}\|w-f^k\|_2^2\\
    &=\arg\min \frac{\beta_3}{2}\|w-f^k - \frac{\lambda_4^k}{\beta_3}\|_2^2.\\
  \end{array}
\end{equation*}
The minimizer is given explicitly by
\begin{equation}\label{COTVSUB3}
  w^{k+1}=\mathcal{P}_{\Omega}\left[f^k + \frac{\lambda_4^k}{\beta_3}\right].
\end{equation}

4.  Compute $f^{k+1}$ by solving the normal equation
\begin{equation}\label{COTVSUB4}
 \begin{array}{l}
 (\beta_1 (\nabla_x^* \nabla_x+\nabla_y^* \nabla_y)+\beta_2 H^* H + \beta_3 I)f^{k+1} \\
=\nabla_x^*(\beta_1 v_x^{k+1}-\lambda^k_1)+{\nabla_y}^*(\beta_1 v_y^{k+1}-\lambda_2^k)+H^* (\beta_2 z^{k+1} -\lambda_3^k) +\beta_2 H^* g + \beta_3 (w^{k+1} - \frac{\lambda_4^k}{\beta_3}),\\
\end{array}
\end{equation}
where ``$*$'' denotes the conjugate transpose, see \cite{CX2010} for more details.
Since all the parameters are positive, the coefficient matrix in
(\ref{COTVSUB4}) are always invertible and symmetric positive definite. In addition, note that
$H$, $\nabla_x$, $\nabla_y$ have BCCB
structure under periodic boundary conditions. We
know that the computations with BCCB matrix can be very
efficient by using fast Fourier transforms (FFTs).

5.  Update the multipliers via
\begin{equation}\label{COTVSUB5}
  \left\{
  \begin{array}{*{20}{l}}
  \lambda_1^{k+1}&=&\lambda_1 ^k - \gamma\beta_1(v_x^{k+1}-\nabla_x f^{k+1}),\\
  \lambda_2^{k+1}&=&\lambda_2 ^k - \gamma\beta_1(v_x^{k+1}-\nabla_x f^{k+1}),\\
  \lambda_3^{k+1}&=&\lambda_3 ^k - \gamma\beta_2(z^{k+1}-(Hf^{k+1}-g)),\\
  \lambda_4^{k+1}&=&\lambda_4 ^k - \gamma\beta_3(w^{k+1} -f^{k+1}).\\
  \end{array} \right.
\end{equation}

Based on the discussions above, we present the ADMM algorithm using inner MM iteration for solving the convex CL1-OGS-ATV model (\ref{L1COGSTV})
shown as Algorithm 4.\\\\
\begin{tabular}{l}
\hline
\hline
{\textsc{\textbf{Algorithm 4}}} \textup{CL1-OGS-ATV-ADM4 for the minimization problem (\ref{L1COGSTV})}\\
\hline
 \textit{\textbf{initialization}}:  \\ \ Starting point $v_x^0=v_y^0=g$, $k=0$, $\beta_1$, $\beta_2$, $\beta_3$, $\gamma$, $\mu$,
 group size $K\times K$, \\\ $\lambda_i^0=0$, $i=1,2,3,4$,
 Maximum inner iterations $NIt$.\\
 \textit{\textbf{iteration}}:\\
\end{tabular}\\
\begin{tabular}{l}
$\begin{array}{l}
 1.\ Compute\ v_x^{k+1}\ and\ v_y^{k+1}\  according\ to\ (\ref{COTVSUB1})\ and\ (\ref{COTVSUB11}).\\
 2.\ Compute\ z^{k+1}\ according\ to\ (\ref{COTVSUB2}).\\
3.\ Compute\ w^{k+1}\ according\ to\ (\ref{COTVSUB3}).\\\end{array}$\\
 \end{tabular}\\
\begin{tabular}{l}
$\begin{array}{l}
4.\ Compute\ f^{k+1}\ by\ solving\ (\ref{COTVSUB4}).\\
5.\ update\ \lambda_i^0=0, i=1,2,3,4\ according\ to\ (\ref{COTVSUB5}).\quad\quad\quad\quad\quad\quad\quad\ \quad\quad\\
6.\ k=k+1.\\
\end{array}$\\
\textit{\textbf{until a stopping criterion is satisfied.}}\\
\hline
\end{tabular}\\

Since  CL1-OGS-ATV-ADM4 a special case of the general ADMM for the case with two blocks of variables $(v_y, v_x, w, z)$ and $f$, if the Step (1) of Algorithm 4 can be solved exactly, the convergence for CL1-OGS-ATV-ADM4 can be
guaranteed \cite{EB1992}. In this case, if the relax parameter $\gamma \in (0, \frac{\sqrt{5}+1}{2})$, Algorithm
4 is convergent, more details please refer to \cite{RG1984,GM1976,HY1998}.
Besides, although step (1) of Algorithm 4 can not be solved exactly, our numerical experiments will verify the convergence of Algorithm 4. 

\section{Numerical results}

In this section, we present several numerical results to illustrate the performance of the proposed method. We compare our method CL1-OGS-ATV-ADM4 (``Ours'' for short) with other methods, Chan's ADM2CTVL1 proposed in \cite{CTY2013} (``CTY'' for short, Algorithm 2 in \cite{CTY2013} for the constrained TV-L1 model) and Guo's fast $\ell_1$-TV proposed in \cite{GLNG2009} (``GLN'' for short).

All experiments are carried out on Windows 7 32-bit and
Matlab 2010a running on a desktop equipped with an Intel
Core i3-2130 CPU with 3.4 GHz and 3.4 GB of RAM.

The  quality  of  the  restoration  results  is  measured  quantitatively
by  using  the  peak signal-to-noise ratio (PSNR) in decibel (dB) and the  relative  error  (ReE):
\begin{equation*}
\textrm{PSNR} = 10\log_{10} \frac{n^2 {\rm Max}_I^2}{\|f-\bar{f}\|_2^2},\quad
\textrm{ReE} = \frac{\|f -\bar{f}\|_2}{\|\bar{f}\|_2},
\end{equation*}
where $\bar{f}$ and $f$ denote the original and restored images
respectively, and ${\rm Max}_I$ represents the maximum possible pixel value of the image.
In our experiments,  ${\rm Max}_I= 1$.
The stopping criterion used in our work is set to be
\begin{equation}
 \frac{|\mathcal{F}^{k+1}-\mathcal{F}^k|}{|\mathcal{F}^k|}< 10^{-5},
\end{equation}
where $\mathcal{F}^k$ is the objective function value of the respective model in the $k$th iteration, which is
\begin{equation}
 \mathcal{F}^k= \varphi(\nabla_x f^k) +\varphi(\nabla_y f^k) + \mu \|Hf^k-g\|_1.
\end{equation}\\
The stopping criterions of CTY (same as ours) and GLN are set to default as their literature mentioned.

\begin{figure}
  \centering
  \subfigure{\label{originalcame}
    \includegraphics[width=0.2\textwidth,clip]{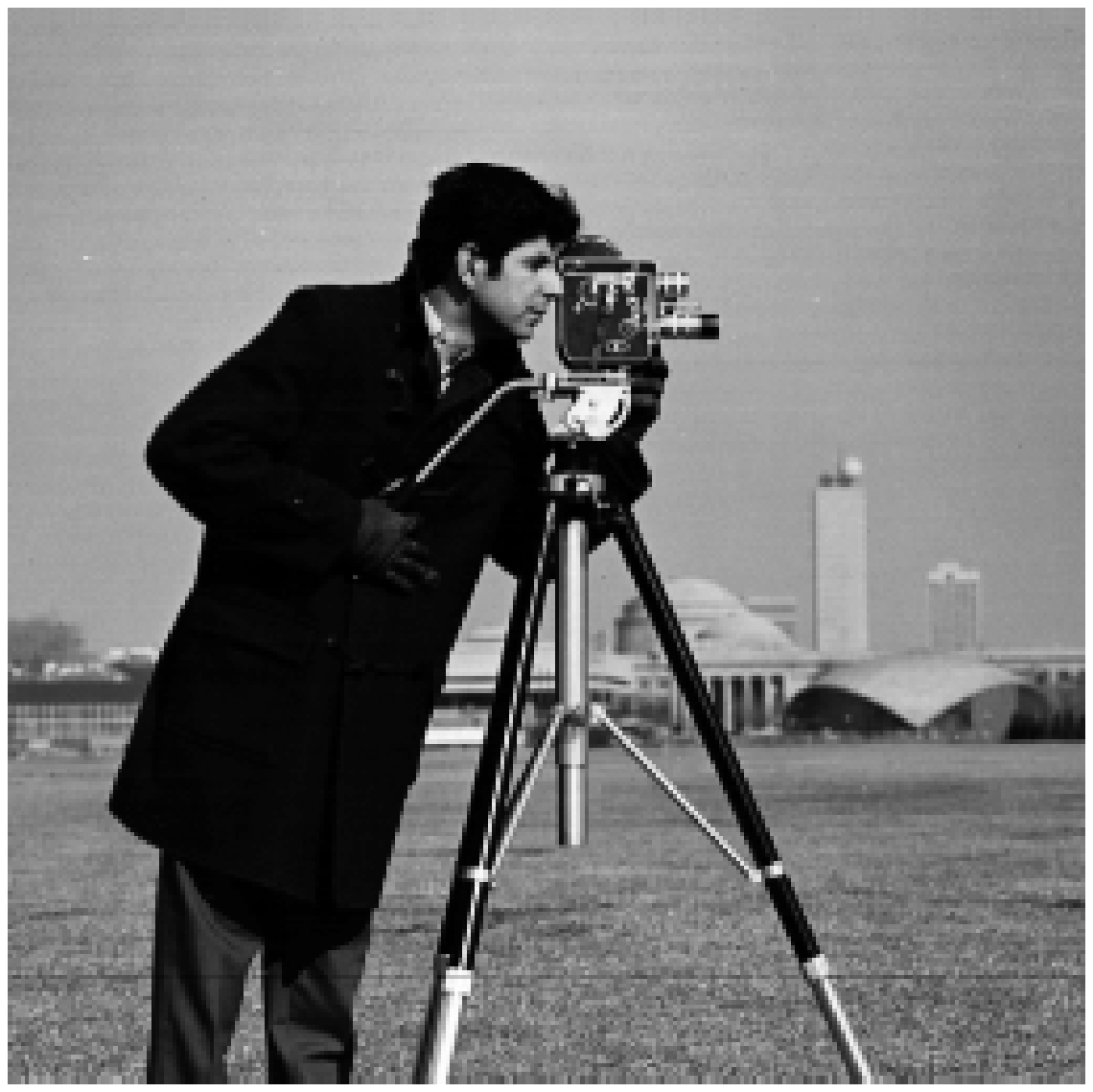}}
    \hspace{0.001in}
   \subfigure{\label{originalsate}
    \includegraphics[width=0.2\textwidth,clip]{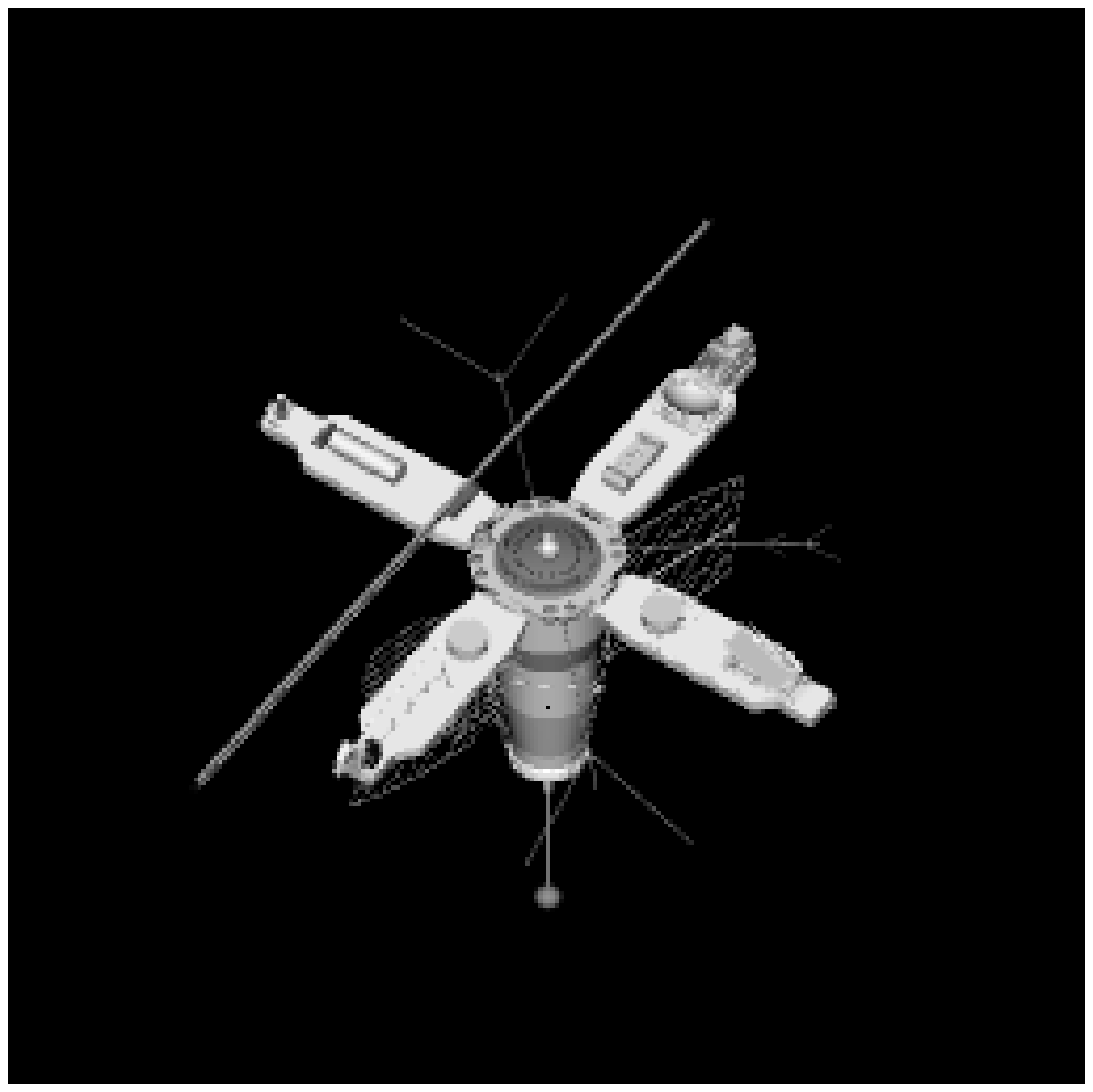}}
    \hspace{0.001in}
    \subfigure{\label{originalhouse}
    \includegraphics[width=0.2\textwidth,clip]{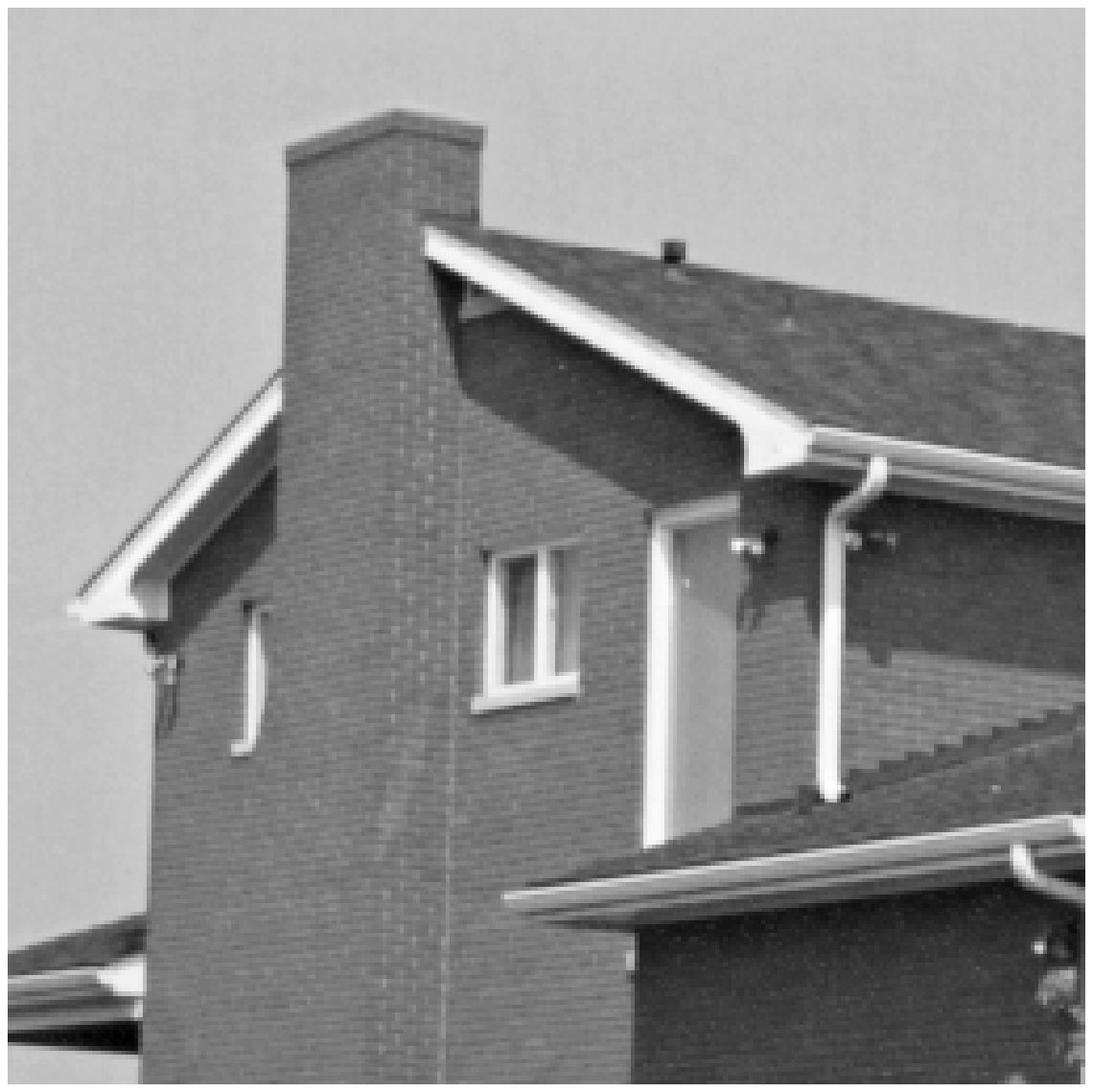}}
    \hspace{0.001in}
    \subfigure{\label{originalboat}
    \includegraphics[width=0.2\textwidth,clip]{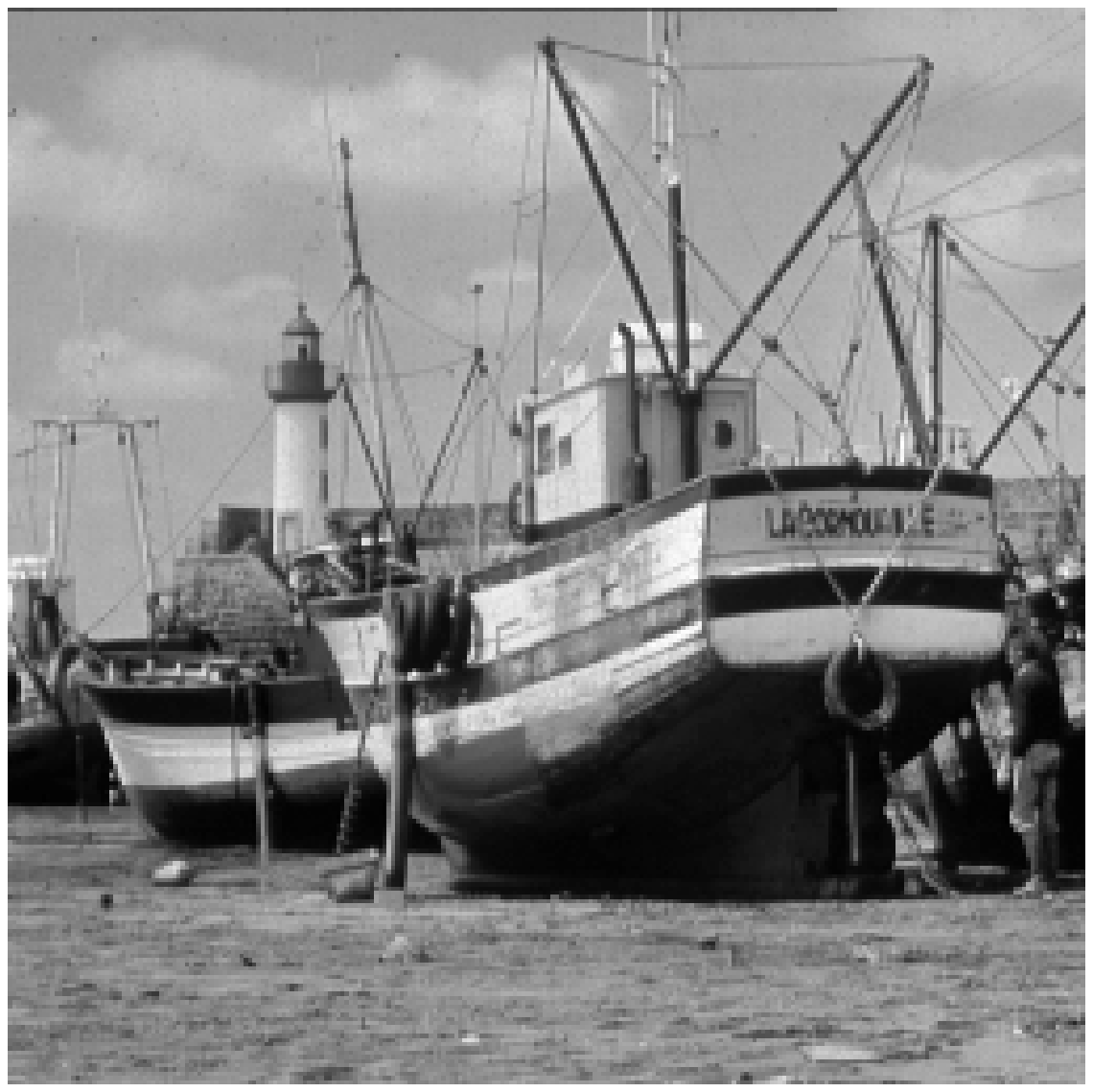}}
    \hspace{0.001in}
   \subfigure{\label{originalbarbar}
    \includegraphics[width=0.2\textwidth,clip]{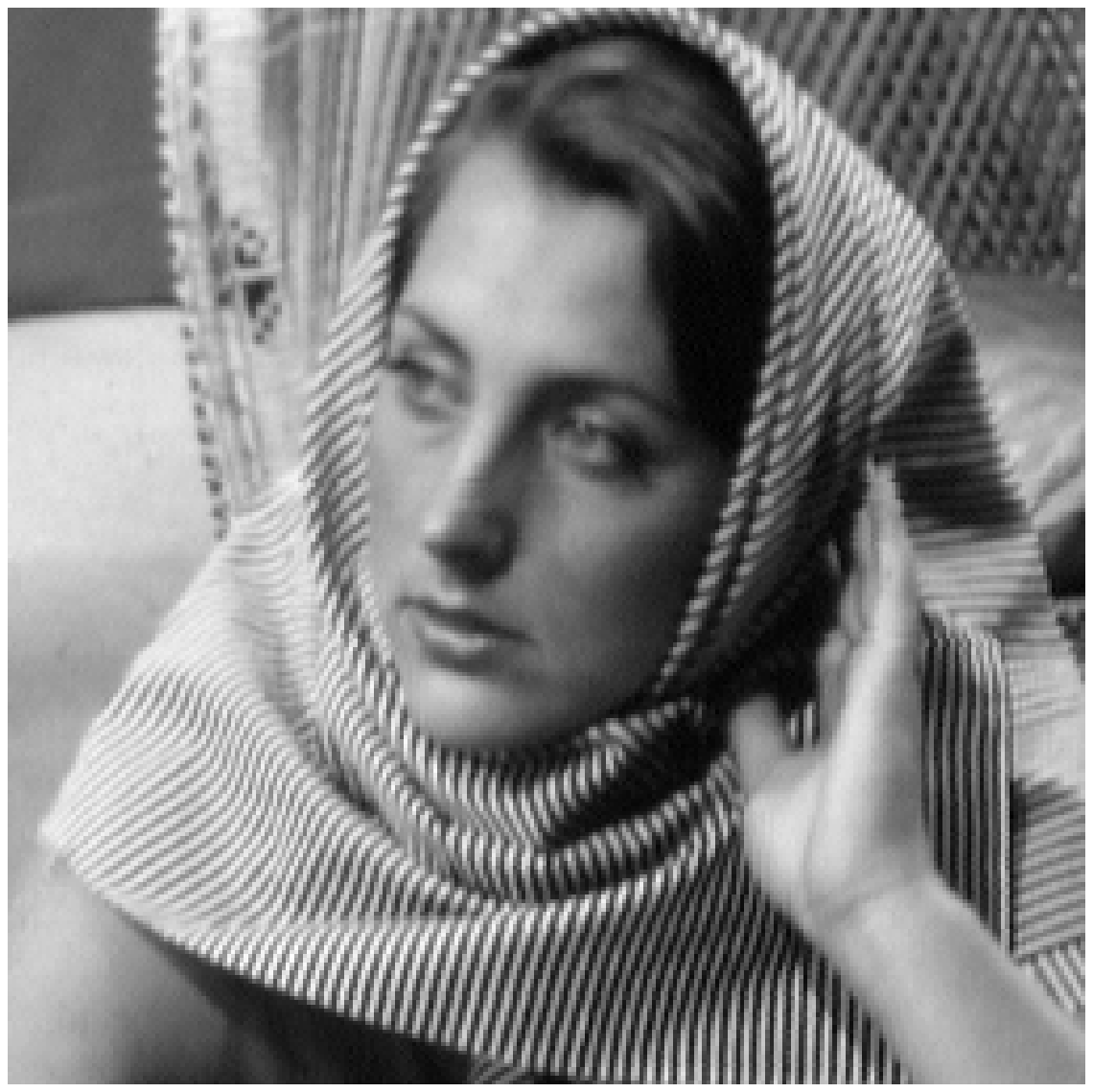}}
    \hspace{0.001in}
    \subfigure{\label{originaleinsta}
    \includegraphics[width=0.2\textwidth,clip]{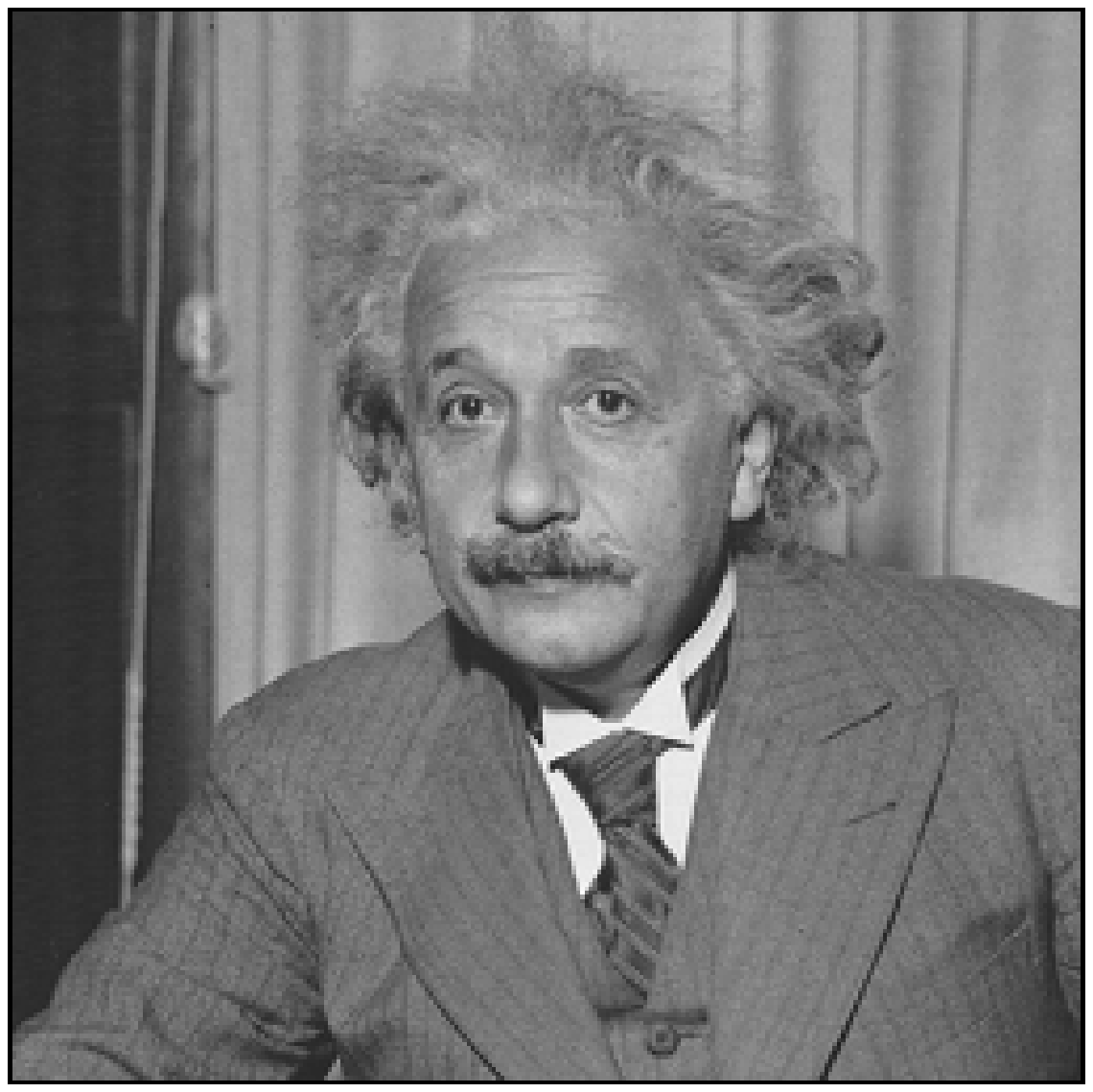}}
    \hspace{0.001in}
    \subfigure{\label{originalpepper}
    \includegraphics[width=0.2\textwidth,clip]{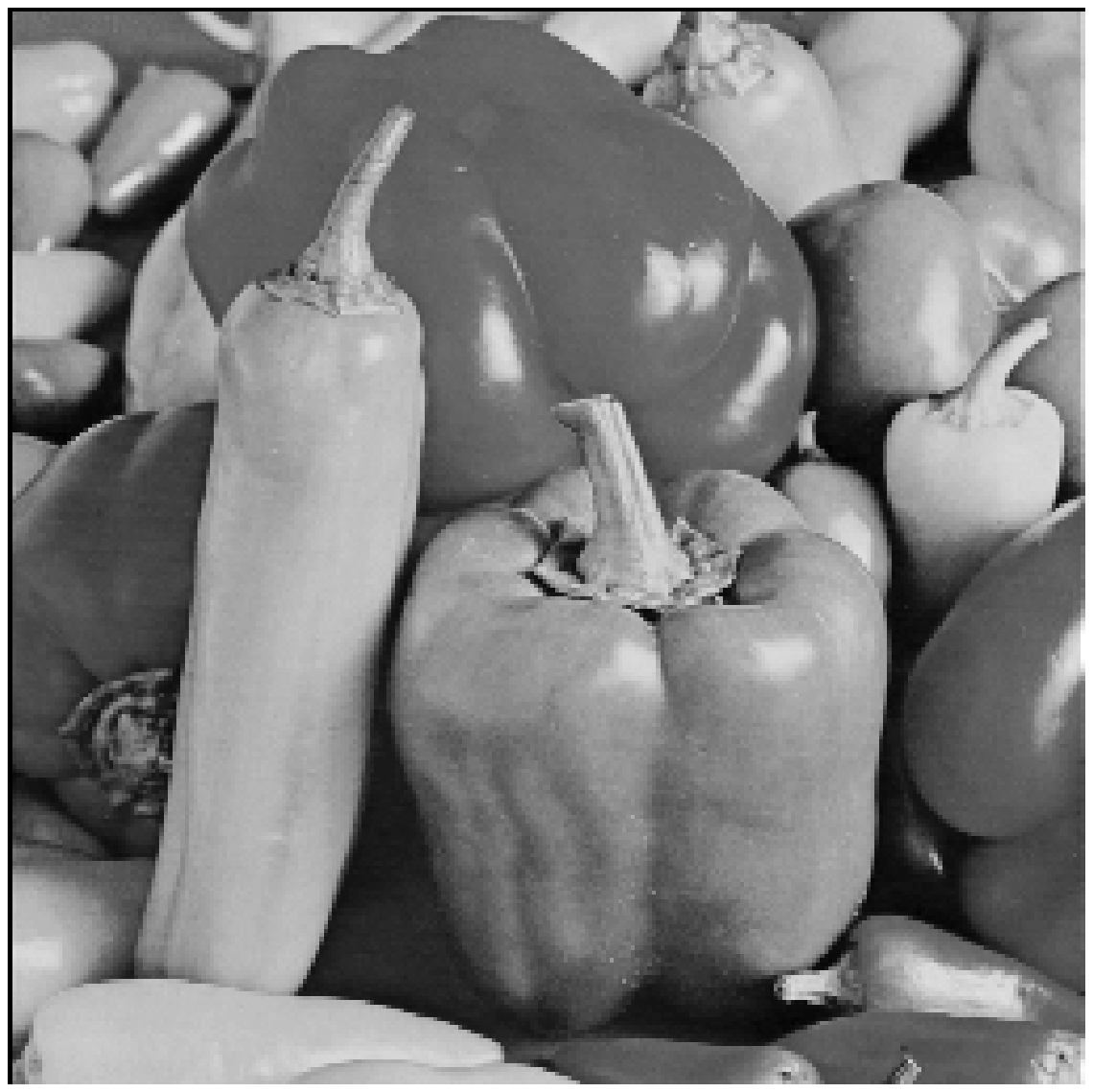}}
    \hspace{0.001in}
   \subfigure{\label{originalweats}
    \includegraphics[width=0.2\textwidth,clip]{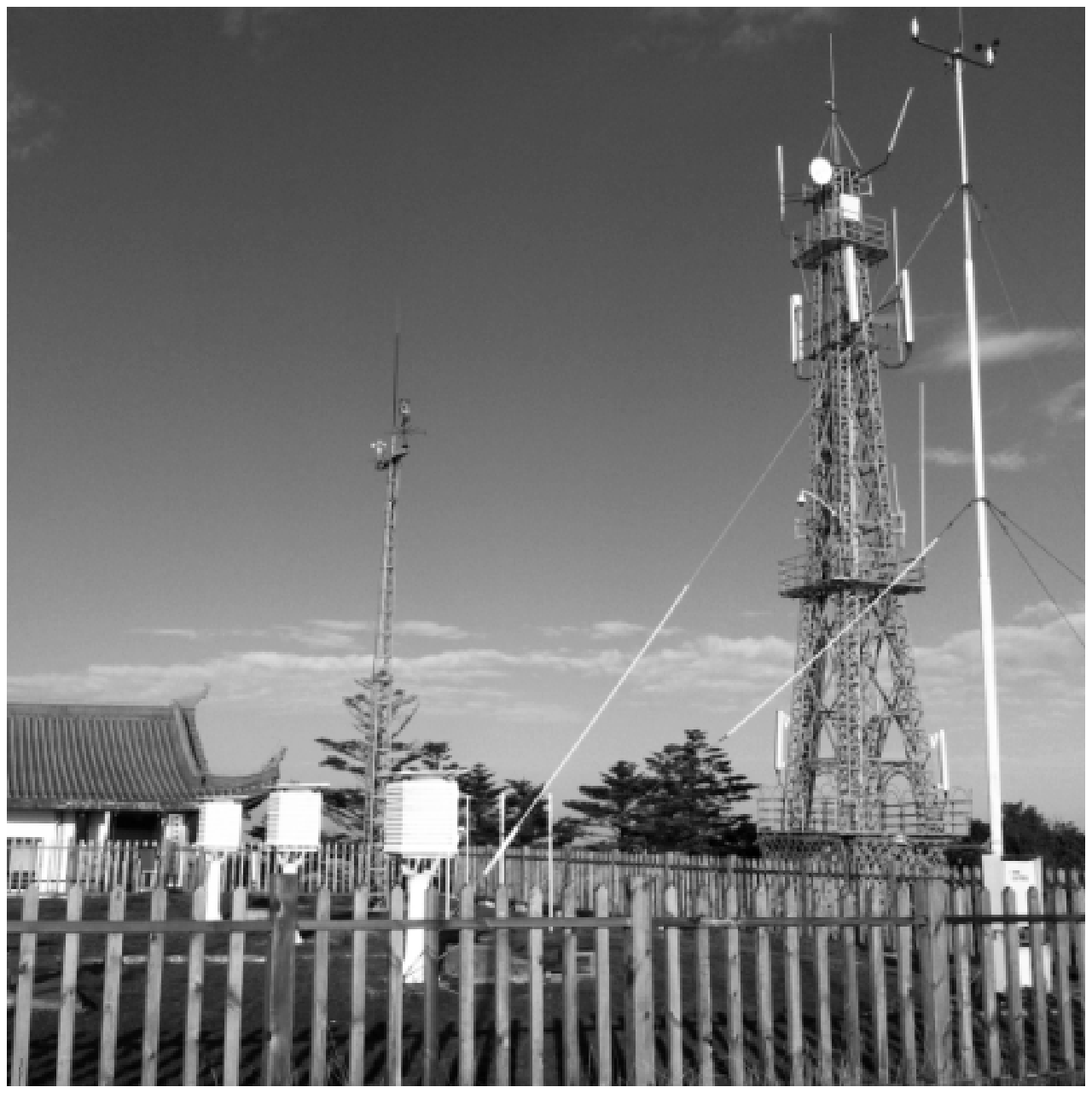}}
         \caption{\textit{Original images. Top row: from left to right, (a) Cameraman, (b) Satellite, (c) House, (d) Boat. Bottom row: from left to right, (e) Barbara, (f) Einstein, (g) Peppers, (h) Weather-station}.}
\label{oringinal}\end{figure}

All the test images are shown in Fig.~\ref{oringinal}, seven 256-by-256 images as: (a) Cameraman.tif, (b) Satellite.pgm, (c) House.png,
(d) Boat.pgm, (e) Barbara.tiff, (f) Einstein.pgm, (g) Peppers.png and one 460-by-460 image (h) Weatherstation.tif. For the sake of simplicity, the pixel values in all of our tests are lied in [0,1] which have been explained above.

We set the penalty parameters $\beta_1=1$, $\beta_2=500$, $\beta_3=1$,
and relax parameter $\gamma=1.618$
throughout all the experiments. And three blur kernels are generated by Matlab built-in function (i) \texttt{fspecial('gaussian',7,5)} for $7\times 7$ Gaussian blur with standard deviation 5, (ii) \texttt{fspecial('gaussian',15,5)} for $15\times 15$ Gaussian blur with standard deviation 5 and (iii) \texttt{fspecial('average',7)} for $7\times 7$  average blur. We
generate all blurring effects using the Matlab built-in function
\texttt{imfilter(I,psf,} \texttt{'circular','conv')} under periodic boundary conditions with ``\texttt{I}'' the original image and ``\texttt{psf}'' the blur kernel.
We generate all noise effects by Matlab built-in function \texttt{imnoise(B,'salt \& pepper',level)} with ``\texttt{B}'' the blurred image and fix the same random matrix for different methods. We only consider the salt-and-pepper noise in our experiments, since the variation method is easy to extend to the random value noise case.

\subsection{Study on the rest parameters}

\begin{figure}
  \centering
  \subfigure{\label{gau75noi4CPUK}
    \includegraphics[width=0.32\textwidth,clip]{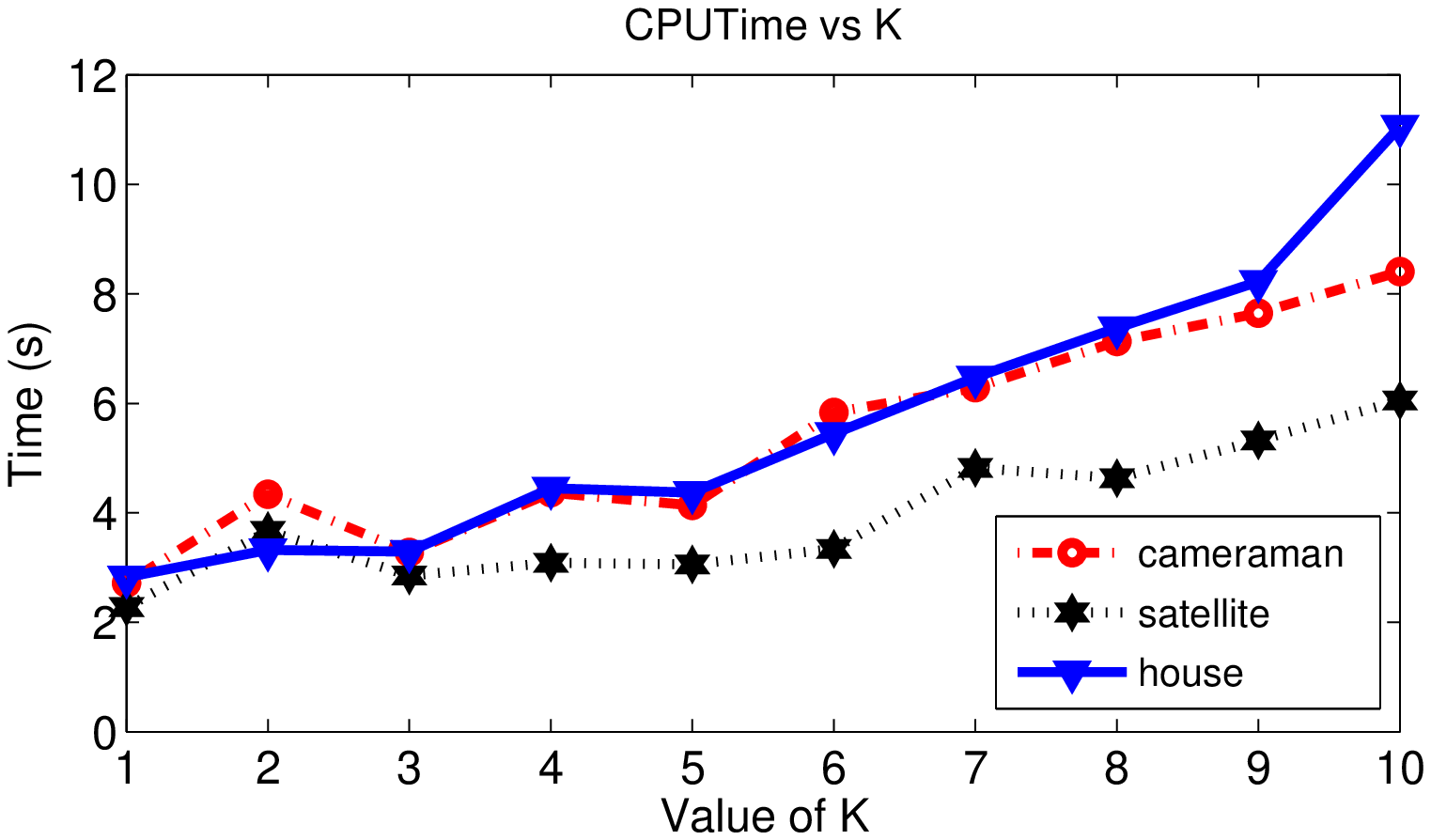}}
    \hspace{0.001in}
   \subfigure{\label{gau75noi4PSNRK}
    \includegraphics[width=0.3\textwidth,clip]{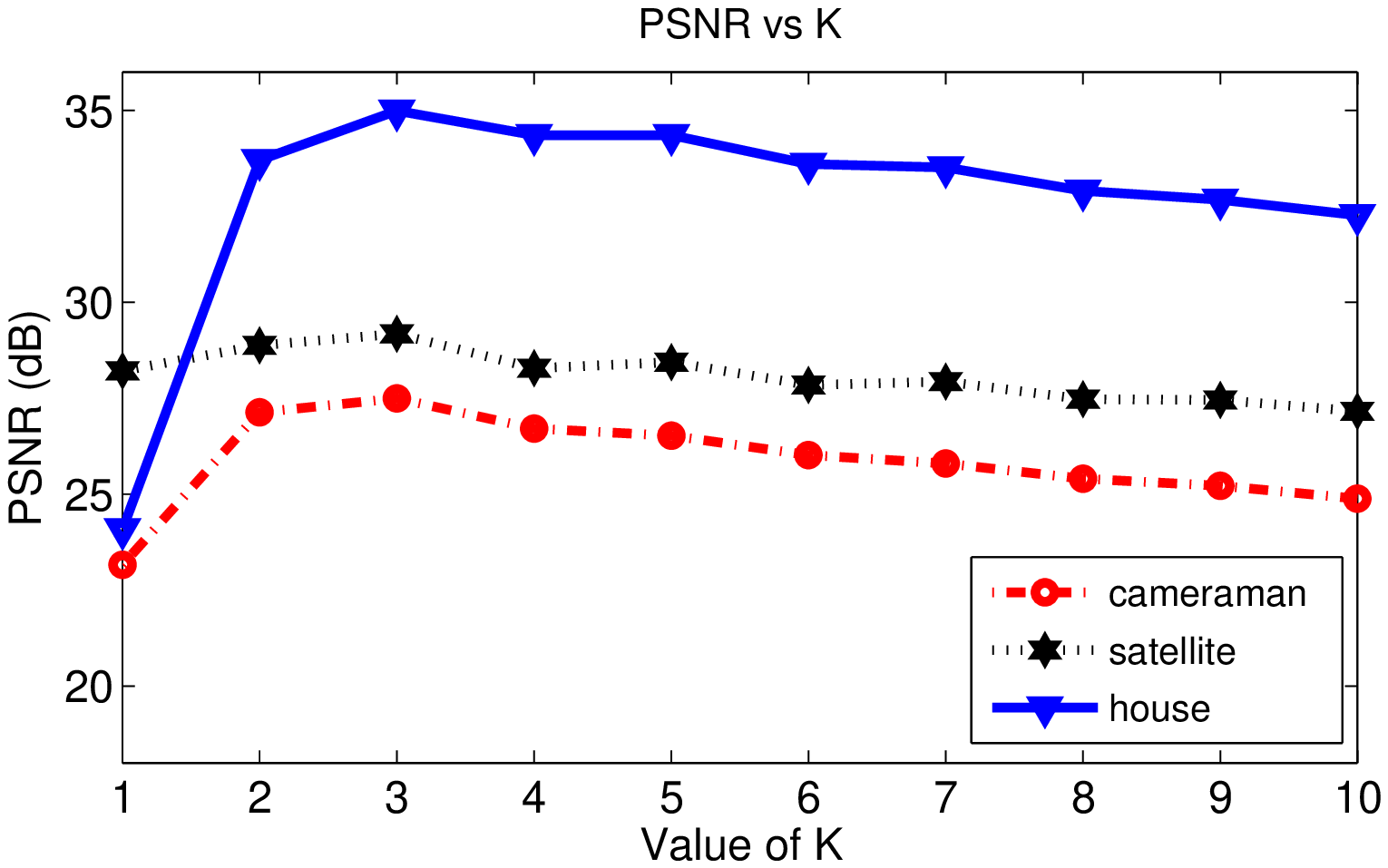}}
    \hspace{0.001in}
    \subfigure{\label{gau75noi4REFK}
    \includegraphics[width=0.32\textwidth,clip]{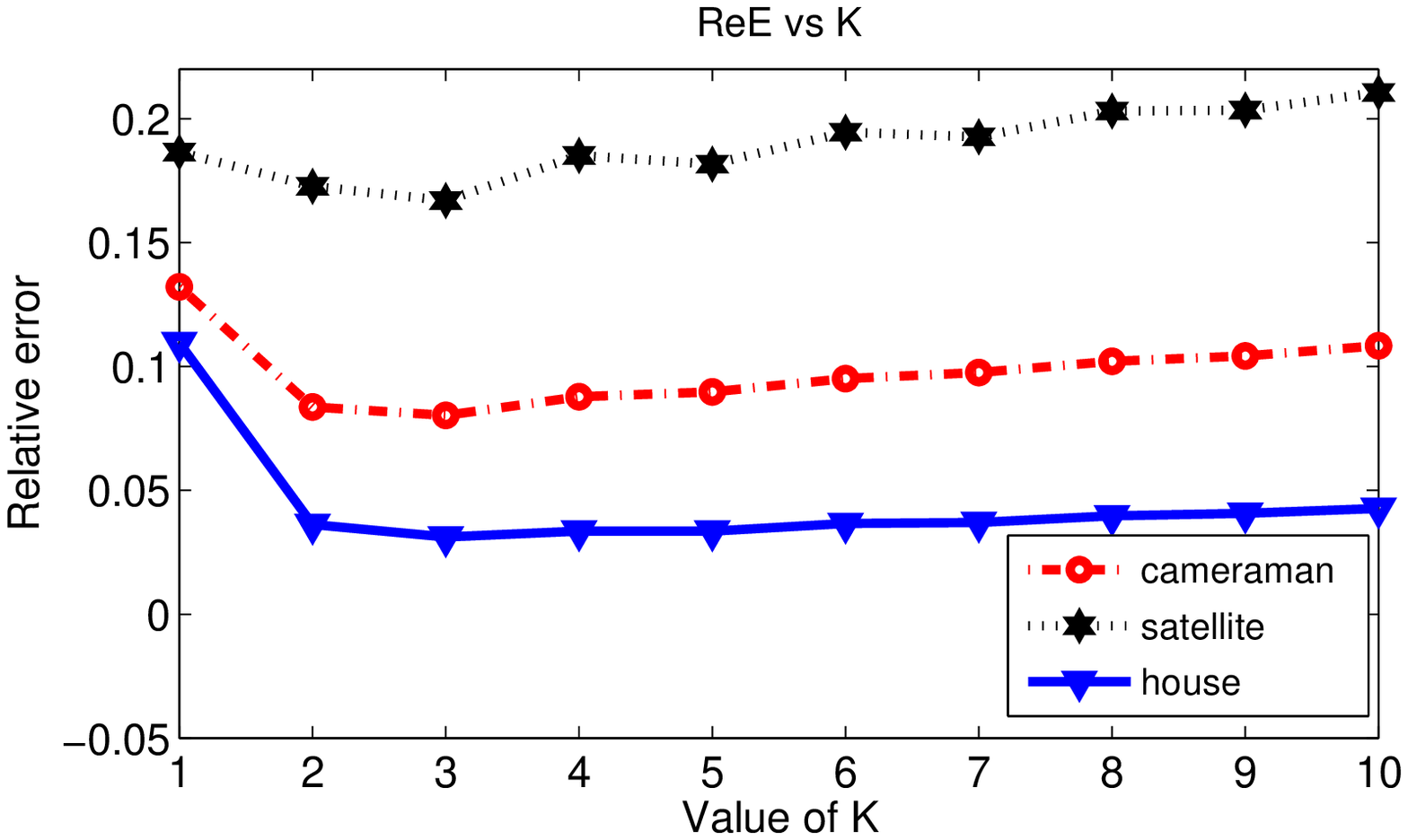}}
     \caption{\textit{Results of our proposed method depending on group size parameter $K$ on the image ``Cameraman'' blurred by
Gaussian blur kernel with $7\times 7$ and standard deviation 5 and corrupted by 40\% salt-and-pepper noise}.}
\label{gau75noi4K}\end{figure}

Firstly, we set the group size parameter $K = 3$ to find a good maximum inner iterations $NIt$. Our experiments are on the image ``Cameraman'' blurred by
Gaussian blur kernel with $7\times 7$ and standard deviation 5 and corrupted by 40\% salt-and-pepper noise. The results are shown in Table~\ref{gau75noi4NIt}.
From Table~\ref{gau75noi4NIt}, we can choose maximum inner iterations $NIt=5$ for the best. Then we fix $NIt=5$ and repeat more
experiments for choosing a good group size parameter $K$. We operate the three 256-by-256 images (a) ``Cameraman'', (b) ``Satellite'', and (c) ``House'' for this best option of parameter $K$.
The results are shown in Fig.~\ref{gau75noi4K}.
From the figure, we can see that $K=3$ is good for all the tests both on CPU time and PSNR. From now on,
we fix that $NIt=5$ and $K=3$ for the following experiments.
\begin{center}
\begin{table}[!hbp]
\renewcommand{\captionlabeldelim}{.}
\setlength{\abovecaptionskip}{0pt}
\setlength{\belowcaptionskip}{11pt} \centering \caption{\textit{PSNR (dB) and time (s) depending on maximum inner iterations $NIt$ on the
image ``Cameraman'' with Gaussian blur kernel $7\times 7$ and standard deviation 5 and 40\% salt-and-pepper noise}.} \centering
\begin{tabular}{|c|c|c|c|c|c|c|c|c|c|c|}
\hline 
  $NIt$               &  1  & 3 & 5 & 7  & 10 & 20 & 50 & 100 & 200 & 1000    \\
\hline
 PNSR & 21.38&27.44&27.50&27.36 &27.22&27.00&26.87&26.85&26.84& 26.83    \\
\hline
Time& 3.245 &2.868&3.267&4.274&6.022&10.87&23.84&45.54&89.30&1965   \\
\hline
\end{tabular}
\label{gau75noi4NIt}\end{table}
\end{center}
\begin{figure}
  \centering
  \subfigure{
    \includegraphics[width=0.3\textwidth,clip]{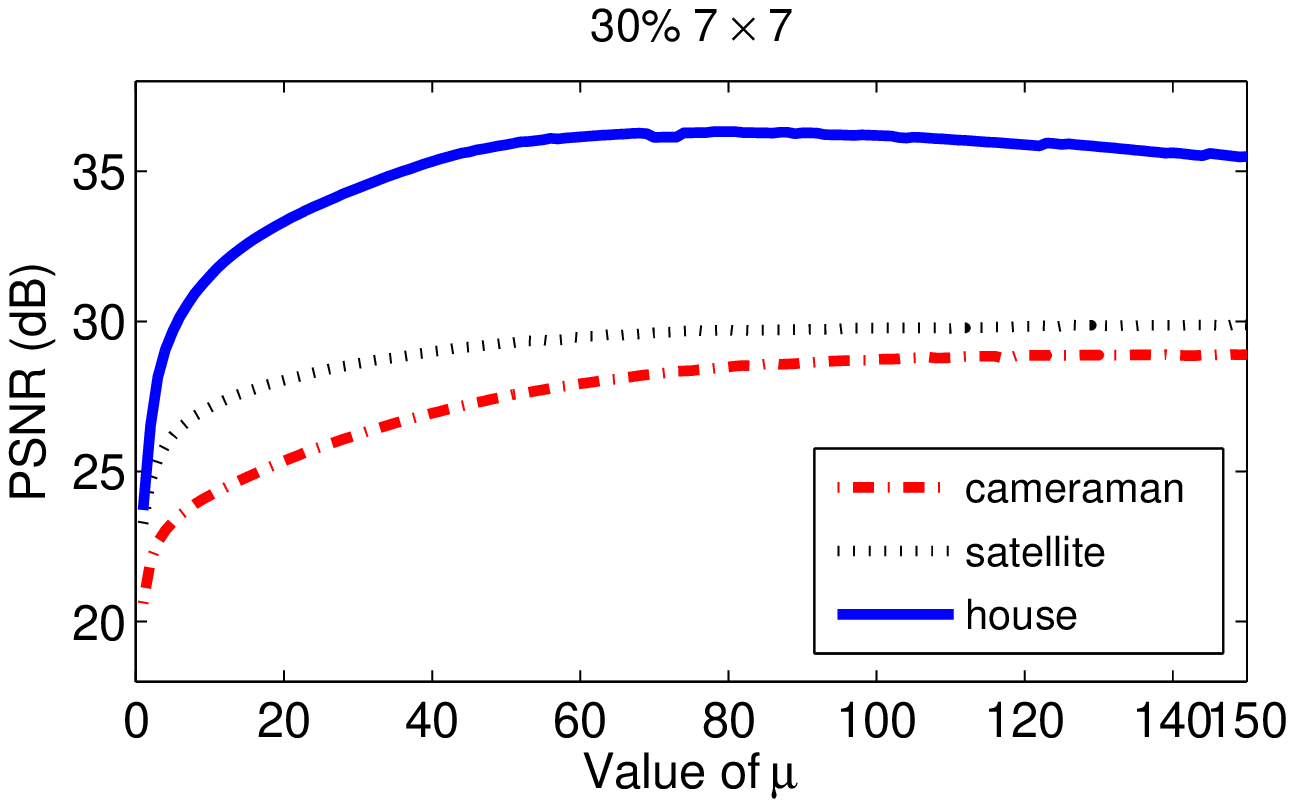}}
    \hspace{0.001in}
    \subfigure{
    \includegraphics[width=0.3\textwidth,clip]{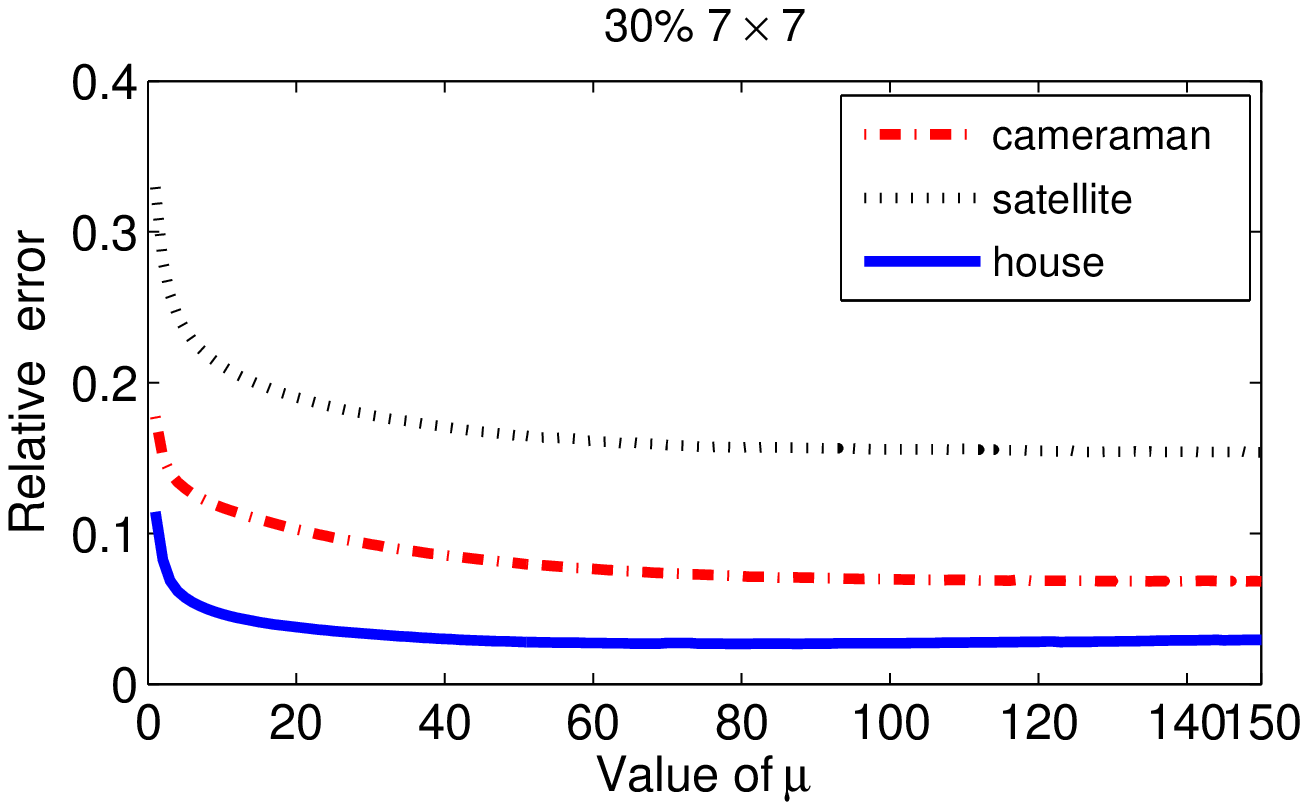}}
    \hspace{0.001in}
   \subfigure{
    \includegraphics[width=0.3\textwidth,clip]{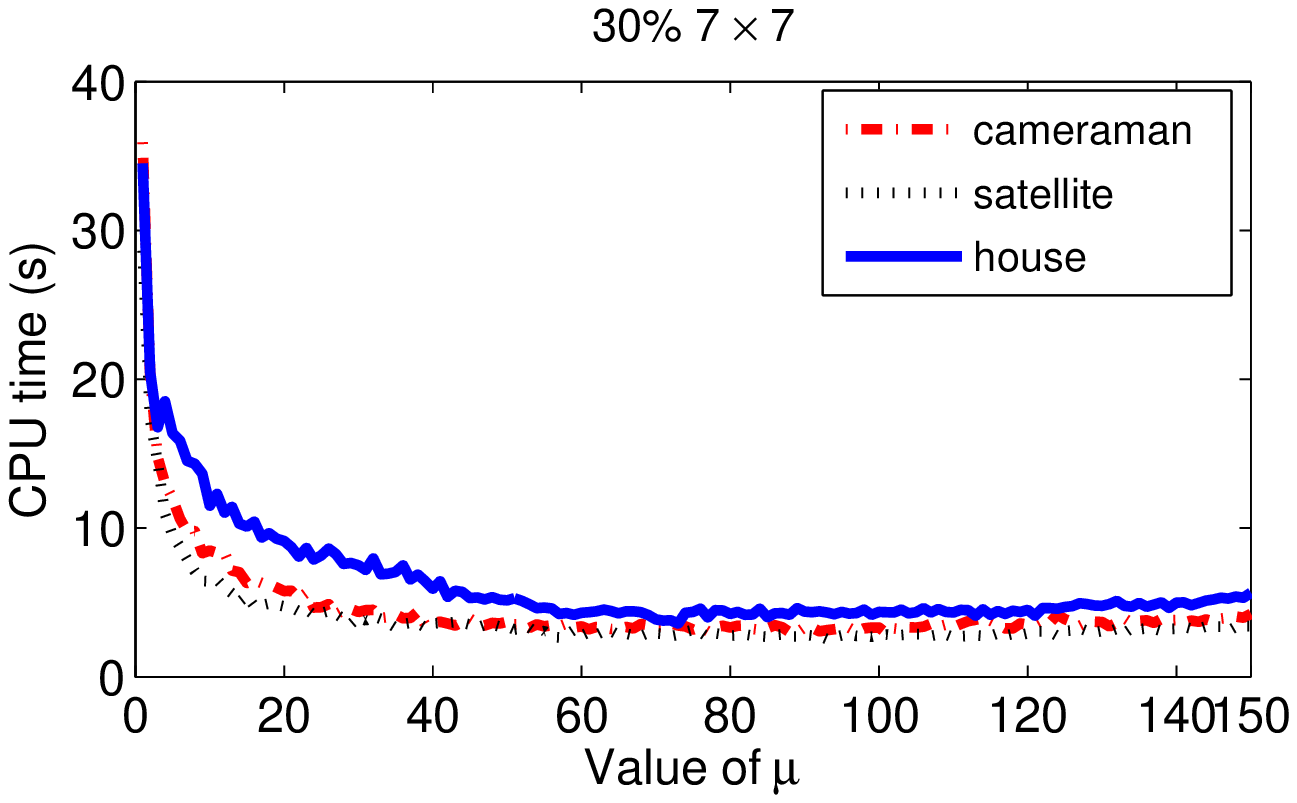}}
    \hspace{0.001in}\\
    \subfigure{
    \includegraphics[width=0.3\textwidth,clip]{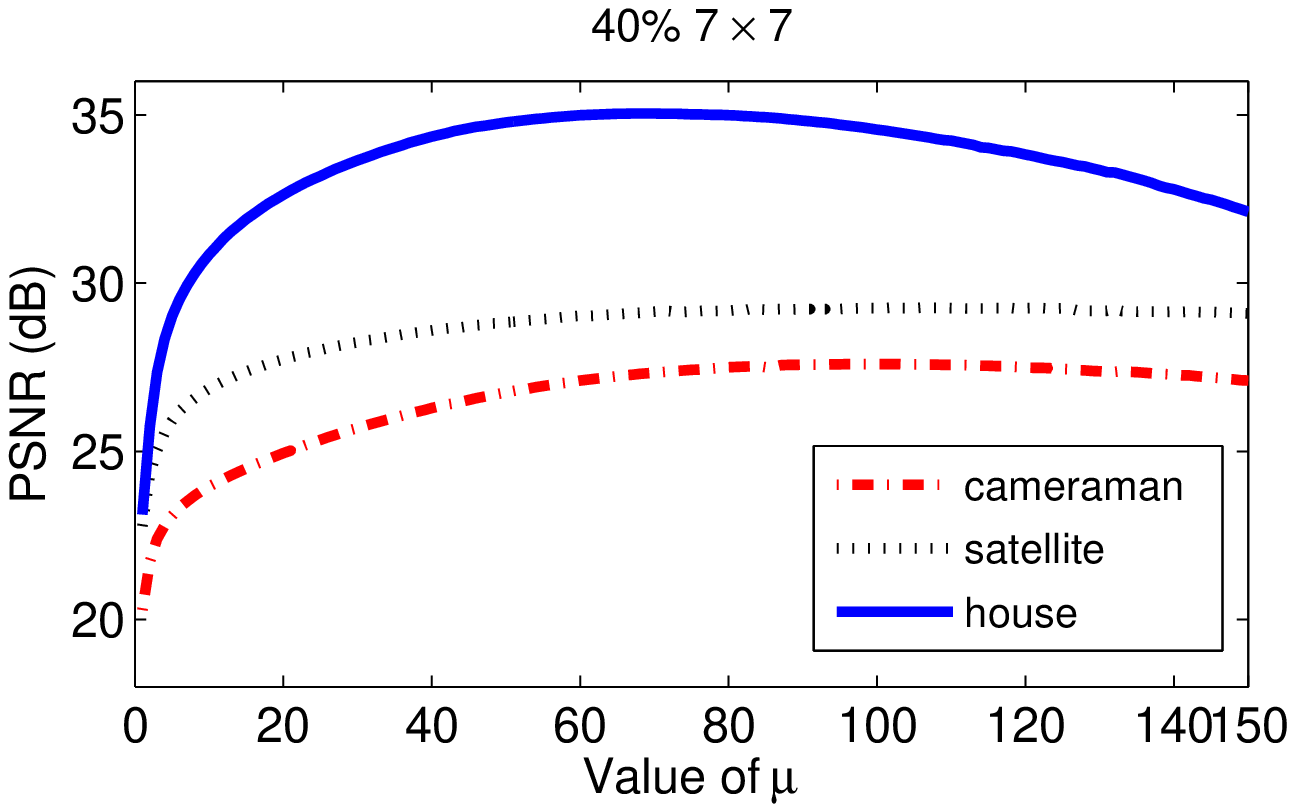}}
    \hspace{0.001in}
     \subfigure{
    \includegraphics[width=0.3\textwidth,clip]{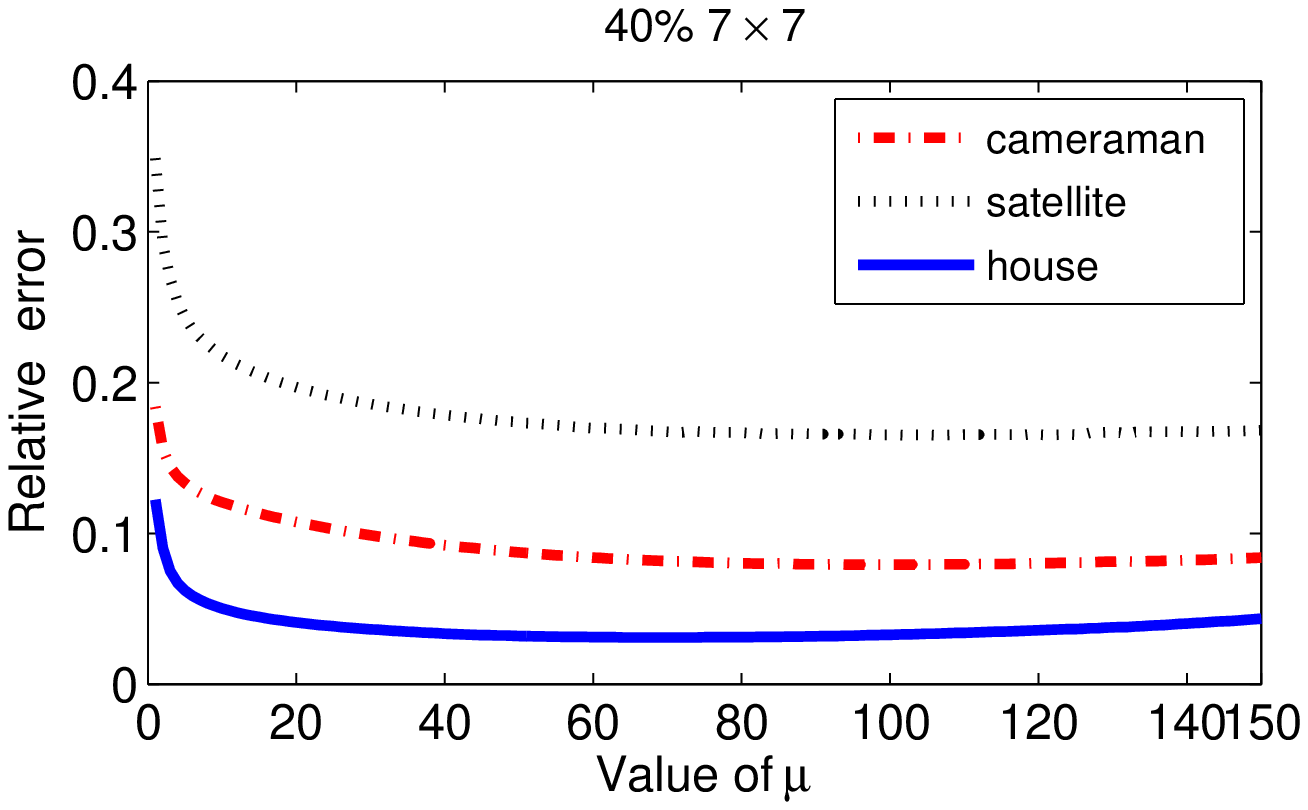}}
    \hspace{0.001in}
   \subfigure{
    \includegraphics[width=0.3\textwidth,clip]{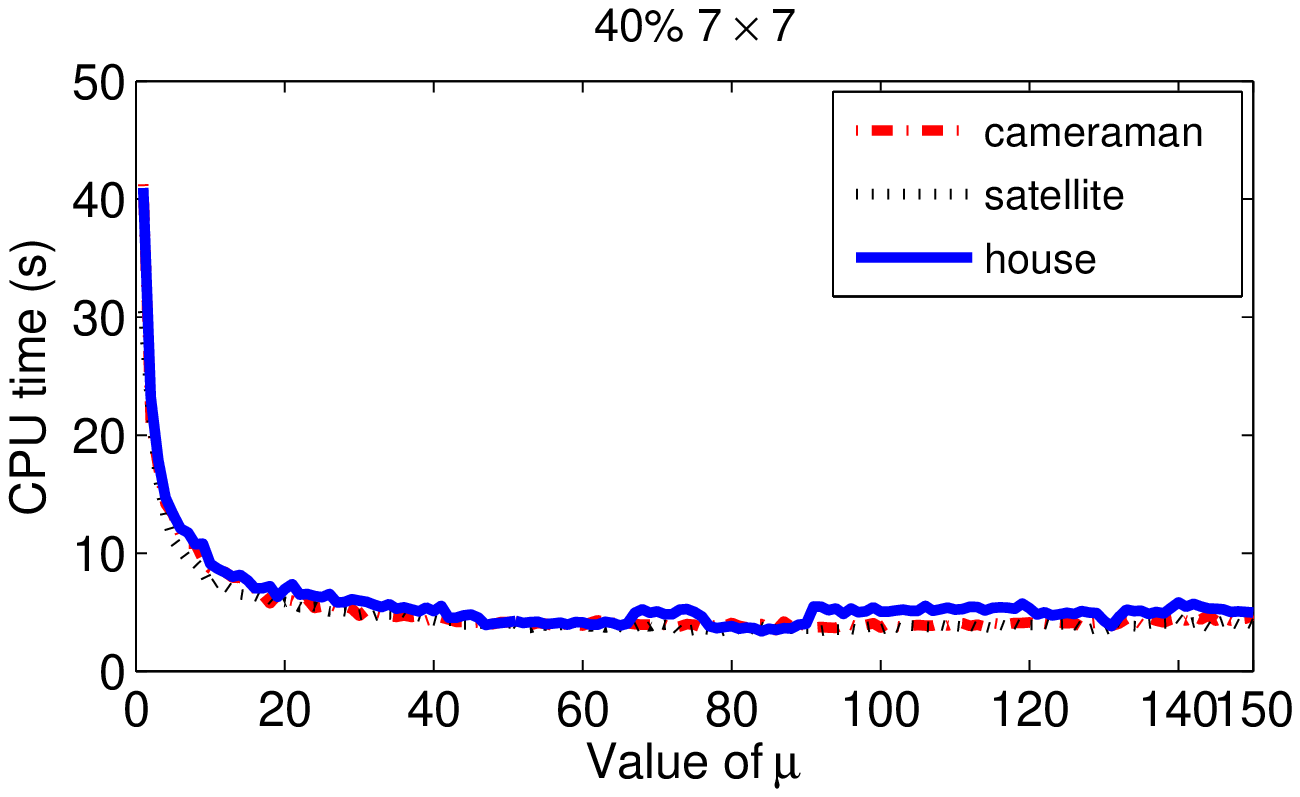}}
    \hspace{0.001in}\\
    \subfigure{
    \includegraphics[width=0.3\textwidth,clip]{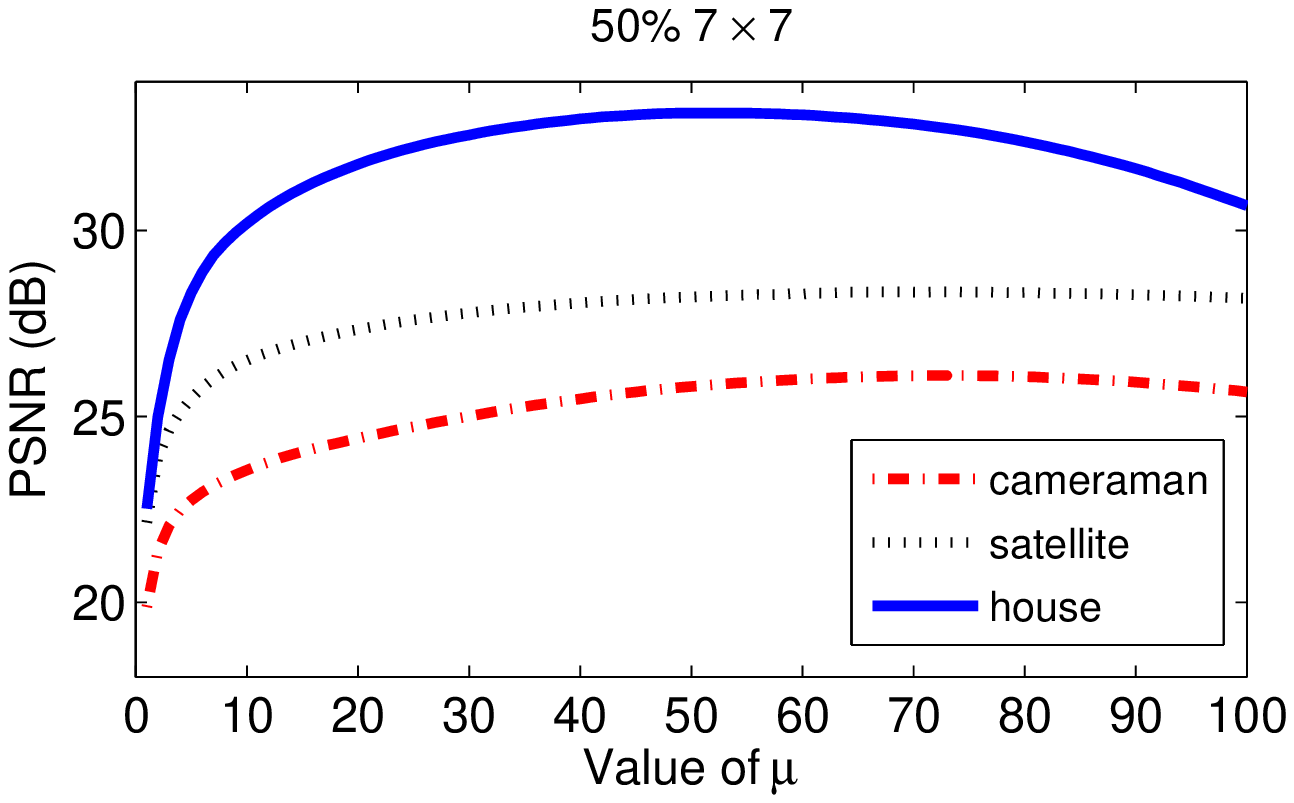}}
    \hspace{0.001in}
        \subfigure{
    \includegraphics[width=0.3\textwidth,clip]{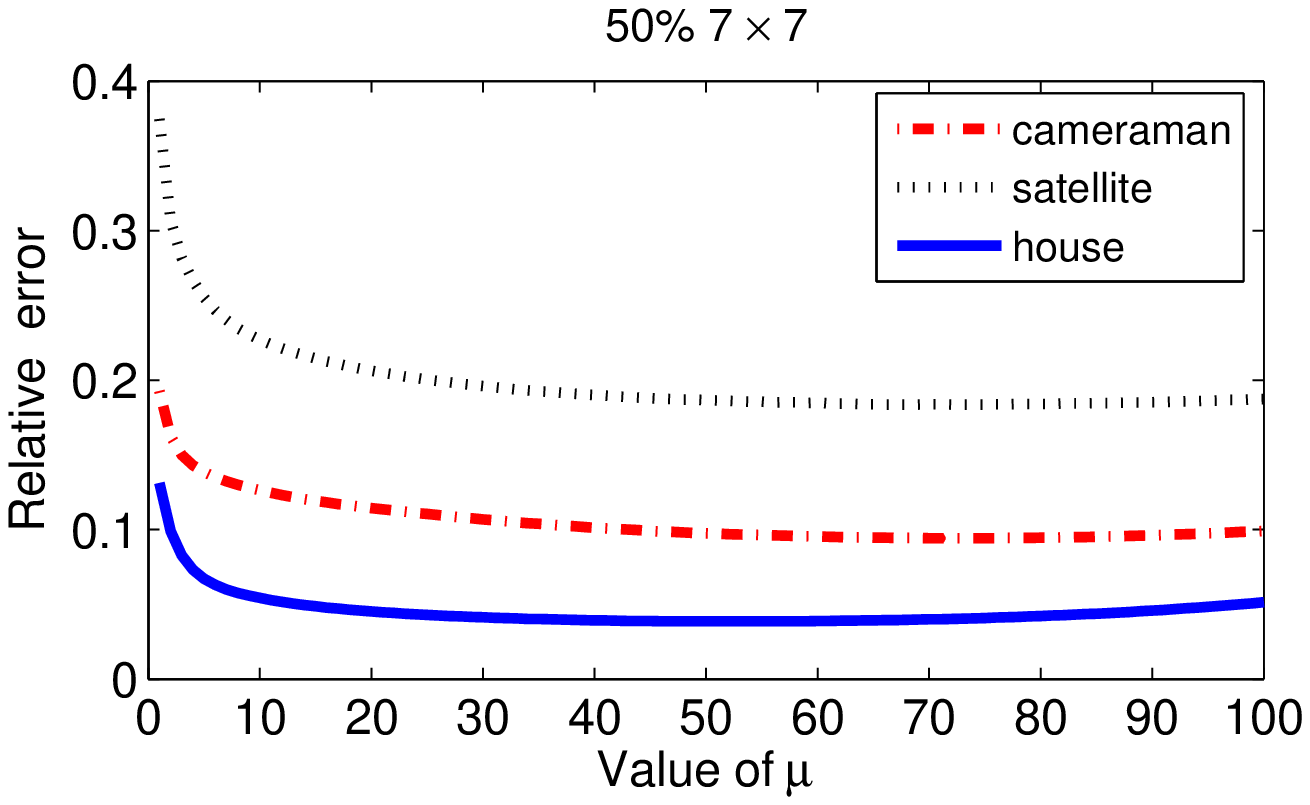}}
    \hspace{0.001in}
    \subfigure{
    \includegraphics[width=0.3\textwidth,clip]{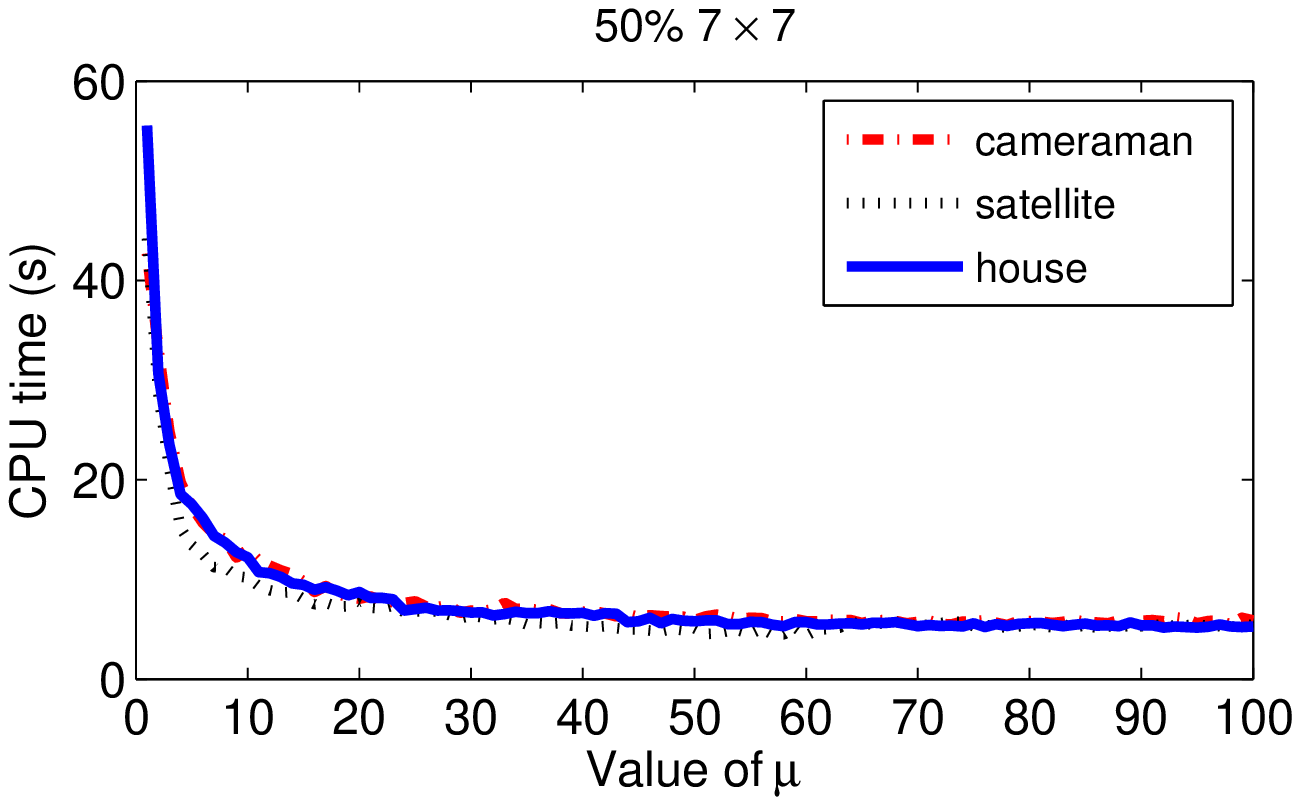}}
    \hspace{0.001in}\\
    \subfigure{
    \includegraphics[width=0.3\textwidth,clip]{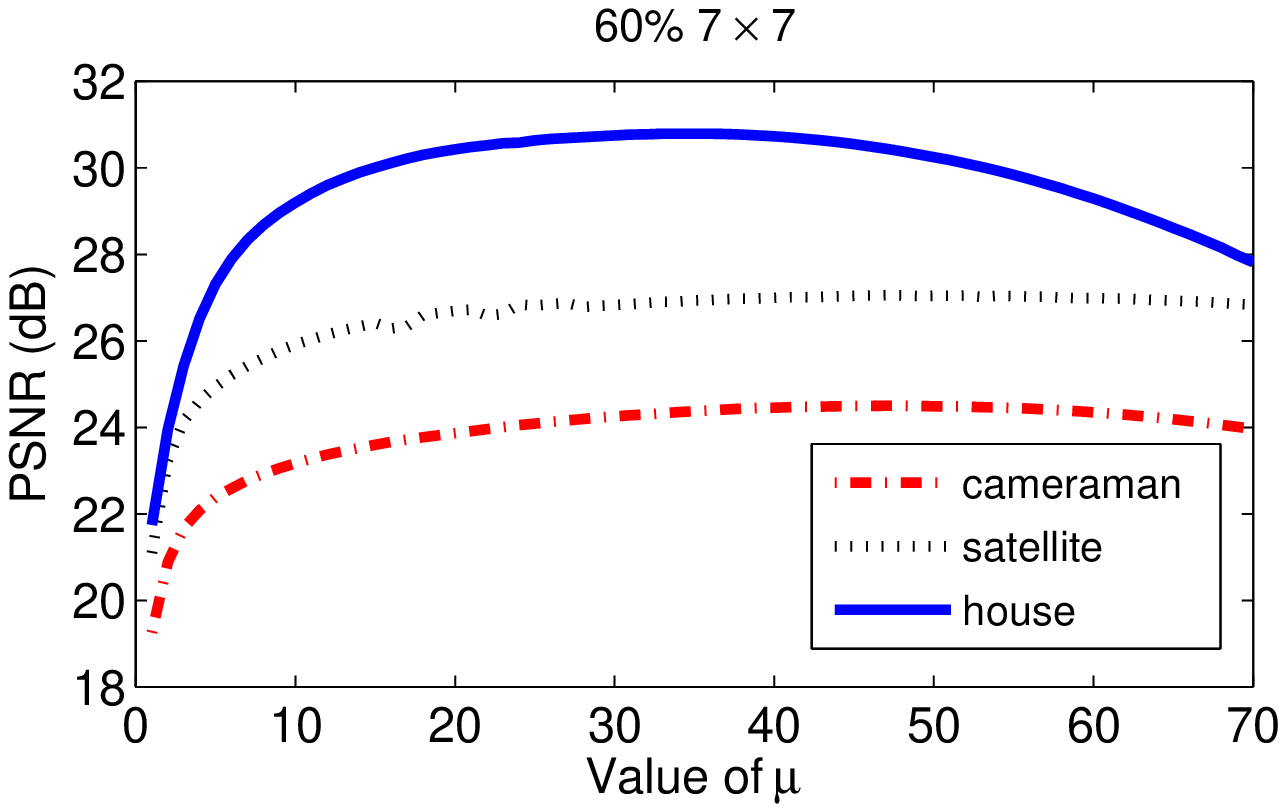}}
    \hspace{0.001in}
       \subfigure{
    \includegraphics[width=0.3\textwidth,clip]{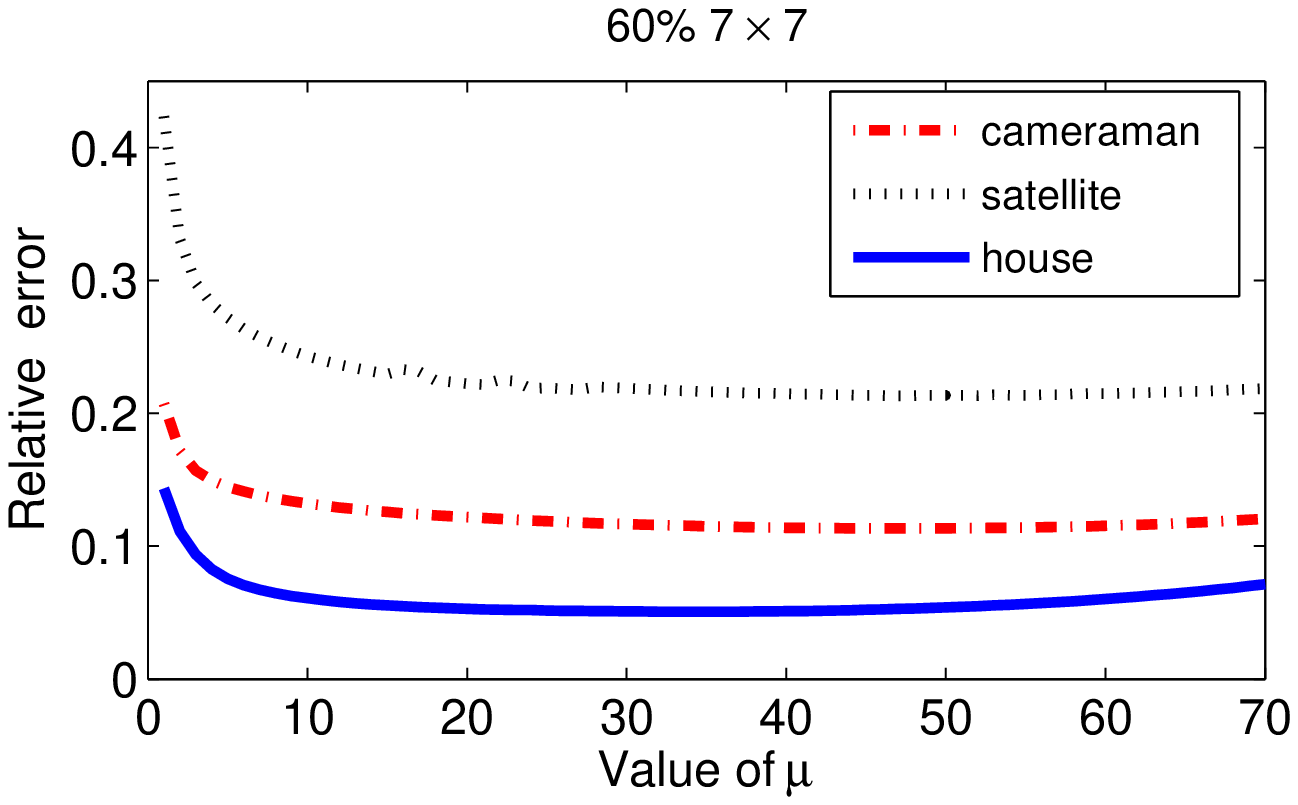}}
    \hspace{0.001in}
   \subfigure{
    \includegraphics[width=0.3\textwidth,clip]{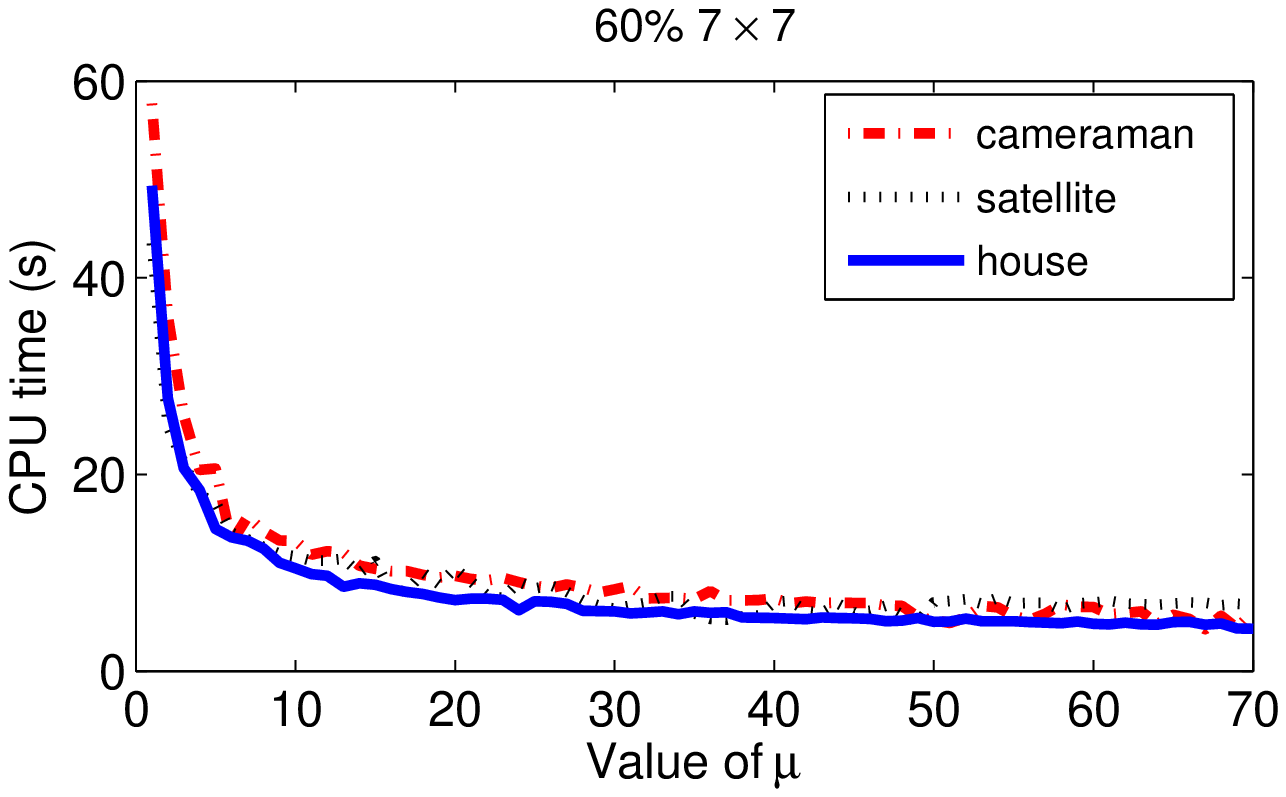}}
    \hspace{0.001in}
    \caption{\textit{Results of our proposed method depending on regularization parameter $\mu$ on the images (a) (b) and (c) that are blurred by
 $7\times 7$ Gaussian blur kernel with standard deviation 5 and corrupted by salt-and-pepper noise from 30\% to 60\%}.}
\label{gau75noi3-6mu}\end{figure}

Then, we test how to select a good regularization parameter $\mu$ for different images. We will point out several important advantages of our method in the following experiments. For the sake of simplicity, we focus on the above three test 256-by-256 images (a) (b) and (c). Under the Gaussian blur with $7\times 7$ window size and standard deviation 5, the images corrupted by added salt-and-pepper noise from 30\% to 60\% are tested.
In Fig.~\ref{gau75noi3-6mu}, we plot PSNR, ReE, and Time for our algorithm against different values of the regularization parameter $\mu$. Each row in Fig.~\ref{gau75noi3-6mu} corresponds to the four salt-and-pepper noise levels. In fact, for all $\mu$, our method always gives high PSNR values. Moreover, the PSNR curves of our method are very flat, which shows that our method is stable for a wide range of $\mu$, which is wider than that in \cite{CTY2013}. That is to say, our method is more robust than the method in \cite{CTY2013}. This is the first advantage of our method.

In addition, to the best of our knowledge, a good image restoration algorithm should satisfy the following two properties. a) It is fast and can reach
good results in term of both numerical values and high visual quality. b) It is not sensitive to parameters. Our method meets the requirement of these two properties, as
 it obtains good restoration results with the same parameters for different images under the same blur and noise. This is the second advantage of our method. Here and in the following experiments, under the same blur and noise, we choose the same parameters for all the test images. Particularly, for the images under the Gaussian blur with $7\times 7$ window size and standard deviation 5 and corrupted by
salt-and-pepper noise from 30\% to 60\%, we set $\mu = 100, 80,60,40$ respectively. After similar tests as Fig.~\ref{gau75noi3-6mu}, we list all the selection rule of $\mu$: for the images under the Gaussian blur with 15-by-15 window size and standard deviation 5 and corrupted by salt-and-pepper noise from 30\% to 60\%, $\mu = 120, 110,100,90$ respectively; for the images under the average blur with $7\times 7$ window size and corrupted by salt-and-pepper noise from 30\% to 60\%, $\mu = 100, 80,60,40$ respectively.

\subsection{Comparison with CTY and GLN for the test image ``Cameraman''}

In this subsection, we mainly compare our proposed method to CTY and GLN for deblurring problems under salt-and-pepper noise. We use image ``Cameraman'' for the experiments in this subsection. Our purposes are (1) to show the improvement of PSNR and to demonstrate the efficiency of our proposed method mainly via a comparison with CTY and GLN, and (2) to illustrate that our proposed method can overcome the staircase effects effectively and get better visual quality than CTY and GLN, which is the third advantage of our method.

Firstly, we generate the blurred images by two Gaussian blurs (i) and (ii) with periodic boundary conditions as mentioned above, and then corrupt the blurred images by salt-and-pepper noise from 30\% to 60\%. For CTY and GLN, we have tuned the parameters manually to give the best PSNR improvement. The numerical results by the three methods are shown in Table~\ref{cameranmantol1CA}. From the table, we see that both our proposed method and CTY are much faster and can get higher PSNR than GLN. Our proposed method needs the fewest iterations than the other two methods, and the time is always close to CTY. 
\begin{table}[!hbp]
\setlength{\abovecaptionskip}{0pt}
\setlength{\belowcaptionskip}{11pt} \centering \caption{\textit{Numerical comparison of
the fast $\ell_1$-TV method (GLN) \cite{GLNG2009}, the ADM2CTVL1 method (CTY) \cite{CTY2013}, and our proposed method (Ours) under two Gaussian blurs (Bls) (i) {\rm \texttt{fspecial('gaussian',7,5)}} and (ii) {\rm \texttt{fspecial('gaussian',15,5)}} and corrupted by salt-and-pepper noise from 30\% to 60\%}.} \centering
\begin{tabular}{|c|c|r|r|r|}
\hline 
  \multirow{2}{*}{Bls} & Noise & \multicolumn{1}{|c|}{GLN} & \multicolumn{1}{|c|}{ CTY }&  \multicolumn{1}{|c|}{Ours} \\
 \cline{3-5}
     &                 level              & \multicolumn{1}{|r|}{Itrs/PSNR/Time/ReE} & \multicolumn{1}{|r|}{Itrs/PSNR/Time/ReE}& \multicolumn{1}{|r|}{Itrs/PSNR/Time/ReE}  \\
 \hline
  \multirow{4}{*}{(i)}&  30\% & 200/27.34/8.30/0.0893&122/27.66/4.91/0.0787&38/28.73/3.24/0.0696\\
    & 40\%&200/25.93/8.11/0.0961	&102/26.63/3.91/0.0887	&43/27.50/3.59/0.0802 \\
        & 50\%&200/24.73/8.13/0.1103	&77/25.42/2.90/0.1019	&49/26.00/3.92/0.0953\\
       & 60\%&200/23.43/8.38/0.1281	&56/24.20/2.51/0.1178	&62/24.50/4.96/0.1137\\
\hline
  \multirow{4}{*}{(ii)} &  30\% &200/24.37/8.16/0.1155	&113/24.10/3.68/0.1183	&37/24.52/2.92/0.1132 \\
       & 40\%&200/23.76/8.05/0.1235	&95/23.86/3.57/0.1221	&35/24.22/2.88/0.1169 \\
        & 50\%&200/23.46/8.11/0.1277	&61/23.55/2.56/0.1246	&35/23.93/2.91/0.1210   \\
       & 60\%&200/22.05/8.51/0.1502&55/23.12/2.17/0.1327	&36/23.48/2.73/0.1274 \\
\hline
\end{tabular}
\label{cameranmantol1CA}\end{table}
\begin{figure}
  \centering
  \subfigure{
    \includegraphics[width=0.2\textwidth,clip]{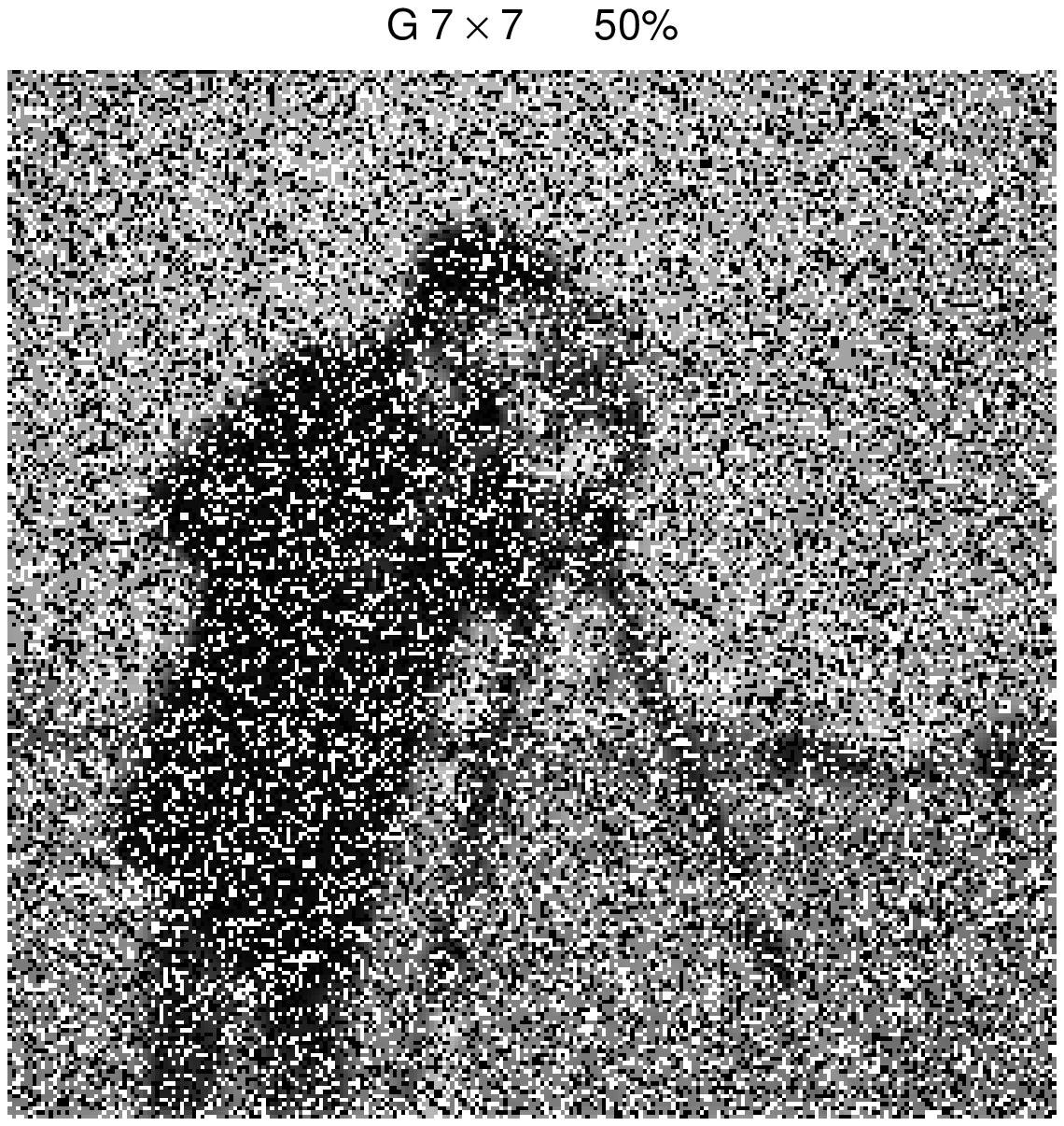}}
    \hspace{0.001in}
   \subfigure{
    \includegraphics[width=0.2\textwidth,clip]{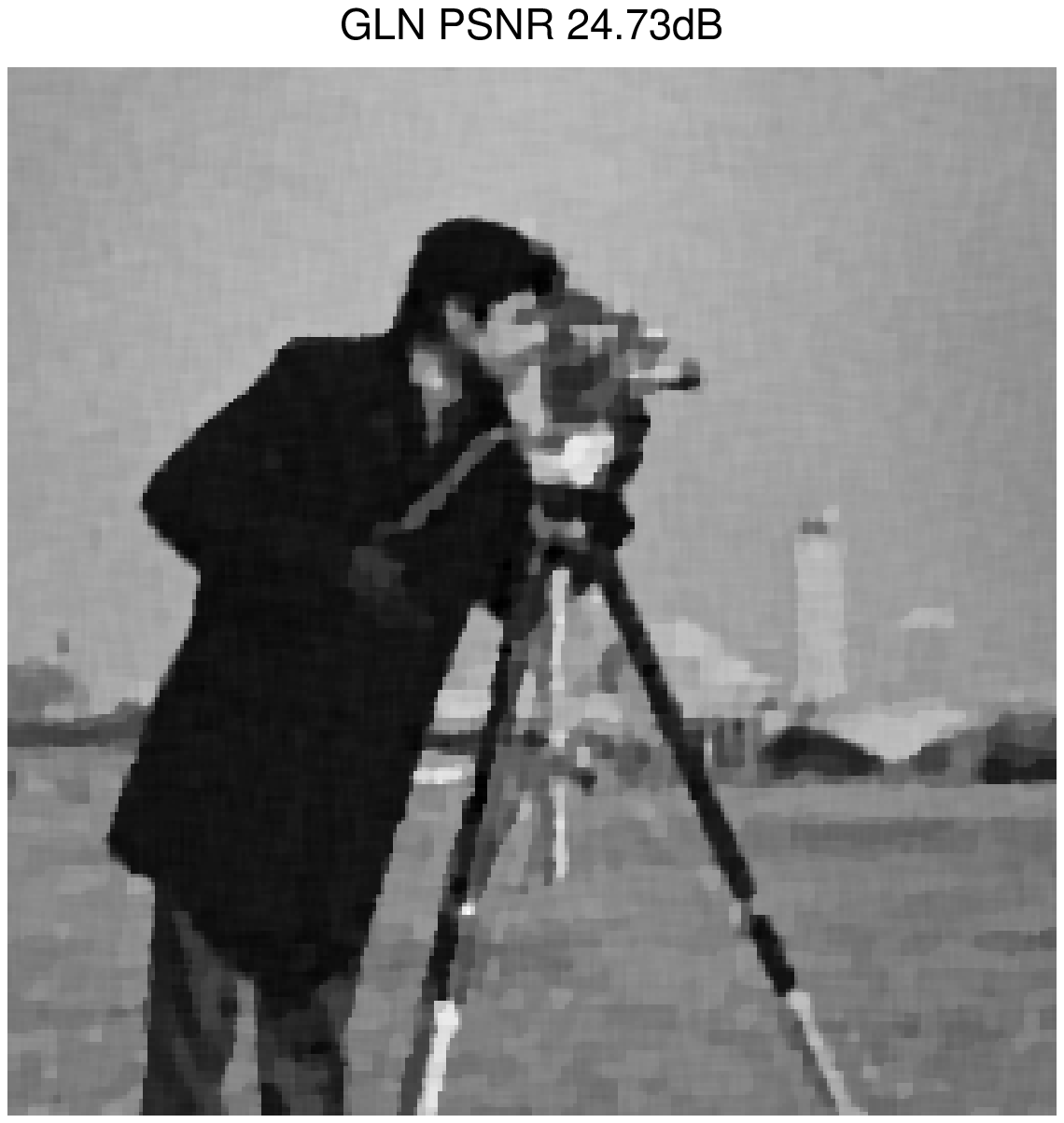}}
    \hspace{0.001in}
    \subfigure{
    \includegraphics[width=0.2\textwidth,clip]{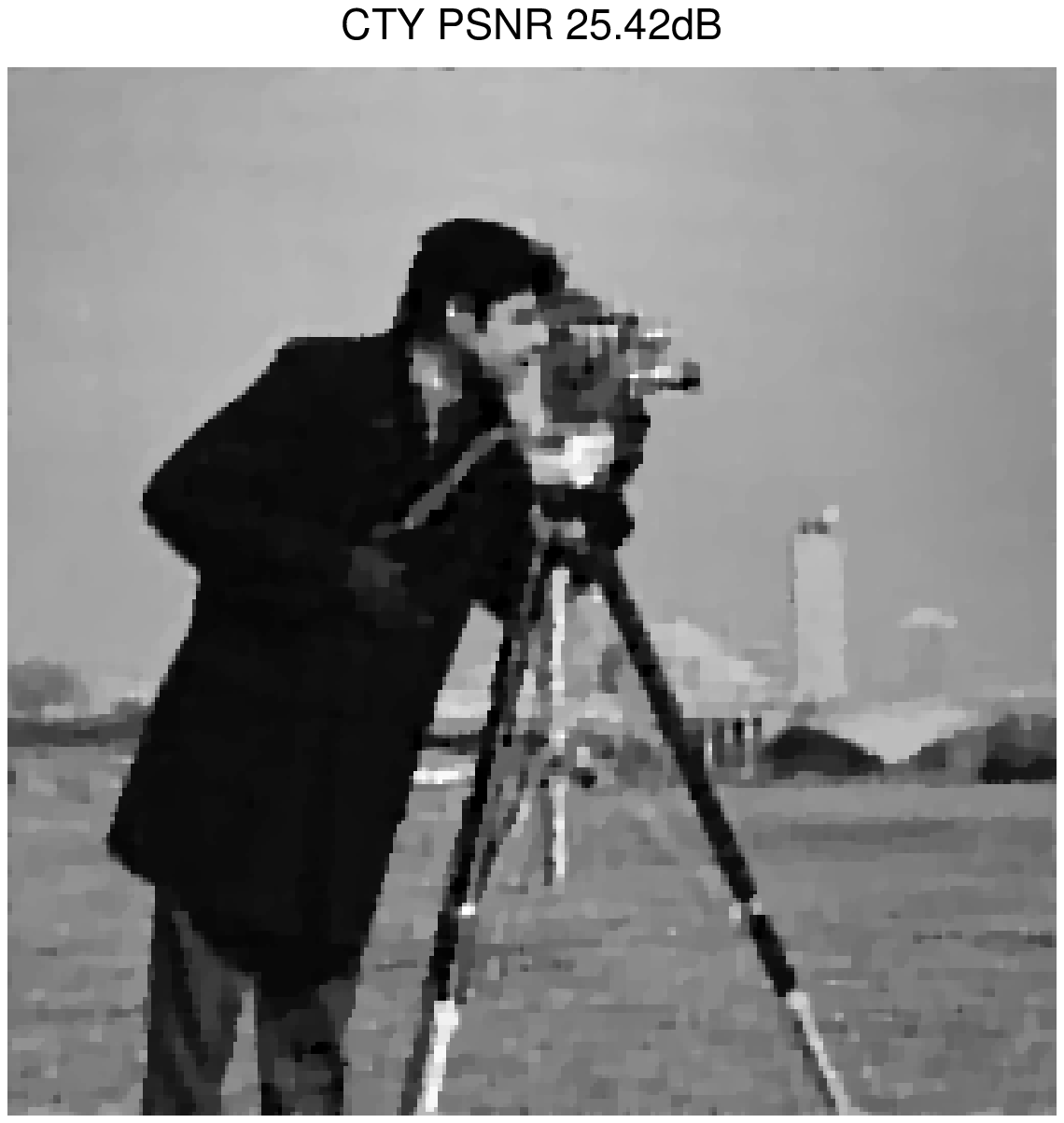}}
    \hspace{0.001in}
    \subfigure{
    \includegraphics[width=0.2\textwidth,clip]{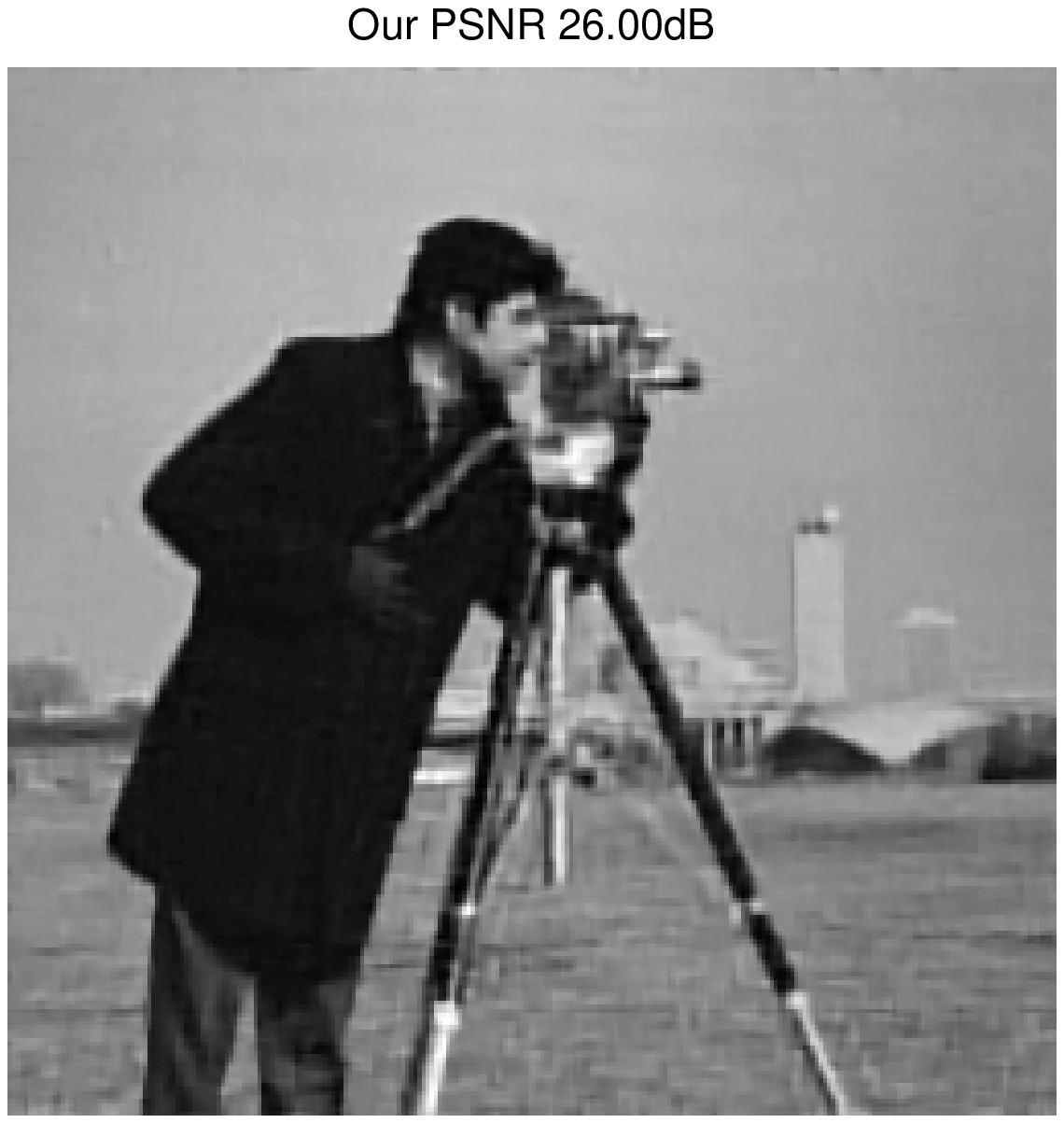}}
    \hspace{0.001in}
   \subfigure{
    \includegraphics[width=0.2\textwidth,clip]{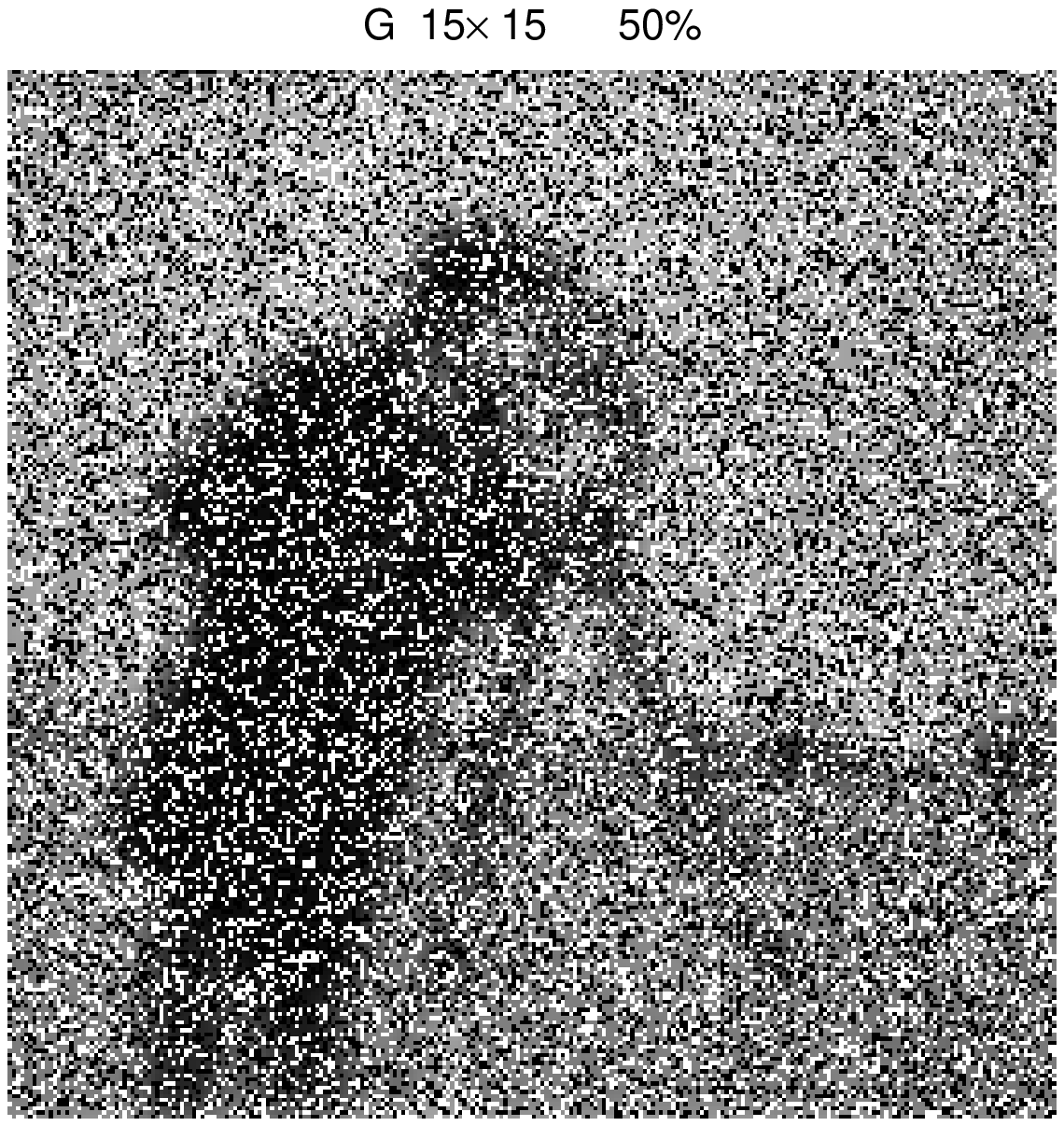}}
    \hspace{0.001in}
    \subfigure{
    \includegraphics[width=0.2\textwidth,clip]{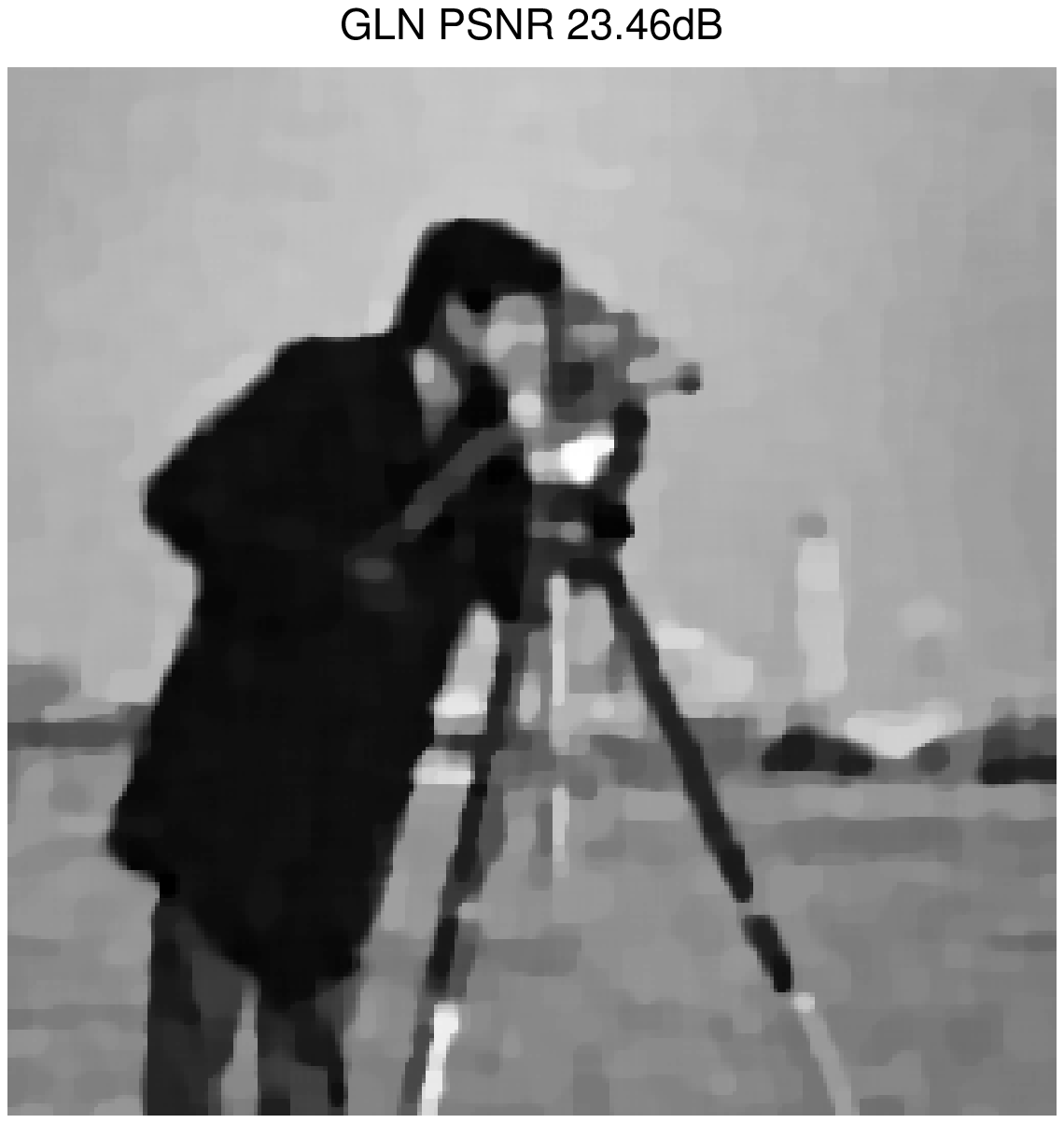}}
    \hspace{0.001in}
    \subfigure{
    \includegraphics[width=0.2\textwidth,clip]{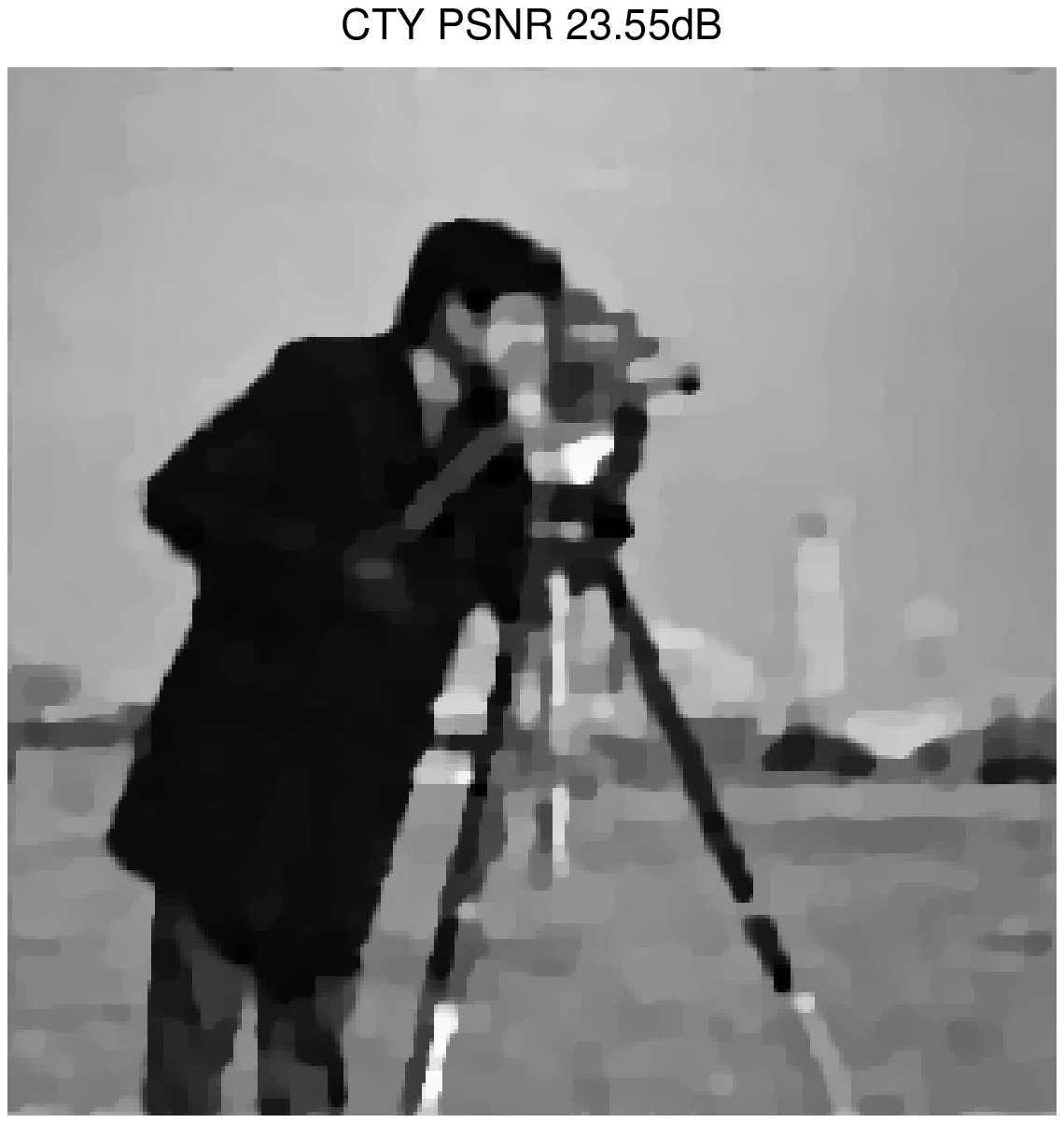}}
    \hspace{0.001in}
   \subfigure{
    \includegraphics[width=0.2\textwidth,clip]{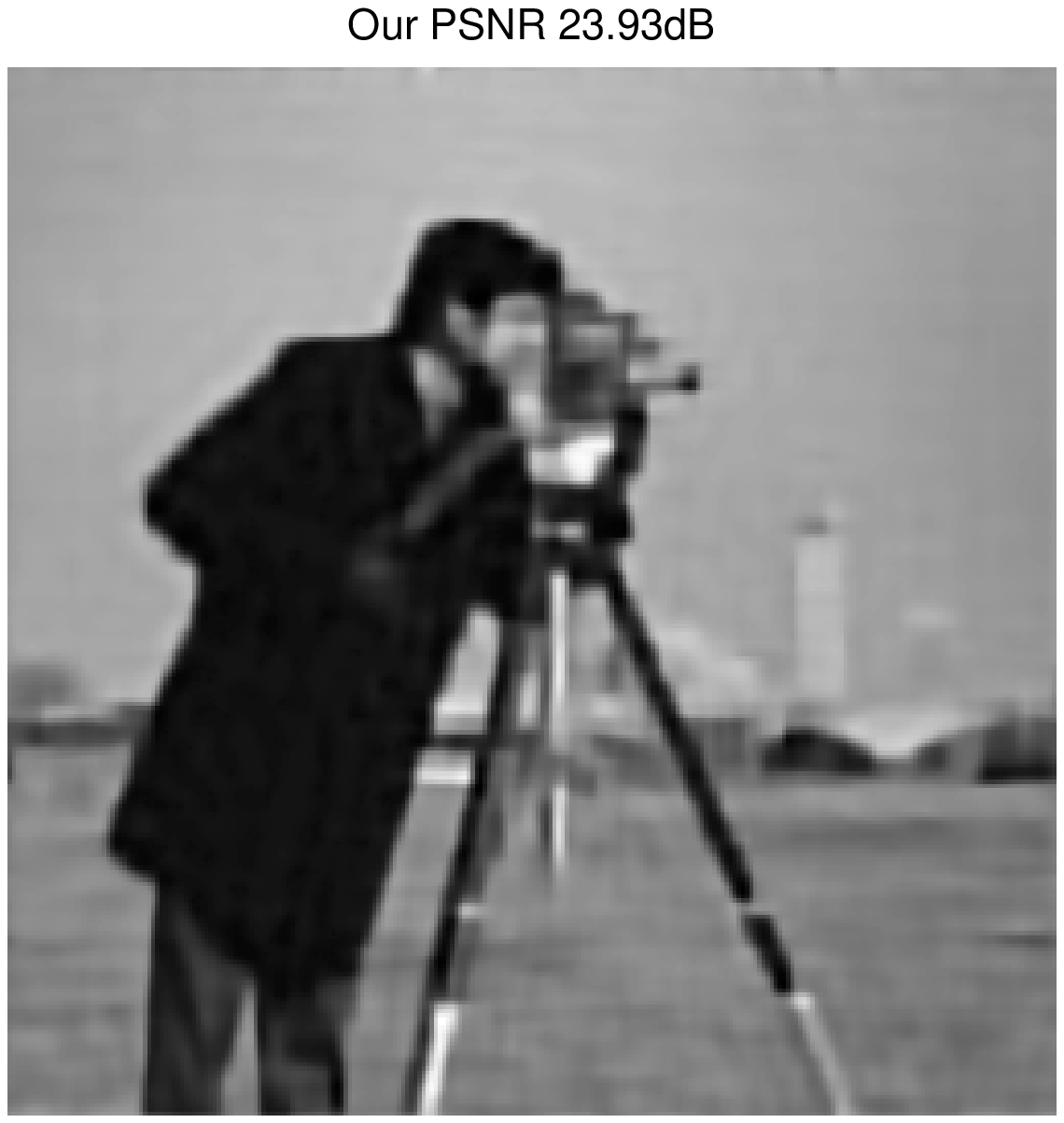}}
         \caption{\textit{Left column: blurred and noisy image. Right columns: restored images by GLN, CTY, and our proposed method respectively}.}
\label{camtoctygln}\end{figure}
\begin{figure}
  \centering
  \subfigure{
    \includegraphics[width=0.4\textwidth,clip]{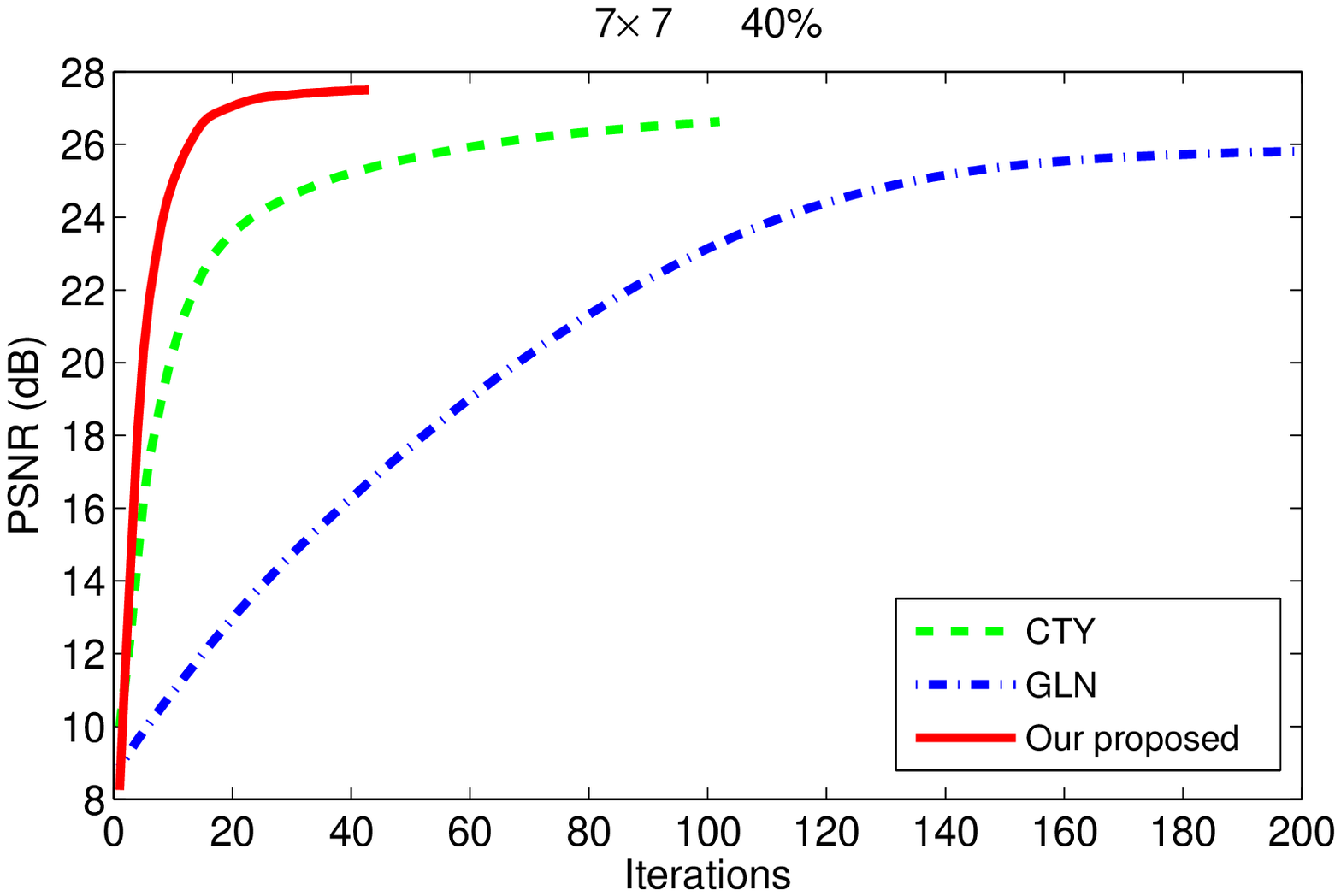}}
    \hspace{0.001in}
   \subfigure{
    \includegraphics[width=0.4\textwidth,clip]{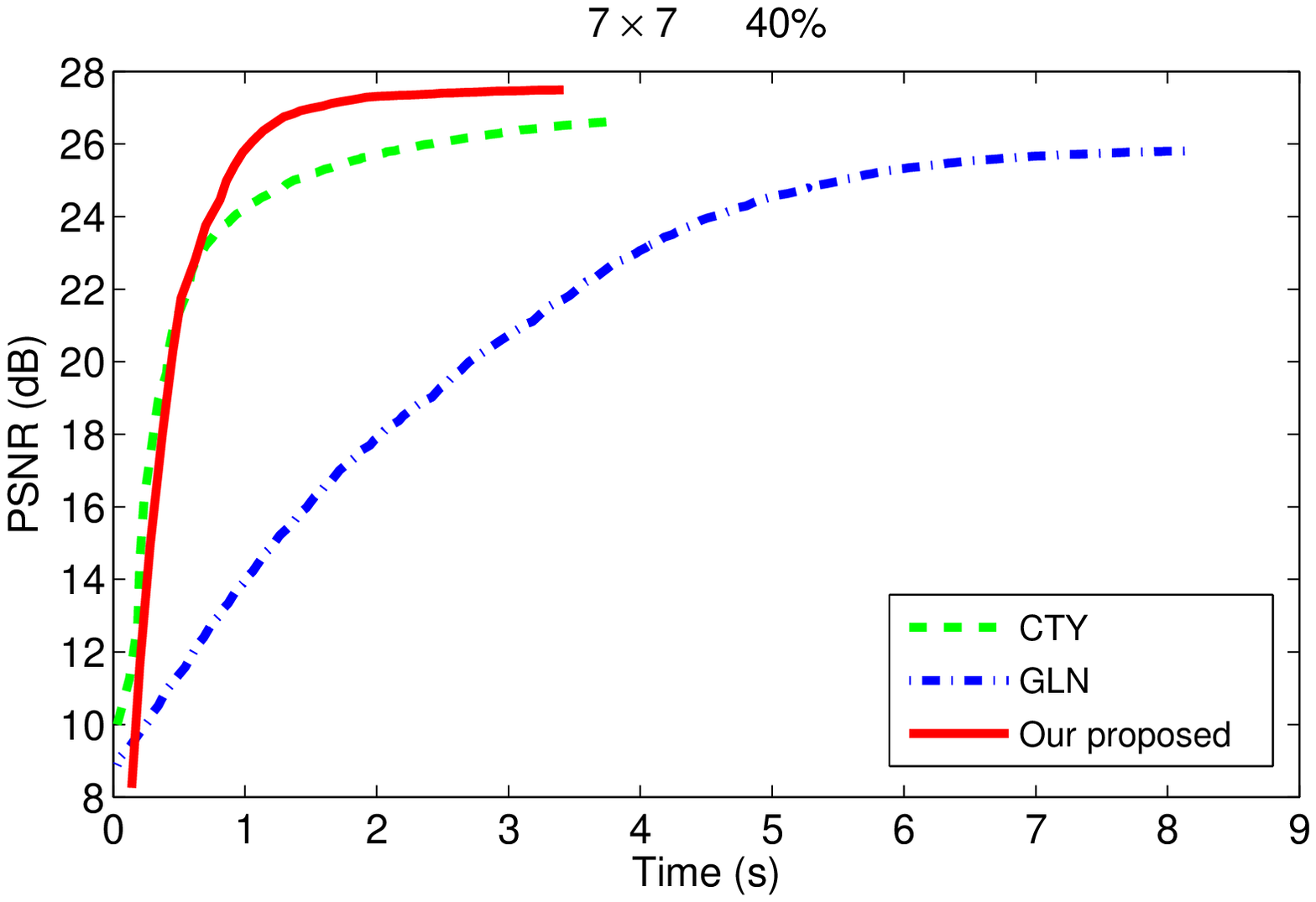}}
    \caption{ \textit{Restoration of the ``Cameraman'' image under $7\times 7$ Gaussion blur with standard deviation 5 and corrupted by 40\% salt-and-pepper noise:
evolution of the PSNR over time and iterations of GLN, CTY, and our proposed method}.}
\label{gau75noi4PSNR2Itrtime}\end{figure}
\begin{table}[!hbp]
\renewcommand{\captionlabeldelim}{.}
\setlength{\abovecaptionskip}{0pt}
\setlength{\belowcaptionskip}{11pt} \centering \caption{\textit{Numerical comparison of CTY and our proposed method under $7\times 7$ Gaussian blur with standard deviation 5 and corrupted by salt-and-pepper noise from 30\% to 60\%}.} \centering
\begin{tabular}{|c|c|r|r|}
\hline 
  \multirow{2}{*}{Images} &\multirow{1}{*}{Noise}                        &  \multicolumn{1}{|c|}{CTY} &  \multicolumn{1}{|c|}{Ours}  \\
 \cline{3-4}
           &                 level                 &\multicolumn{1}{|r|}{$\mu$/Itrs/PSNR/Time/ReE}&\multicolumn{1}{|r|}{$\mu$/Itrs/PSNR/Time/ReE}  \\
 \hline
  \multirow{2}{*}{(b) } &30\% &25/75/29.26/3.01/0.1653&100/31/\text{29.78}/2.65/0.1558\\
        &40\%&24/64/28.59/2.36/0.1787&80/36/\text{29.18}/2.76/0.1669\\
   \multirow{2}{*}{ Satellite}      &50\%&24/54/27.70/1.98/0.1978&60/40/\text{28.28}/3.32/0.1850\\
        &60\%&11/45/26.80/1.93/0.2196&40/66/\text{27.04}/4.52/0.2136\\
\hline
  \multirow{2}{*}{(c)} &30\%&16/90/33.26/3.45/0.0381&100/37/\text{36.20}/2.96/0.0272 \\
        & 40\%&12/75/32.53/3.23/0.0414&80/42/\text{34.99}/3.45/0.0312\\
   \multirow{2}{*}{ House}      & 50\%&13/64/31.56/2.69/0.0463&60/49/\text{33.11}/3.87/0.0288\\
        & 60\%&9/61/30.16/2.58/0.0545&40/62/\text{30.79}/4.27/0.0506\\
        \hline
  \multirow{2}{*}{(d)}&30\%&25/129/28.30/5.65/0.0720&100/35/\text{29.95}/2.70/0.0595\\
        &40\%&24/111/27.24/4.27/0.0813&80/36/\text{28.77}/2.71/0.0682\\
    \multirow{2}{*}{ Boat }     &50\%&18/82/26.13/3.21/0.0924&60/42/\text{27.23}/3.20/0.0814\\
        &60\%&10/61/24.87/2.54/0.1069&40/58/\text{25.57}/4.46/0.0985\\
        \hline
  \multirow{2}{*}{(e)} &30\%&35/125/25.57/5.21/0.0985&100/33/\text{27.00}/3.24/0.0835\\
        &40\%&30/87/24.77/4.77/0.1080&80/41/\text{25.92}/3.95/0.0946\\
    \multirow{2}{*}{ Barbara}     &50\%&21/84/24.14/3.70/0.1162&60/34/\text{24.87}/3.24/0.1067\\
        &60\%&6/58/23.62/2.37/0.1234&40/61/\text{23.94}/5.01/0.1189\\
        \hline
        \multirow{2}{*}{(f)}&30\%&18/106/31.67/4.60/0.0583&100/36/\text{32.87}/3.60/0.0509\\
        &40\%&16/80/30.87/3.42/0.0640&80/36/\text{32.04}/3.28/0.0559\\
    \multirow{2}{*}{ Einstein}     &50\%&16/75/29.88/3.21/0.0718&60/47/\text{30.80}/4.04/0.0645\\
        &60\%&12/64/28.13/2.54/0.0878&40/62/\text{28.59}/5.27/0.0832\\
        \hline
        \multirow{2}{*}{(g)} &30\%&30/131/30.50/5.86/0.0568&100/37/\text{32.54}/3.51/0.0449\\
        &40\%&26/80/28.66/4.41/0.0702&80/36/\text{31.13}/3.54/0.0528\\
     \multirow{2}{*}{ Peppers}    &50\%&23/91/26.66/3.63/0.0883&60/43/\text{28.38}/3.84/0.0725\\
        &60\%&15/66/24.57/2.68/0.1124&40/61/\text{25.83}/5.01/0.0972\\
        \hline
  \multirow{2}{*}{(h)} &30\%&30/142/28.74/17.70/0.0751&100/31/\text{30.69}/11.08/0.0599\\
        &40\%&19/117/27.72/16.44/0.0844&80/37/\text{29.55}/13.28/0.0684\\
    \multirow{2}{*}{ Weather-station}    &50\%&17/84/26.70/13.56/0.0949&60/34/\text{28.00}/15.97/0.0817\\
        &60\%&10/58/25.30/10.65/0.1115&40/61/\text{25.94}/20.06/0.1035\\
\hline
\end{tabular}
\label{allimagecty2ourforgauss}
\end{table}
\begin{table}[!hbp]
\renewcommand{\captionlabeldelim}{.}
\setlength{\abovecaptionskip}{0pt}
\setlength{\belowcaptionskip}{11pt} \centering \caption{\textit{Numerical comparison of CTY and our proposed method under $7\times 7$ average blur and corrupted by salt-and-pepper noise from 30\% to 60\%}.} \centering
\begin{tabular}{|c|c|r|r|}
\hline
\multirow{2}{*}{Images} &Noise               &  \multicolumn{1}{|c|}{CTY} &  \multicolumn{1}{|c|}{Ours}   \\
 \cline{3-4}
           &                 level                      &\multicolumn{1}{|r|}{$\mu$/Itrs/PSNR/Time/ReE}&\multicolumn{1}{|r|}{$\mu$/Itrs/PSNR/Time/ReE}    \\
 \hline
\multirow{2}{*}{(b)} &30\% &21/72/29.45/2.71/0.1617&100/34/30.05/2.84/0.1510\\
        &40\%&18/57/28.66/2.23/0.1772&80/40/29.23/3.14/0.1660\\
     \multirow{2}{*}{ Satellite}    &50\%&14/48/27.69/1.92/0.1981&60/49/28.24/3.71/0.1860\\
        &60\%&9/41/26.74/1.84/0.2209&40/72/27.08/5.35/0.2126\\
\hline
  \multirow{2}{*}{(c)} &30\%&20/105/33.61/3.90/0.0366&100/38/36.40/2.90/0.0265 \\
        & 40\%&12/75/32.63/2.93/0.0410&80/43/35.02/3.48/0.0311 \\
     \multirow{2}{*}{ House}    & 50\%&10/65/31.69/2.66/0.0456&60/50/32.89/3.87/0.0397\\
        & 60\%&7/60/30.18/2.51/0.0543&40/59/30.65/4.46/0.0515\\
        \hline
  \multirow{2}{*}{(d)}&30\%&27/138/28.54/5.43/0.0700&100/36/30.14/3.04/0.0582\\
        &40\%&21/104/27.40/4.23/0.0799&80/37/28.86/2.96/0.0675\\
      \multirow{2}{*}{Boat}   &50\%&17/82/26.19/3.10/0.0917&60/43/27.24/3.51/0.0813\\
        &60\%&13/65/24.87/2.50/0.1068&40/59/25.52/4.52/0.0991\\
        \hline
  \multirow{2}{*}{(e)} &30\% &31/142/25.65/5.44/0.0976&100/32/27.16/2.56/0.0820\\
        &40\%&28/113/24.81/5.09/0.1075&80/39/26.02/3.37/0.0935 \\
     \multirow{2}{*}{Barbara}    &50\%&19/84/24.17/3.31/0.1157&60/44/24.94/3.53/0.1060\\
        &60\%&6/57/23.60/2.43/0.1237&40/61/23.97/4.85/0.1184\\
        \hline
        \multirow{2}{*}{(f)} &30\%&22/117/31.80/4.54/0.0575&100/36/33.01/2.96/0.0500  \\
        &40\%&15/88/31.11/3.45/0.0623&80/38/32.15/3.14/0.0552\\
      \multirow{2}{*}{Einstein}   &50\%&14/74/30.00/2.98/0.0708&60/48/30.79/3.93/0.0646 \\
        &60\%&12/64/28.36/2.40/0.0855&40/61/28.72/4.65/0.0820\\
        \hline
        \multirow{2}{*}{(g)} &30\%&28/134/30.86/5.46/0.0545&100/38/32.71/3.03/0.0440\\
        &40\%&26/109/29.01/4.27/0.0674&80/39/31.30/3.21/0.0518\\
      \multirow{2}{*}{ Peppers}   &50\%&22/90/26.84/3.45/0.0865&60/44/28.63/3.71/0.0704\\
        &60\%&14/70/24.78/2.71/0.1097&40/60/26.04/4.66/0.0948\\
        \hline
  \multirow{2}{*}{(h)} &30\%&25/117/28.87/21.53/0.0739&100/35/30.86/11.65/0.0588\\
        &40\%&23/98/27.84/17.94/0.0834&80/42/29.57/13.63/0.0682\\
      \multirow{2}{*}{Weather-station}   &50\%&15/73/26.70/13.32/0.0949&60/51/27.93/16.43/0.0824\\
        &60\%&10/60/25.26/10.55/0.1120&40/65/25.93/20.31/0.1037\\
\hline
\end{tabular}
\label{allimagecty2ourforaverage}
\end{table}
\begin{figure}
  \centering
  \subfigure{
    \includegraphics[width=0.2\textwidth,clip]{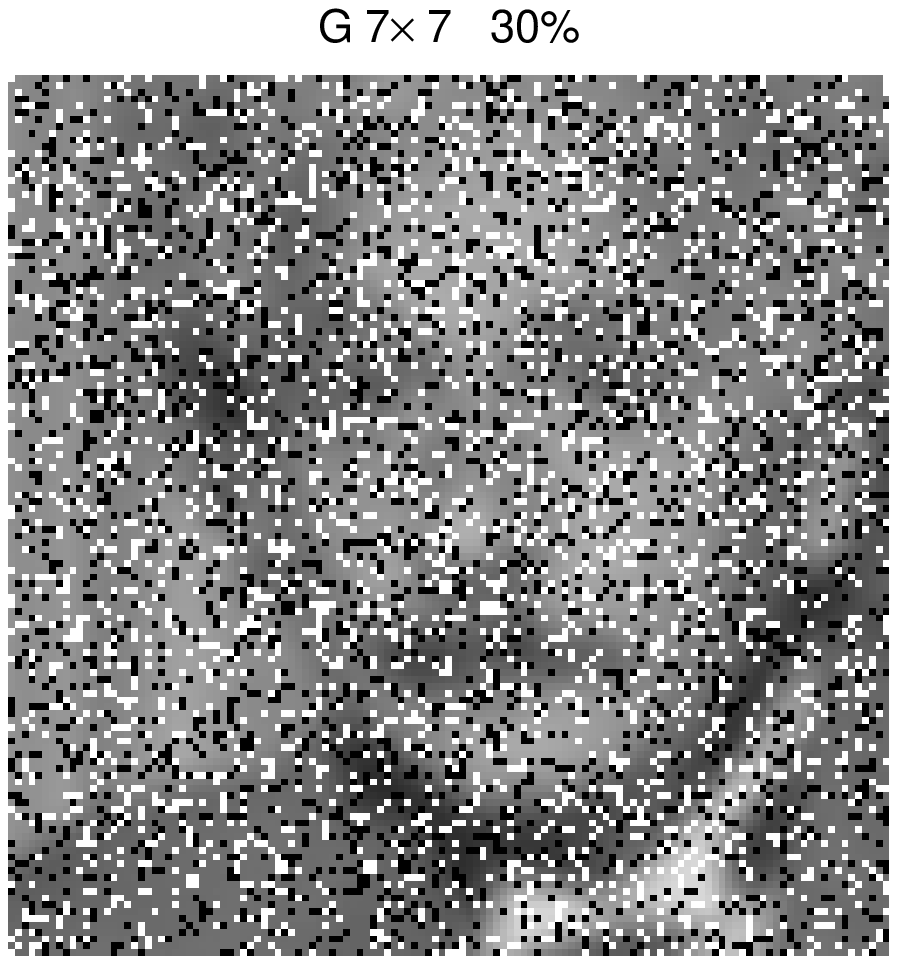}}
    \hspace{0.001in}
   \subfigure{
    \includegraphics[width=0.2\textwidth,clip]{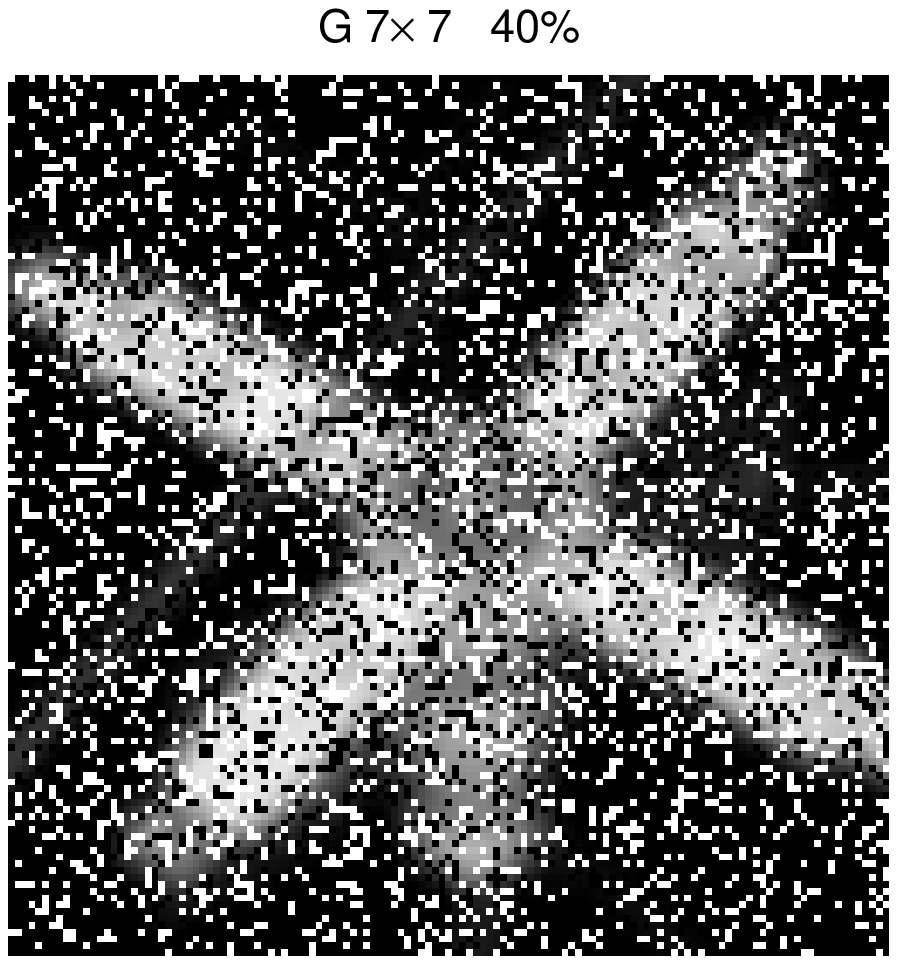}}
    \hspace{0.001in}
    \subfigure{
    \includegraphics[width=0.2\textwidth,clip]{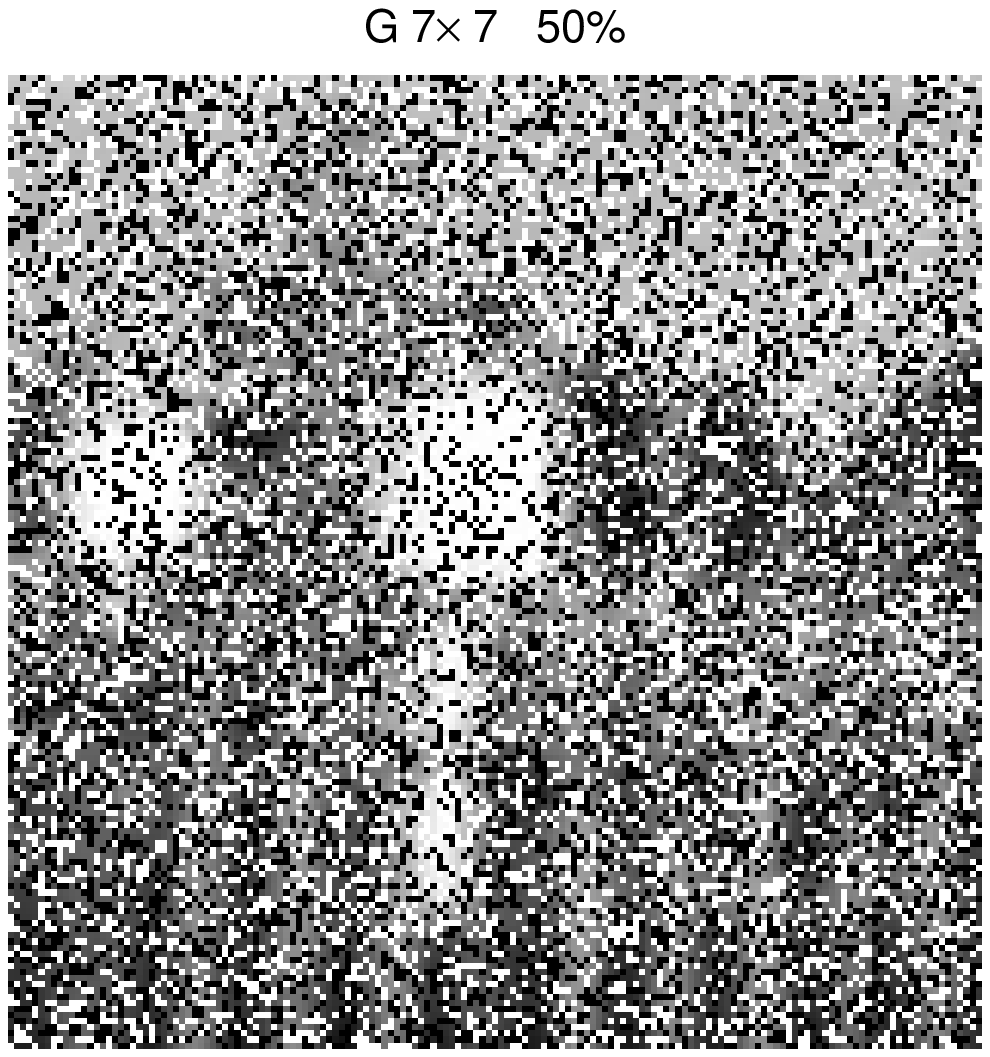}}
    \hspace{0.001in}
    \subfigure{
    \includegraphics[width=0.2\textwidth,clip]{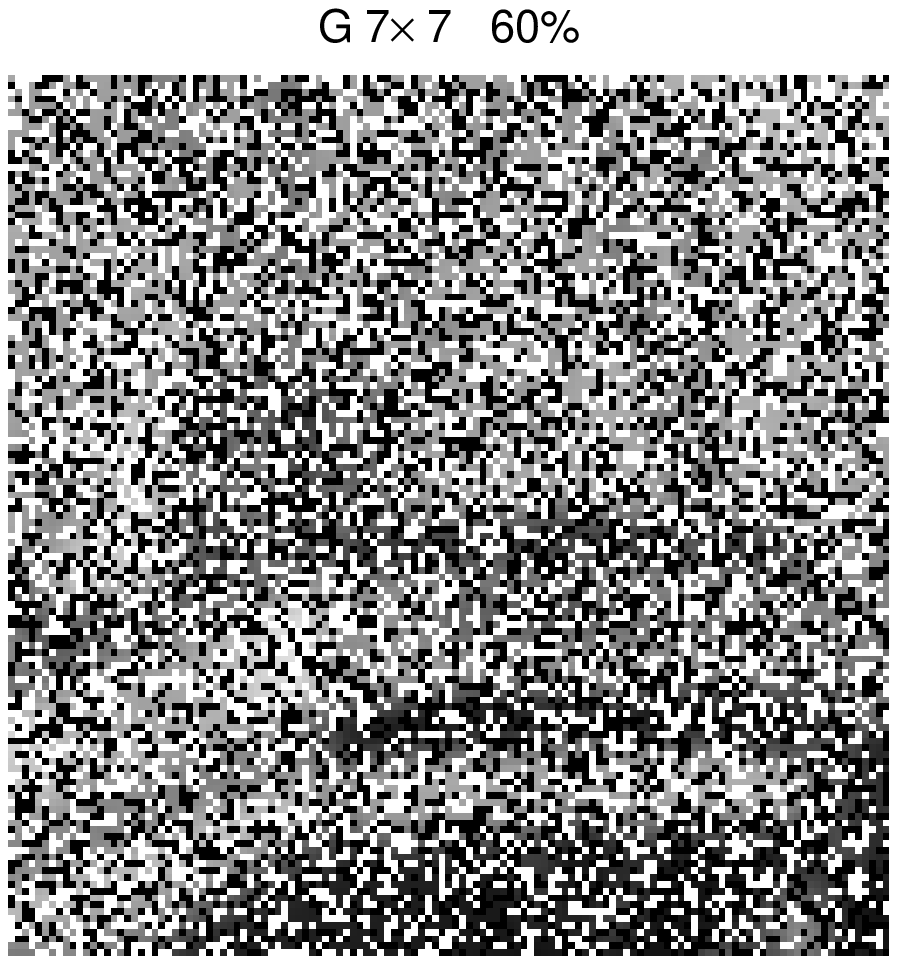}}
    \hspace{0.001in}
   \subfigure{
    \includegraphics[width=0.2\textwidth,clip]{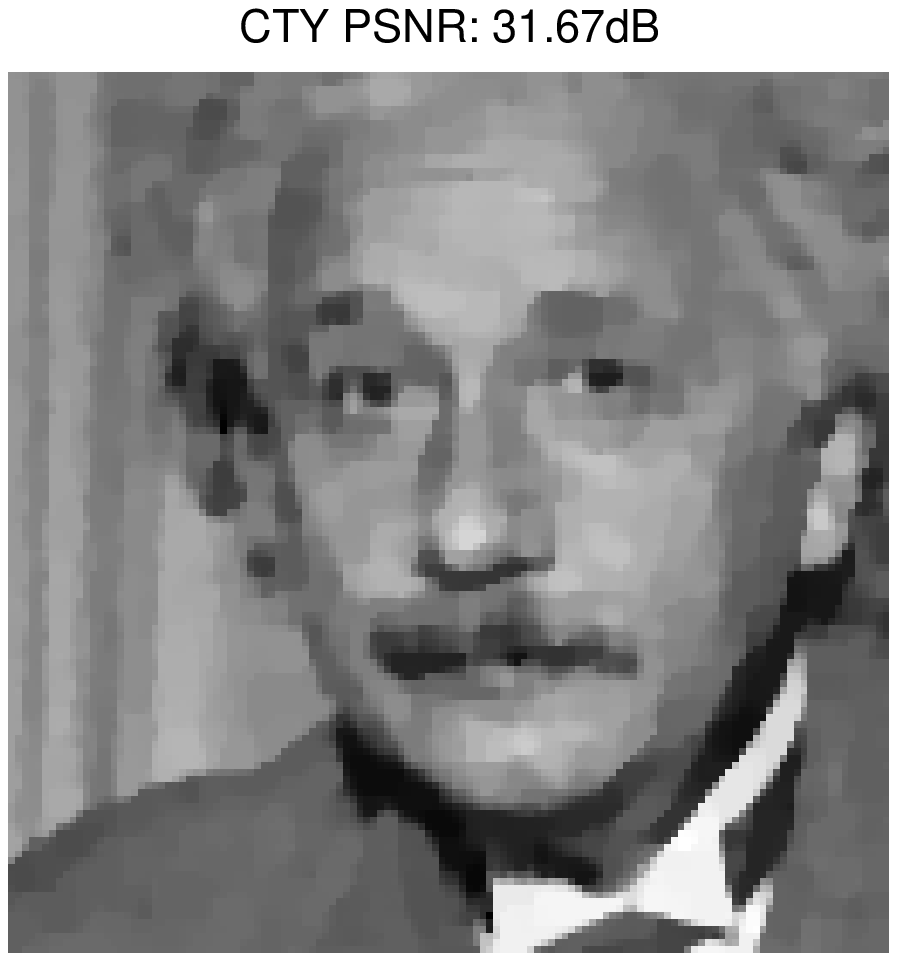}}
    \hspace{0.001in}
    \subfigure{
    \includegraphics[width=0.2\textwidth,clip]{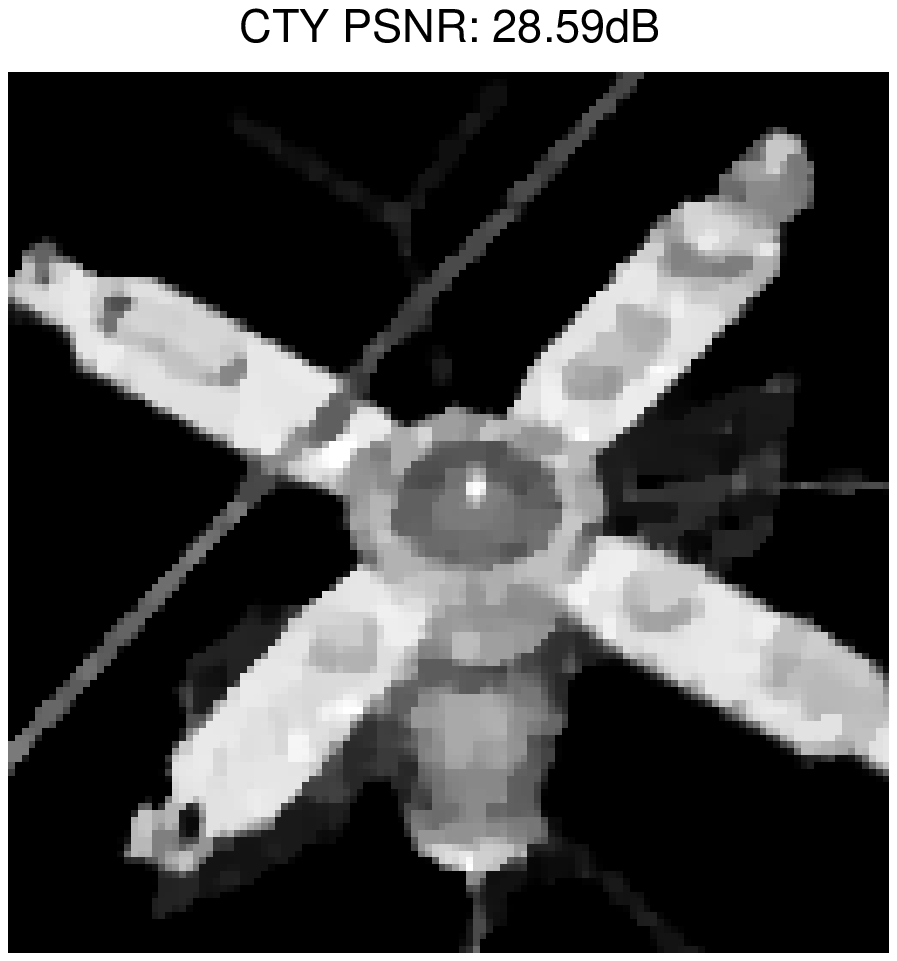}}
    \hspace{0.001in}
    \subfigure{
    \includegraphics[width=0.2\textwidth,clip]{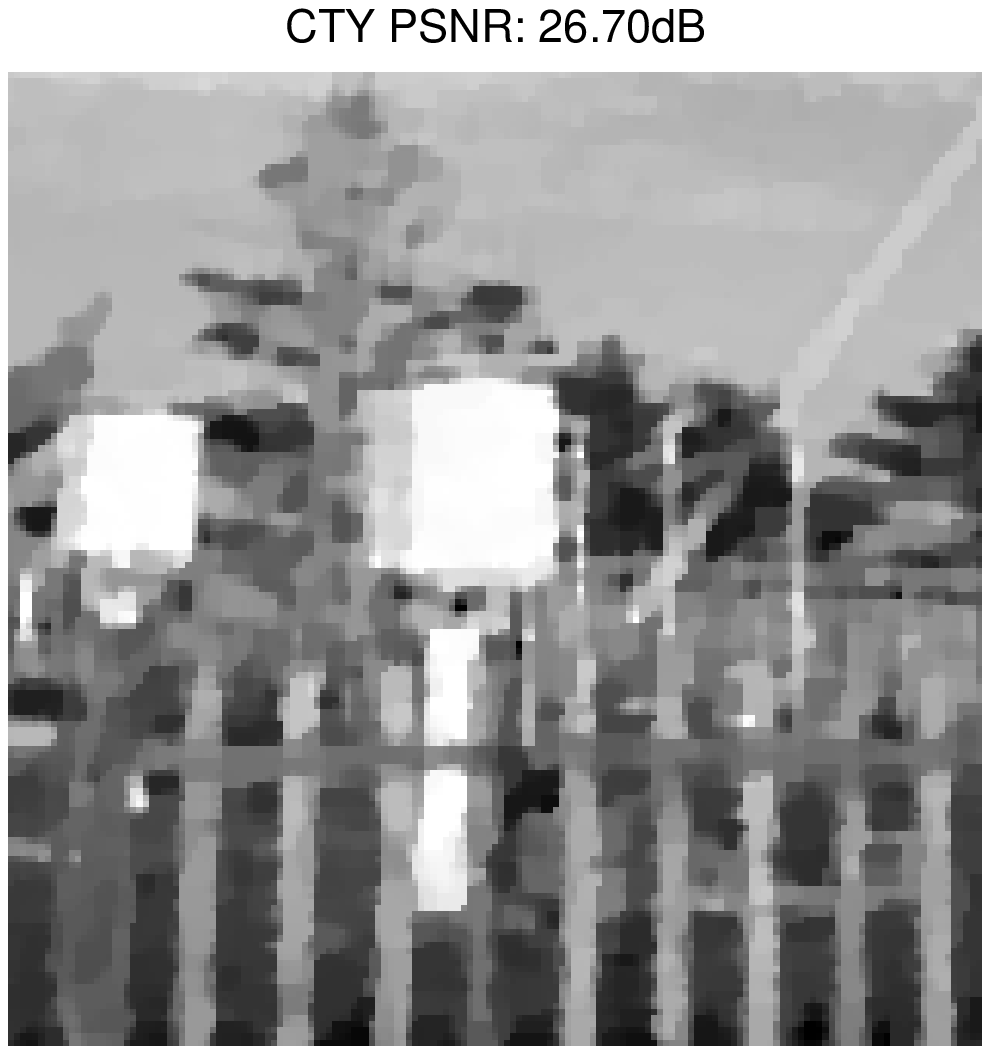}}
    \hspace{0.001in}
    \subfigure{
    \includegraphics[width=0.2\textwidth,clip]{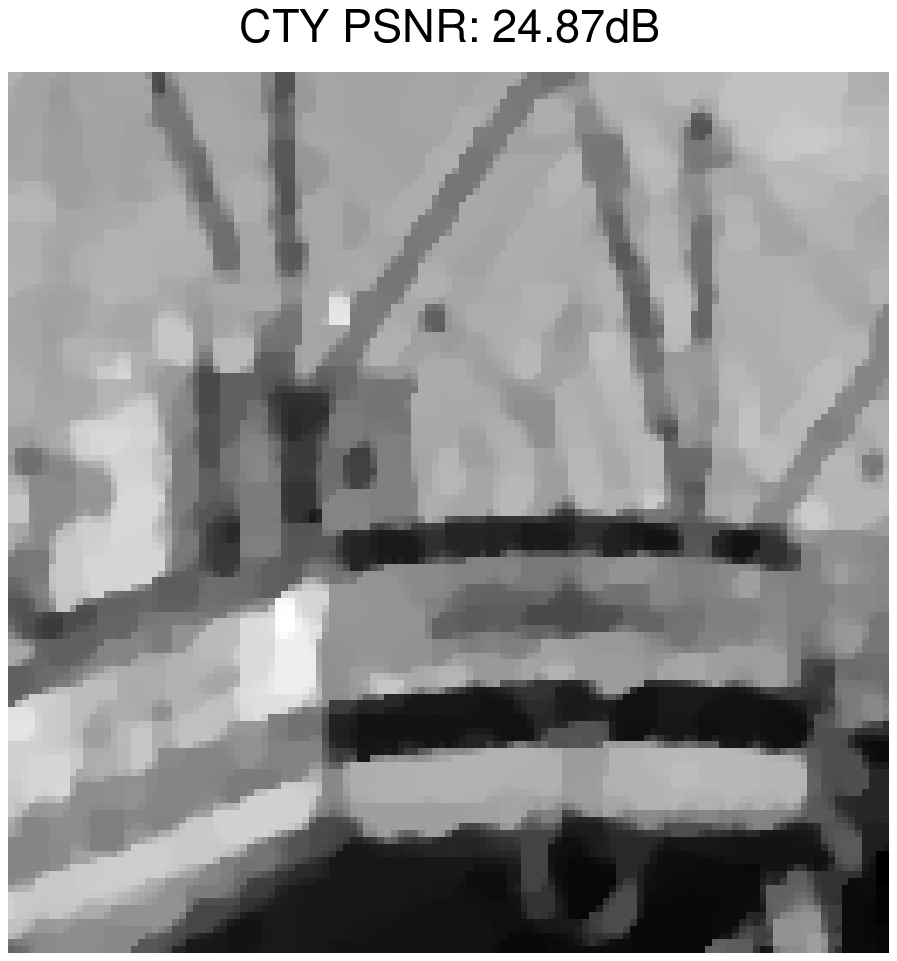}}
    \hspace{0.001in}
   \subfigure{
    \includegraphics[width=0.2\textwidth,clip]{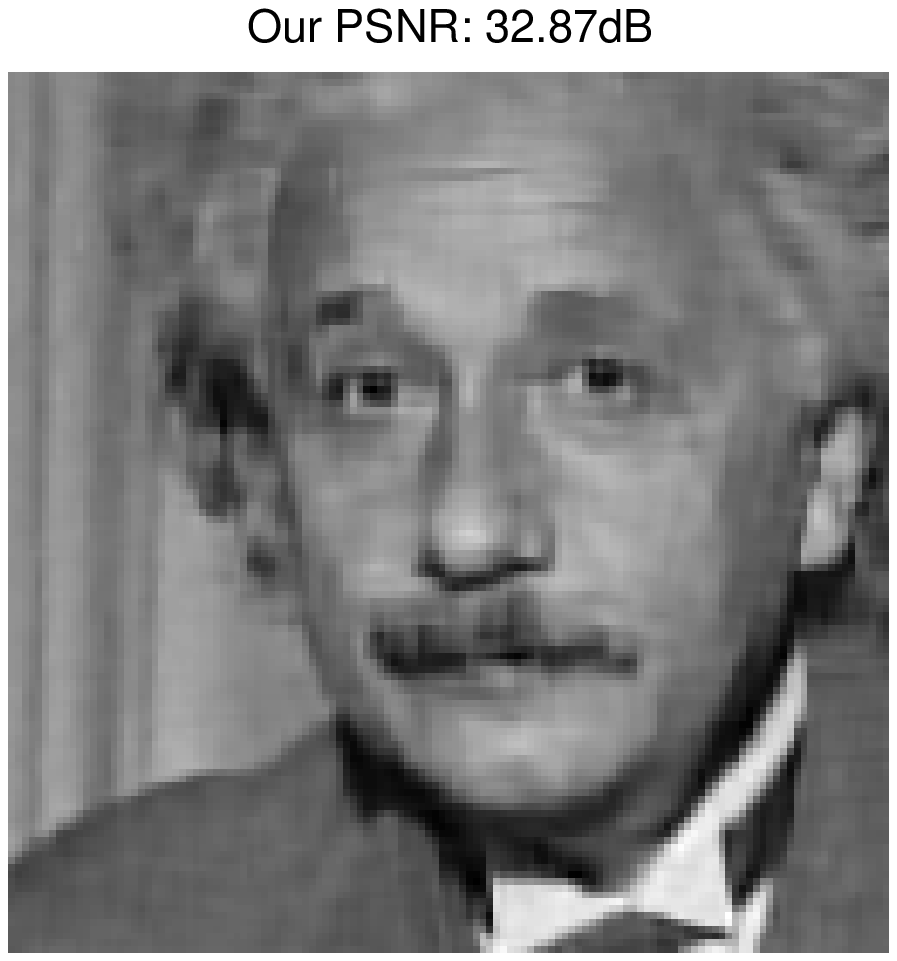}}
    \hspace{0.001in}
    \subfigure{
    \includegraphics[width=0.2\textwidth,clip]{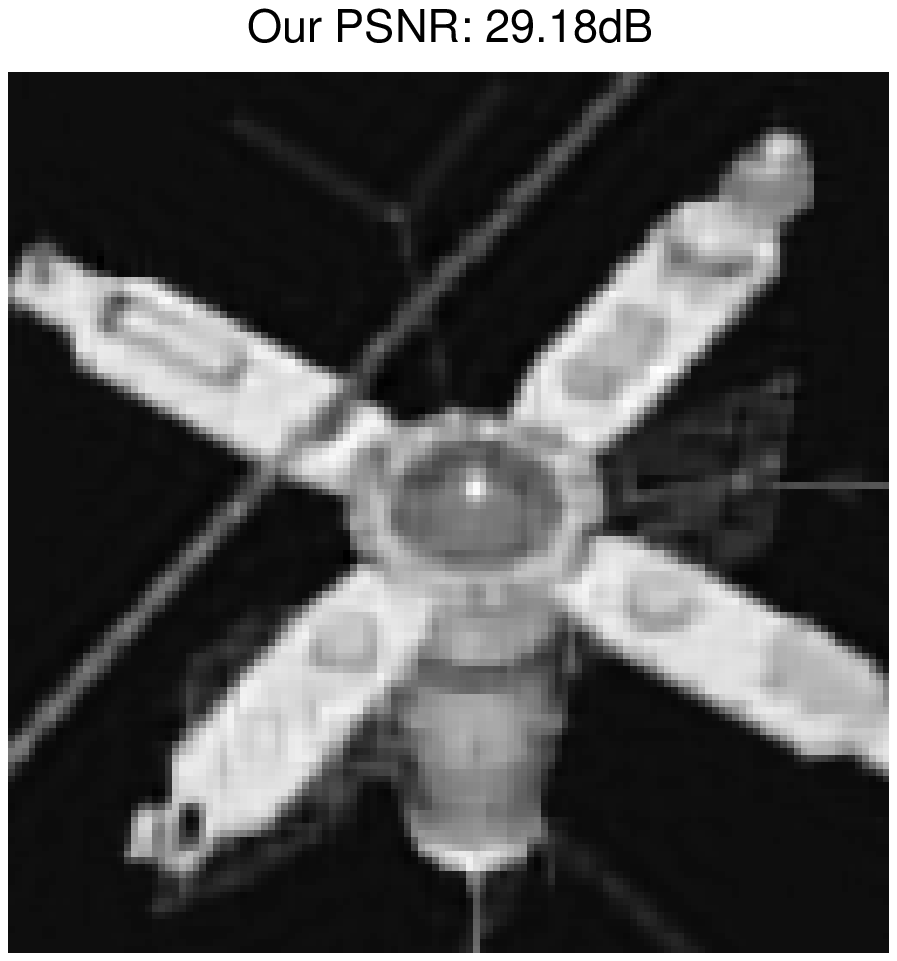}}
    \hspace{0.001in}
    \subfigure{
    \includegraphics[width=0.2\textwidth,clip]{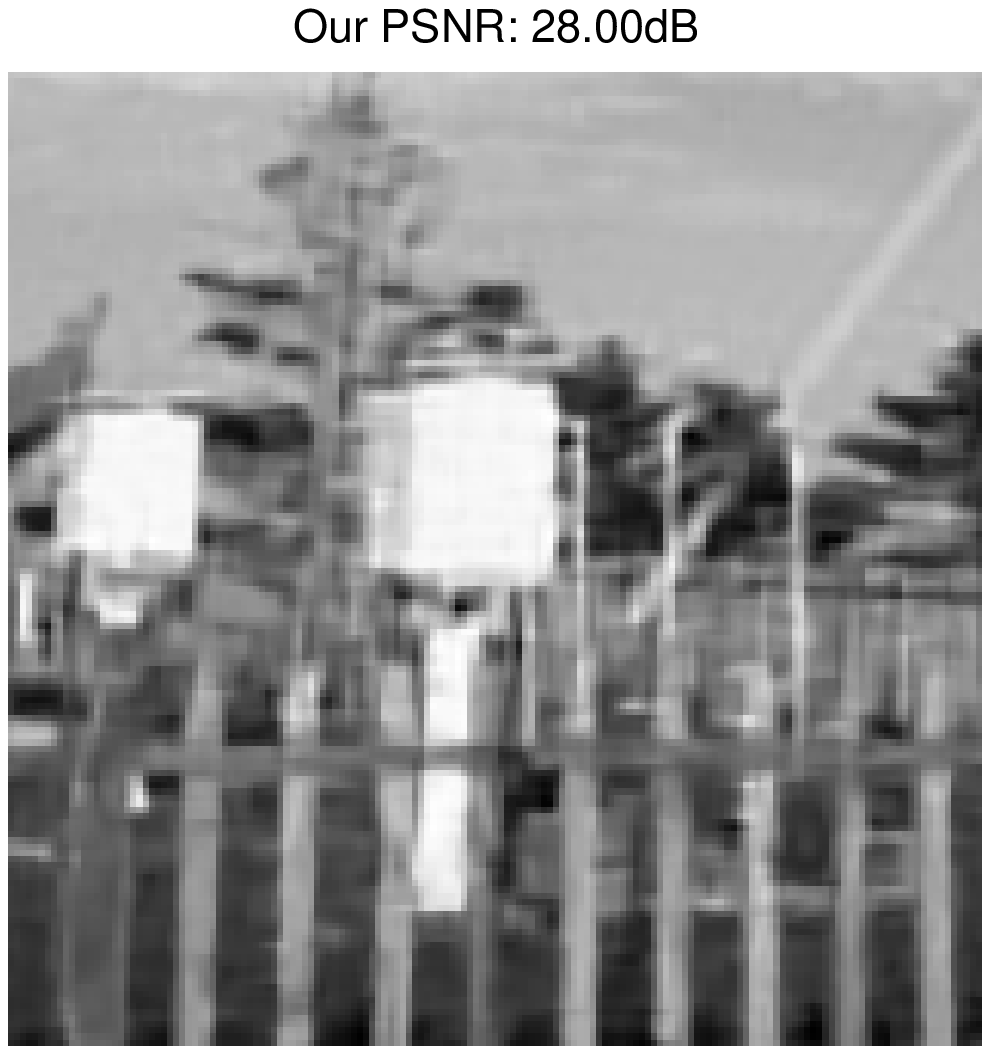}}
    \hspace{0.001in}
   \subfigure{
    \includegraphics[width=0.2\textwidth,clip]{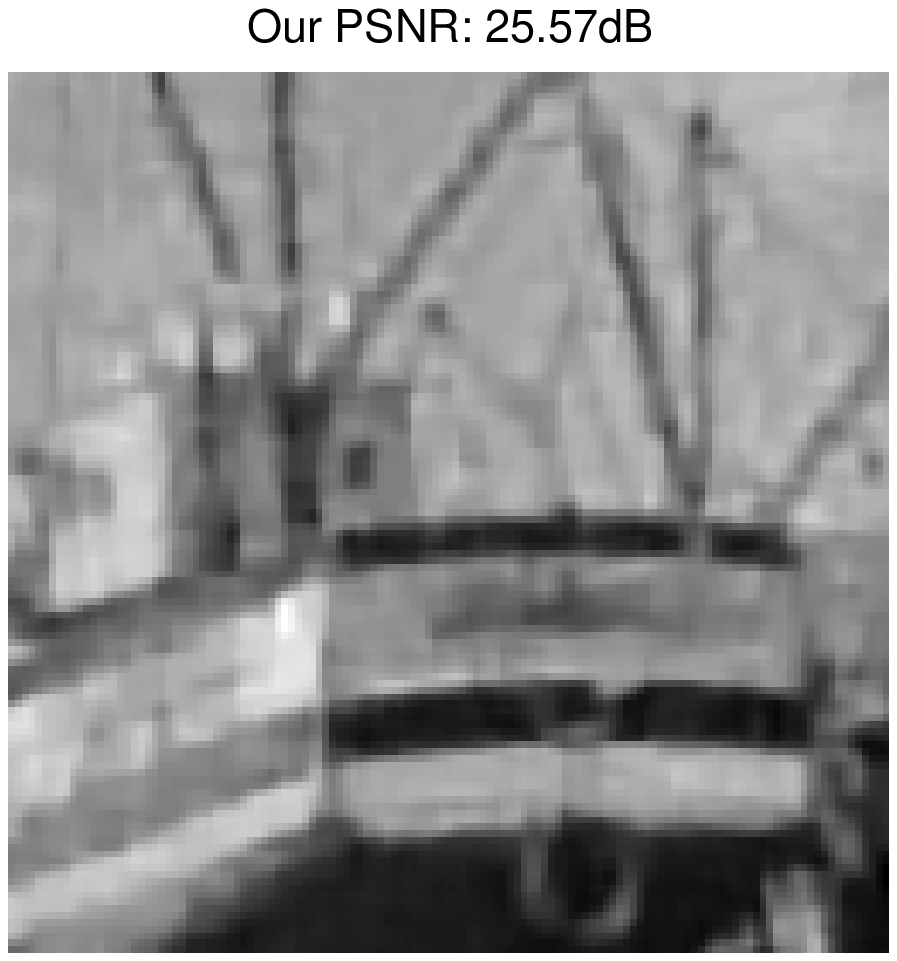}}
         \caption{\textit{Several random examples of degraded and restored images. Top row, zoom parts of blurred and noisy images under $7\times 7$  Gaussian blur with standard deviation 5 and corrupted by salt-and-pepper noise. Second row, zoom parts of restored images by CTY respectively. Third row, zoom parts of restored images by our proposed method respectively}.}
\label{imagestoctyforgauss7}\end{figure}
\begin{figure}
  \centering
  \subfigure{
    \includegraphics[width=0.2\textwidth,clip]{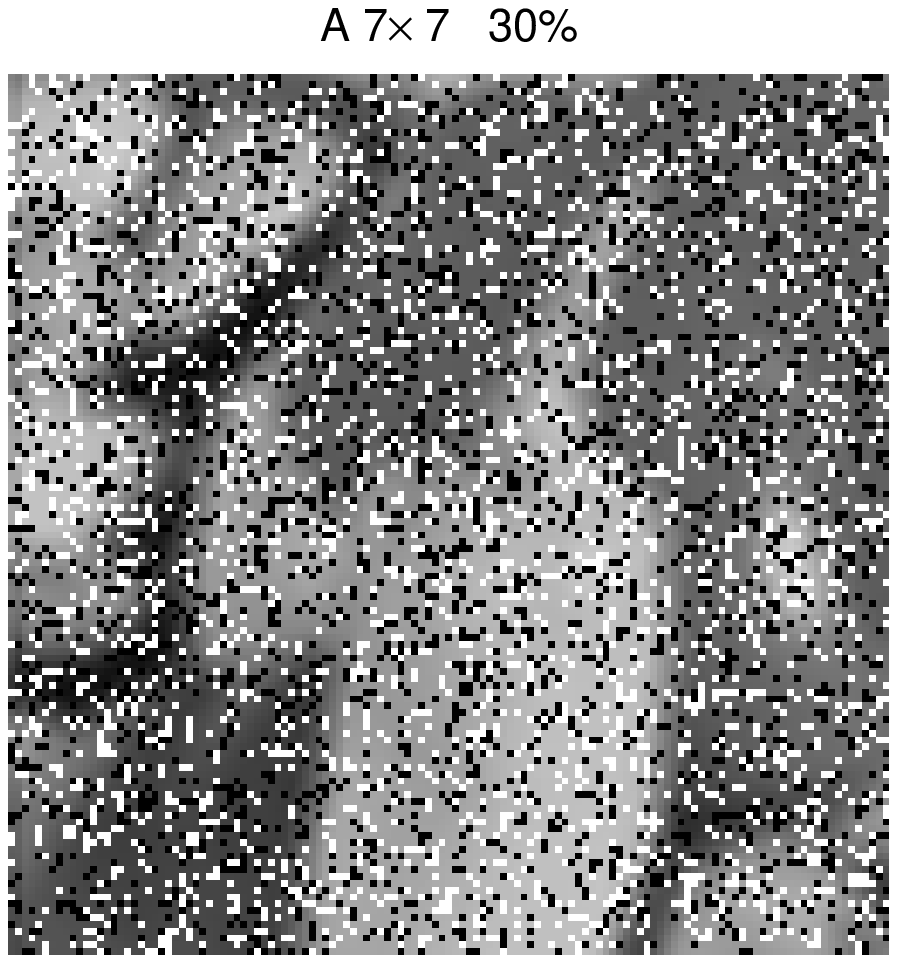}}
    \hspace{0.001in}
   \subfigure{
    \includegraphics[width=0.2\textwidth,clip]{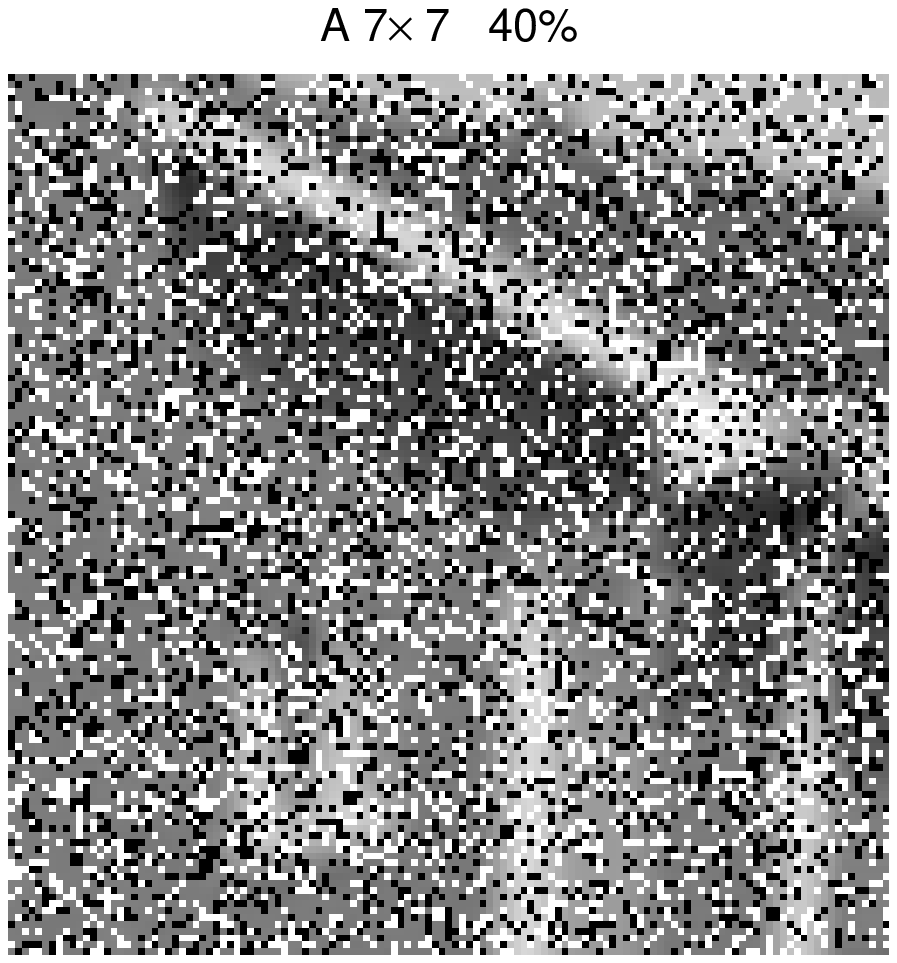}}
    \hspace{0.001in}
    \subfigure{
    \includegraphics[width=0.2\textwidth,clip]{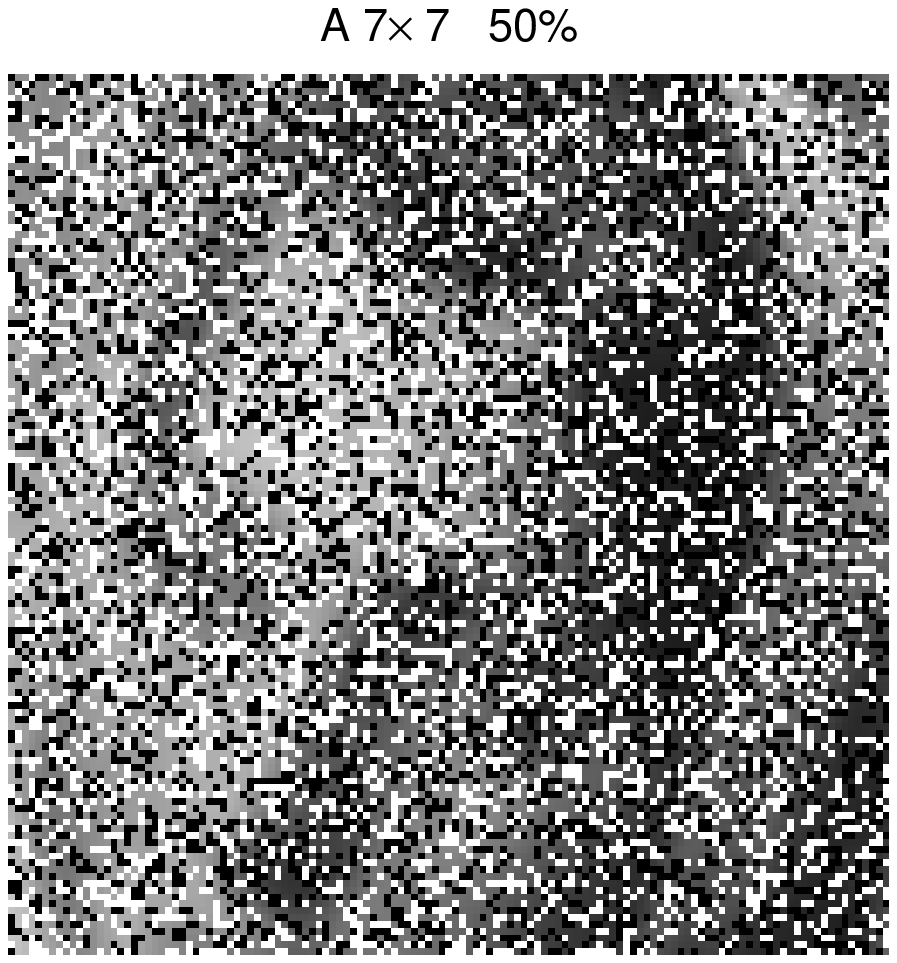}}
    \hspace{0.001in}
    \subfigure{
    \includegraphics[width=0.2\textwidth,clip]{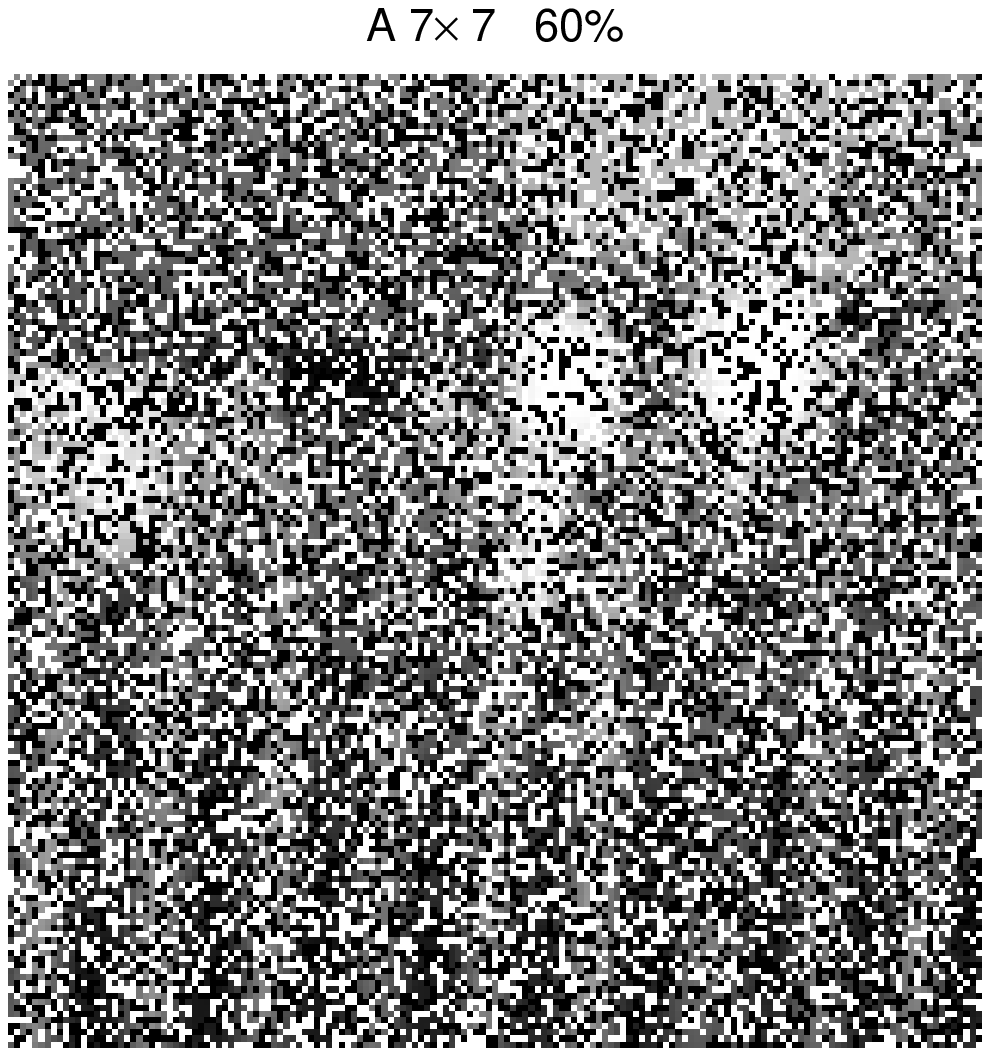}}
    \hspace{0.001in}
   \subfigure{
    \includegraphics[width=0.2\textwidth,clip]{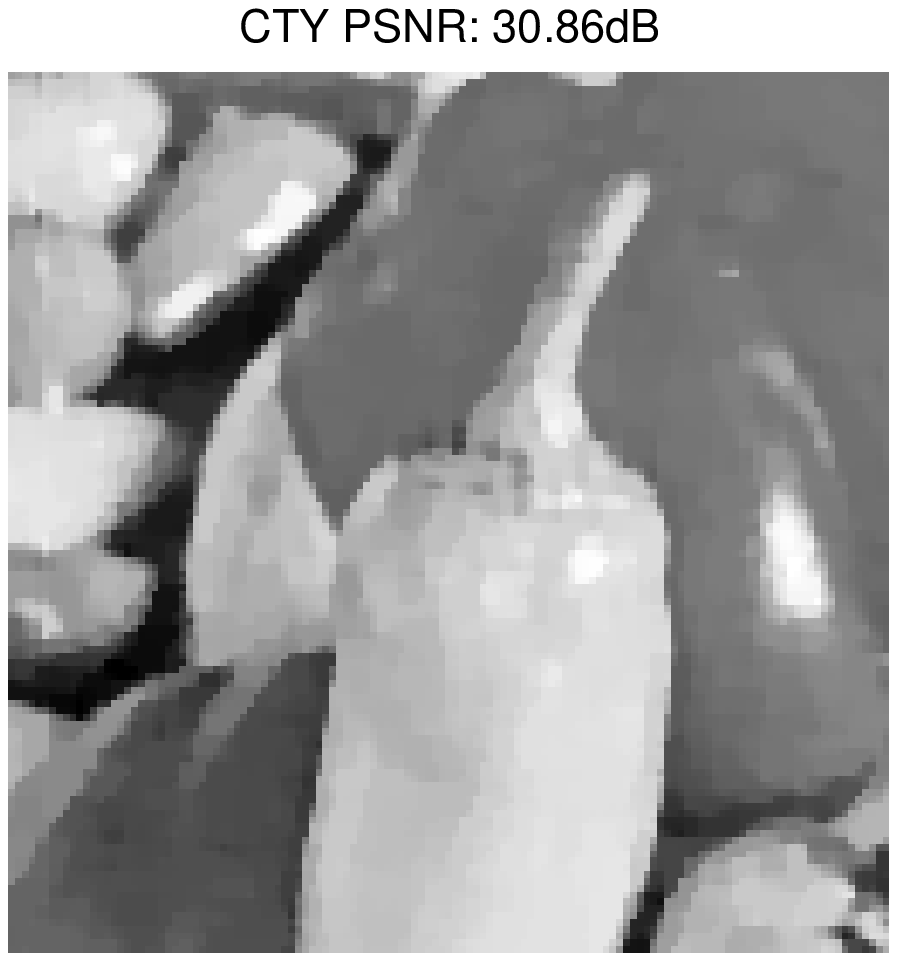}}
    \hspace{0.001in}
    \subfigure{
    \includegraphics[width=0.2\textwidth,clip]{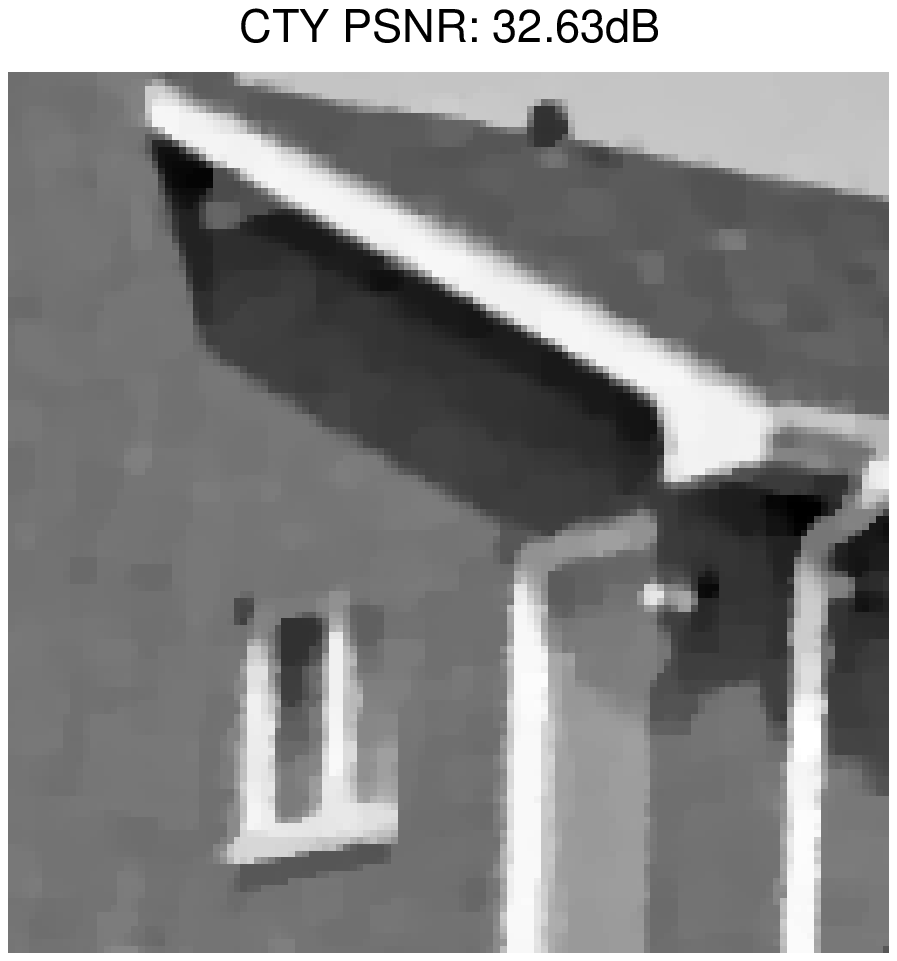}}
    \hspace{0.001in}
    \subfigure{
    \includegraphics[width=0.2\textwidth,clip]{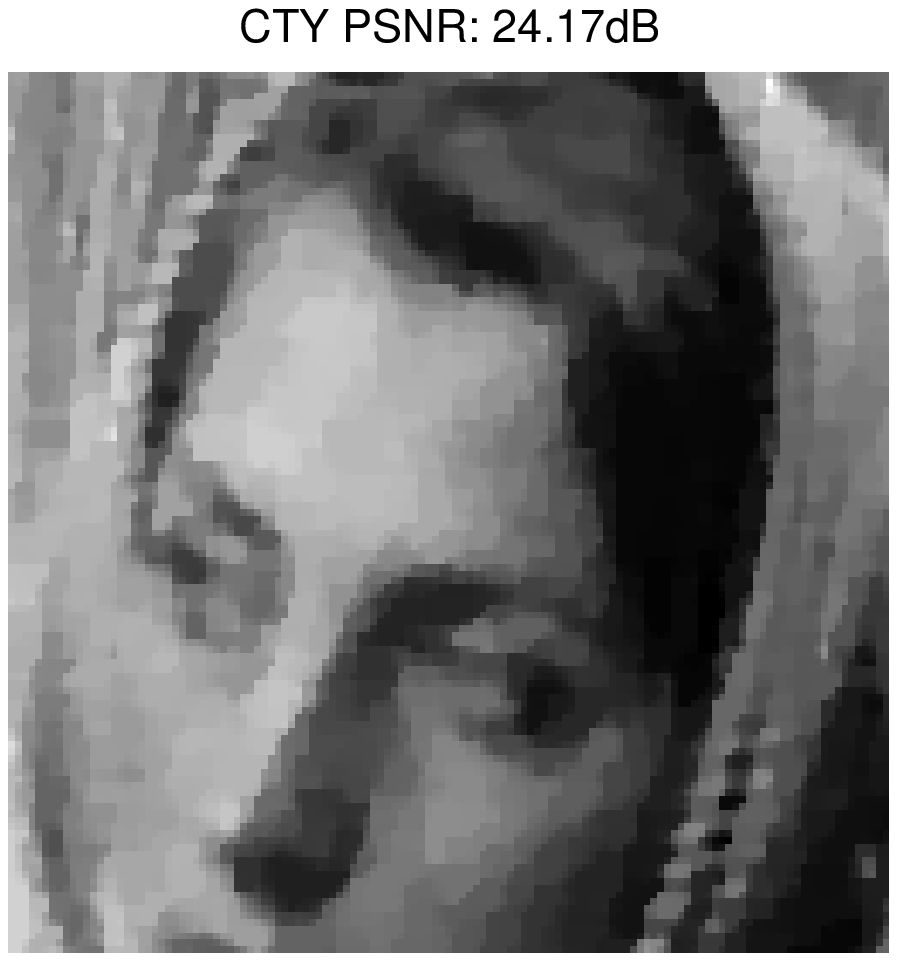}}
    \hspace{0.001in}
    \subfigure{
    \includegraphics[width=0.2\textwidth,clip]{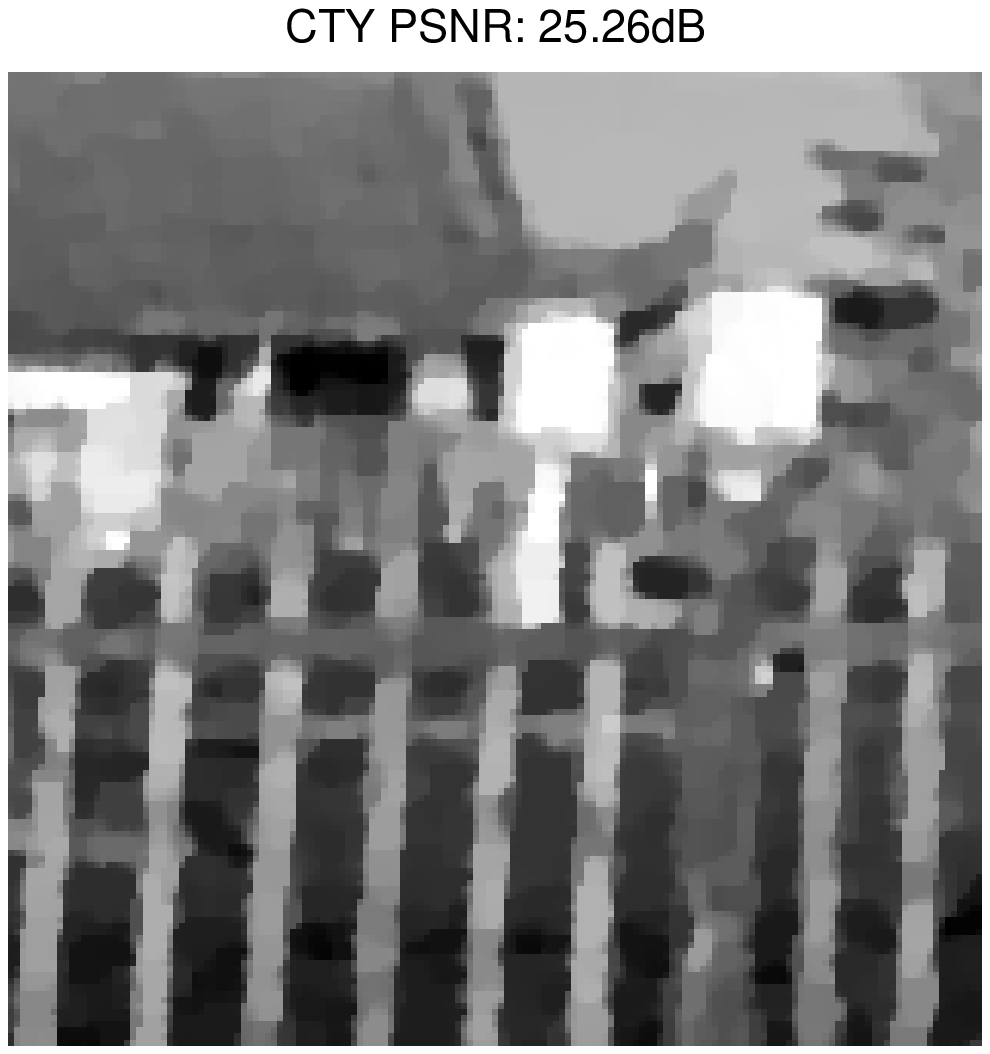}}
    \hspace{0.001in}
   \subfigure{
    \includegraphics[width=0.2\textwidth,clip]{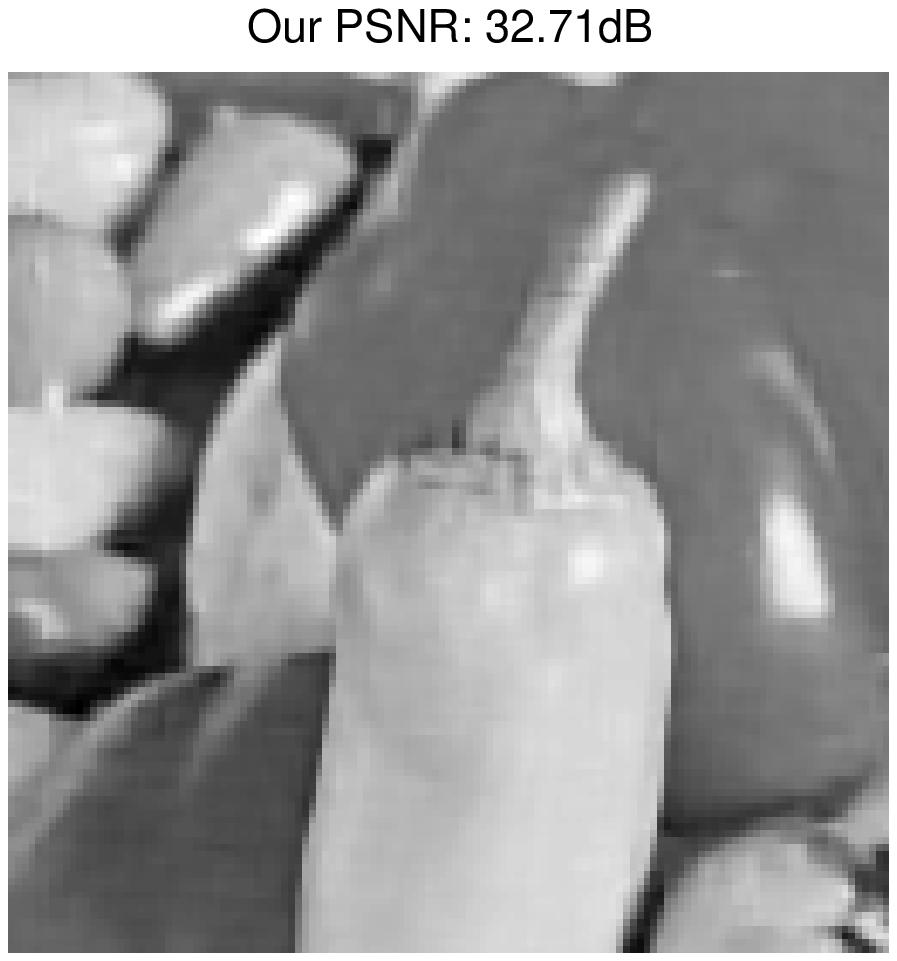}}
    \hspace{0.001in}
    \subfigure{
    \includegraphics[width=0.2\textwidth,clip]{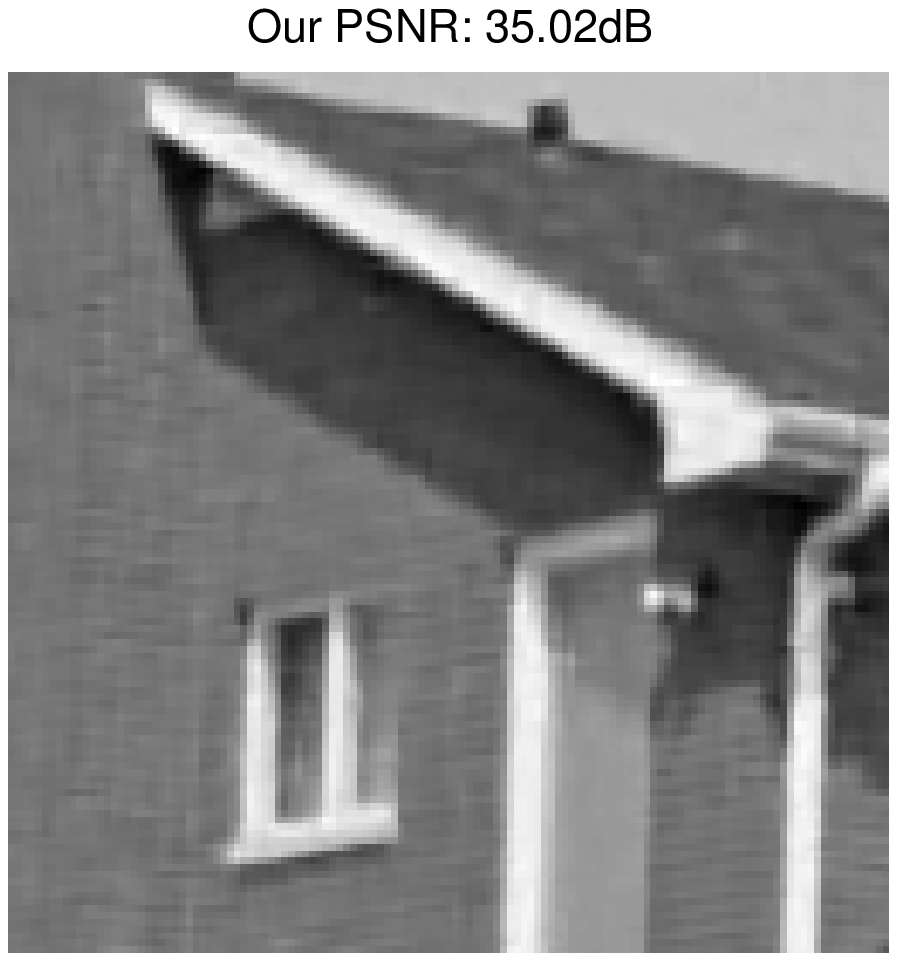}}
    \hspace{0.001in}
    \subfigure{
    \includegraphics[width=0.2\textwidth,clip]{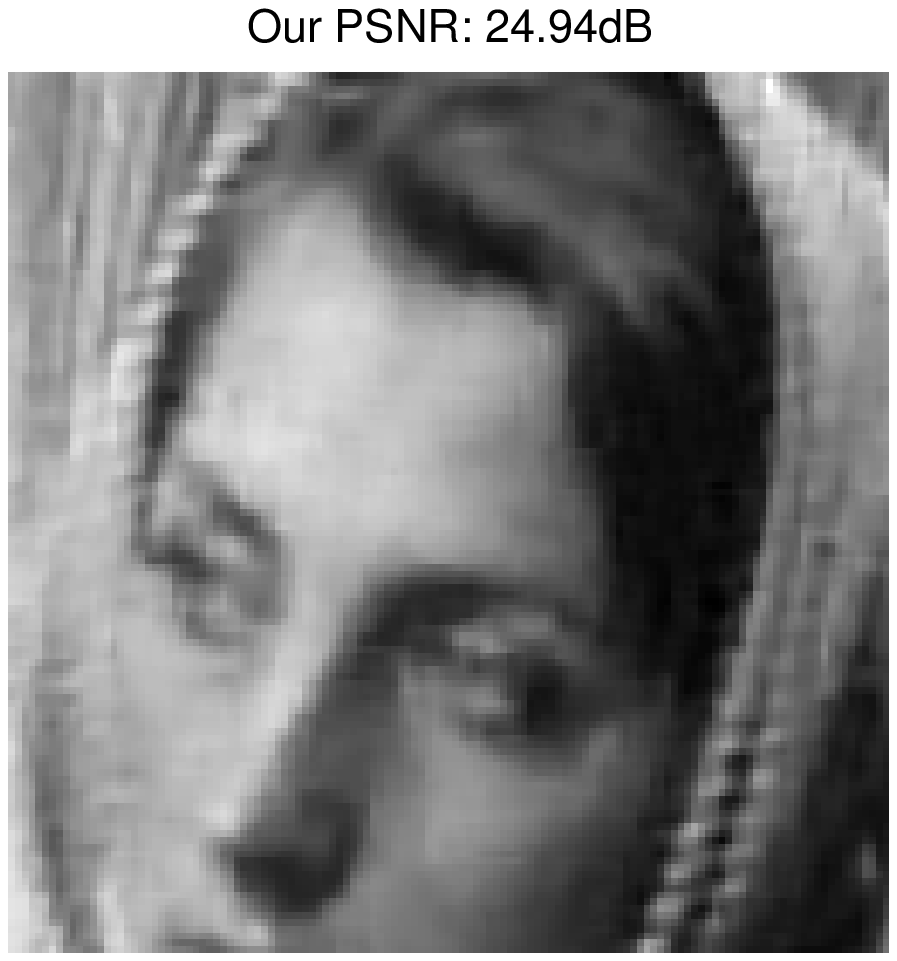}}
    \hspace{0.001in}
   \subfigure{
    \includegraphics[width=0.2\textwidth,clip]{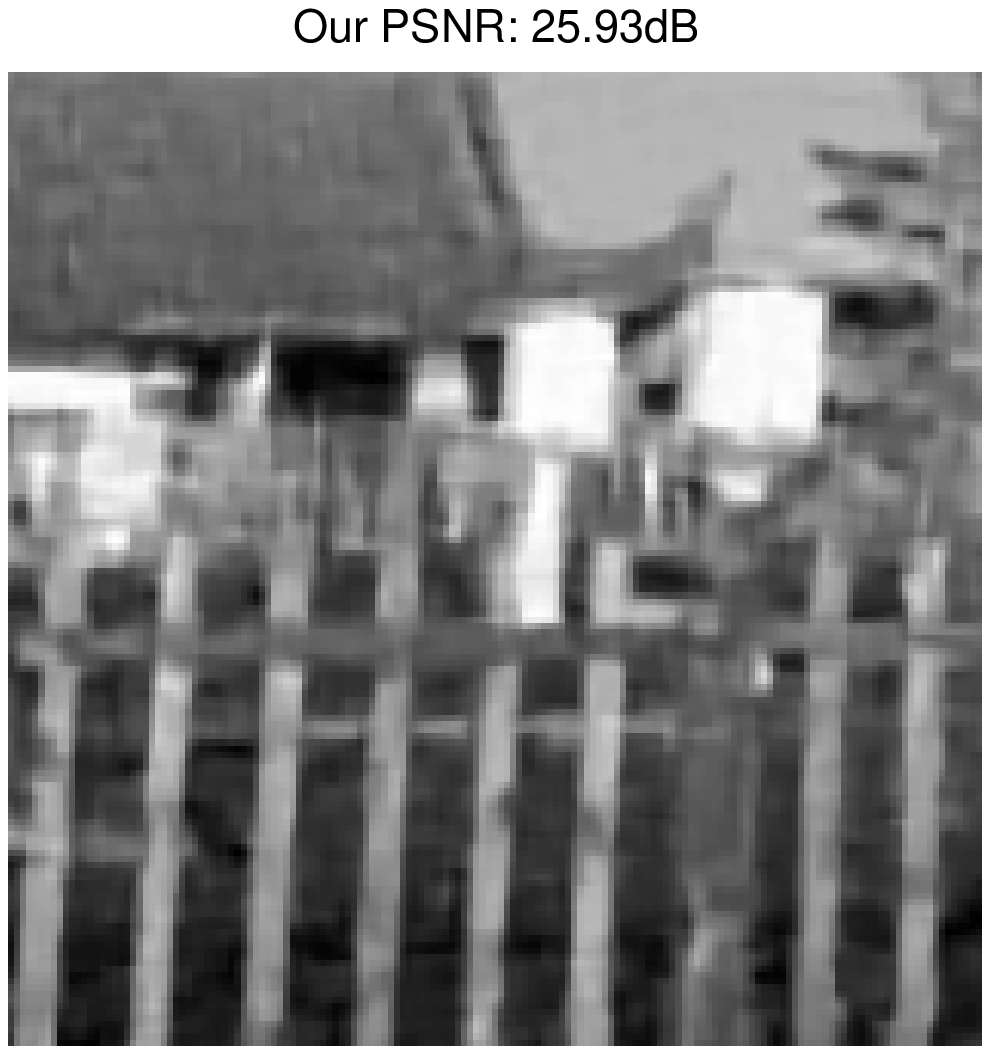}}
         \caption{\textit{Several random examples of degraded and restored images. Top row, zoom parts of blurred and noisy images under $7\times 7$ average blur and corrupted by salt-and-pepper noise. Second row, zoom parts of restored images by CTY respectively. Third row, zoom parts of restored images by our proposed method respectively}.}
\label{imagestoctyforaver7}\end{figure}

Here, we also show the images restored by the three methods. We display the degraded images and the restored images by three methods under two Gaussian blurs (i) and (ii) and 50\% level of noise. The results are show in Fig.~\ref{camtoctygln}. We can easily see the third advantage of our proposed method that our method can overcome the staircase effects effectively and get better visual quality than others. Moreover, we also plot the evolution of the PSNR over time and iterations for the three different methods in Fig.~\ref{gau75noi4PSNR2Itrtime} for the image blurred by $7\times 7$ Gaussian blur and corrupted by 40\% level salt-and-pepper noise.

From the expriments and the description in \cite{GLNG2009}, the GLN method has three sensitive parameters that depend on blur, noise level and test images rather than only one sensitive regular parameter $\mu$ in CTY and our proposed method. Besides, the results of GLN is nearly same as FTVd, while CTY is much better than FTVd in \cite{CTY2013}. From the above tests, we also see that both our proposed method and CTY can get better results than GLN. Moerever, the staircase effects by GLN method is also existent. Therefore, we omit the following comparison with GLN and only list the comparison with CTY.

{\bf {\em Remarks}}. Here and in the following tests for CTY, we tune the regularization parameter $\mu$ to be optimal by checking the highest PSNR and the lowest ReE under a ``for loop'' of Matlab code from 1 to 70 by step length 1 for different images under different blurs and noise levels. Besides, in the experiments, we find that, under the $15\times 15$ Gaussion blur with standard deviation 5, we can increase inner penalty parameters $\beta_2$ to get higher PSNR. We do not change the parameter in our proposed method because of the good properties of a good restoration algorithm introduced above. However, here in Table~\ref{cameranmantol1CA} for the CTY method under $15\times 15$ Gaussion blur, we set the inner penalty parameter $\beta_2=50$ instead of $\beta_2=20$ in \cite{CTY2013} for higher PSNR. For other tests following we also set this inner penalty parameter of CTY to be as default ($\beta_2=20$) in \cite{CTY2013}.

\subsection{Comparison with CTY for other test images}
In this subsection, we focus on comparisons between our proposed method and CTY for deblurring problems under salt-and-pepper noise.
Fifty-six degraded test images are generated in the way similar to that in Section 4.2. That is, we first generated the blurred images operating on images (b)-(h) with the periodic boundary condition by two blurs Gaussian blur (i) (also as G) and average blur (iii) (also as A), then corrupted the blurred images by salt-and-pepper noise from 30\% to 60\%. The parameters of our proposed method are set as above description in Section 4.1 and the parameters of CTY as remarks in Section 4.2.

Conclusions similar to those in Section 4.2 can be made based on the results in Table~\ref{allimagecty2ourforgauss} and Table~\ref{allimagecty2ourforaverage}.
For example, our proposed method is always more accurate, with a possible improvement of more than 2.90 dB in PSNR (see image (c) with Gaussian blur and a 30\% level of noise). For all images, the lower the noise level is, the better the improvement of PSNR will be. Even in high noise level, our method is more accurate than CTY by sacrificing partial time. For a further step, the iterations of our method are fewer than CTY for almost all test.


Finally, in Fig.~\ref{imagestoctyforgauss7} and Fig.~\ref{imagestoctyforaver7}, we display the zoom parts of the degraded images examples and the zoom parts of the restored text images respectively for Gaussian blur and average blur with noise level from 30\% to 60\% by two methods. We can easily see the visual improvement in the images by using our method. More specifically, from the third column in Fig.~\ref{imagestoctyforaver7}, although the numerical result does not improve too much, the image visual quality of our proposed method is much better than CTY since our method can overcome the staircase effects effectively. This superiority is obvious for almost all the test images in our work.

\section{Conclusion}

In this paper, we study a
new regularization term using TV with OGS in the classic $\ell_1$-TV model
for the image deblurring with impulse noise. We propose an efficient solving algorithm under the framework of the general ADMM.
 We use an MM inner iteration to solve the subproblem instead of Shrinkage \cite{YWYW2009} in the classic $\ell_1$-TV model. Based on them, we propose an algorithm called CL1-OGS-ATV-ADM4.
 The numerical results illustrate that our method outperforms CTY \cite{CTY2013} and GLN \cite{GLNG2009} both in numerical results and image visual quality. Particularly, our proposed method can overcome staircase effects effectively while CTY and GLN can not.
\section*{Acknowledgment}
The authors would like to thank Prof. M. Tao for providing us the code ADM2CTVL1 (CTY) in \cite{CTY2013}.

\end{document}